\documentclass[3p]{elsarticle}

\usepackage{amssymb,amsmath}
 \newtheorem{thm}{Theorem}[section]
 \newtheorem{cor}[thm]{Corollary}
 \newtheorem{lem}[thm]{Lemma}
 \newtheorem{prop}[thm]{Proposition}
 \newdefinition{defn}[thm]{Definition}
 \newdefinition{rem}[thm]{Remark}
 \newdefinition{ex}[thm]{Example}
\numberwithin{equation}{section}

\newproof{pf}{Proof}

\newcommand{\N}{{\mathbb N}}
\newcommand{\C}{{\mathbb C}}

\newcommand{\D}{{\mathbb D}}
\newcommand{\E}{{\mathbb E}}
\newcommand{\T}{{\mathbb T}}
\newcommand{\K}{{\mathbb K}}
\newcommand{\cL}{{\mathcal L}}
\newcommand{\cR}{{\mathcal R}}
\newcommand{\cQ}{{\mathcal Q}}
\newcommand{\cS}{{\mathcal S}}
\newcommand{\cT}{{\mathcal T}}
\newcommand{\cA}{{\mathcal A}}
\newcommand{\cB}{{\mathcal B}}
\newcommand{\cP}{{\mathcal P}}
\newcommand{\cC}{{\mathcal C}}
\newcommand{\cD}{{\mathcal D}}
\newcommand{\cE}{{\mathcal E}}
\newcommand{\cM}{{\mathcal M}}
\newcommand{\cN}{{\mathcal N}}
\newcommand{\cNM}{{\mathcal N\mathcal M}}
\newcommand{\cY}{{\mathcal Y}}
\newcommand{\cZ}{{\mathcal Z}}
\newcommand{\gC}{{\mathfrak C}}
\newcommand{\gD}{{\mathfrak D}}

\newcommand{\gK}{{\mathfrak K}}

\newcommand{\fC}{{\mathbf C}}

\newcommand{\fP}{{\mathbf P}}
\newcommand{\fQ}{{\mathbf Q}}

\newcommand{\fpi}{\mbox{\boldmath$\pi$\unboldmath}}
\newcommand{\frho}{\mbox{\boldmath$\rho$\unboldmath}}
\newcommand{\fsigma}{\mbox{\boldmath$\sigma$\unboldmath}}
\newcommand{\ftau}{\mbox{\boldmath$\tau$\unboldmath}}

\newcommand{\eps}{{\varepsilon}}

\newcommand{\rank}{{\rm rank\,}}
\newcommand{\diag}{{\rm diag}}
\newcommand{\sq}{\subseteq}
\newcommand{\ul}{\underline}

\begin{document}

\begin{frontmatter}
\title{On an Interpolation Problem for $J$-Potapov Functions}

\author[]{Bernd Fritzsche}\ead{fritzsche@math.uni-leipzig.de}
\author[]{Bernd Kirstein}\ead{kirstein@math.uni-leipzig.de}
\author[]{Uwe Raabe\fnref{RD}}\ead{raabe@math.uni-leipzig.de}

\address{Mathematisches Institut,
Universit\"at Leipzig,
Augustusplatz~10/11,
04109~Leipzig,
Germany}

%


\fntext[RD]{The work of the third author of the present paper was supported by the
EU project "Geometric Analysis on Lie Groups and Applications" (GALA).}


\begin{abstract}
Let $J$ be an $m\times m$ signature matrix (i.e. $J^*=J$ and $J^2=I_m$) and let $\D:=\{z\in\C:|z|<1\}$. Denote by $\cP_{J,0}(\D)$ the class of all meromorphic $m\times m$ matrix-valued functions $f$ in $\D$ which are holomorphic at $0$ and take $J$-contractive values at all points of $\D$ at which $f$ is holomorphic. 
The central theme of this paper is the study of the following interpolation problem:

Let $n$ be a nonnegative integer, and let $(A_j)_{j=0}^n$ be a sequence of complex $m\times m$ matrices. Describe the set $\cP_{J,0}[\D,(A_j)_{j=0}^n]$ of all matrix-valued functions $f\in\cP_{J,0}(\D)$ such that
$\frac{f^{(j)}(0)}{j!}=A_j$
for each $j\in\N_{0,n}$ where the notation $f^{(j)}$ stands for the $j$-th derivative of $f$. 
In particular, characterize the case that the set $\cP_{J,0}[\D,(A_j)_{j=0}^n]$ is nonempty.

In this paper, we will solve this problem in the most general case. Moreover, in the nondegenerate case we will give a description of the corresponding Weyl matrix balls. Furthermore, we will investigate the limit behaviour of the Weyl matrix balls associated with the functions belonging to some particular subclass of  $\cP_{J,0}(\D)$.

\end{abstract}
\end{frontmatter}


\setcounter{section}{-1}
\section{Introduction}
\label{S0}

This paper deals with an interpolation problem for a particular class of meromorphic matrix-valued
functions in the unit disk $\D:=\{z\in\C:|z|<1\}$. This class is the Potapov class $\cP_J(\D)$ (see section \ref{S1}). 
It originates in the fundamental paper \cite{Po1} by V.P. Potapov whose investigations where initiated by the studies of M.S. Liv\v sic (\cite{L1}, \cite{L2}) on characteristic functions of nonunitary operators. 

Our interest is concentrated on the subclass $\cP_{J,0}(\D)$ of all functions belonging to 
$\cP_J(\D)$ which are holomorphic at $0$. Given a finite sequence $(A_j)_{j=0}^n$ from 
$\C^{m\times m}$, we want to determine all $m\times m$ matrix-valued functions $f$ belonging to
$\cP_{J,0}(\D)$ such that their Taylor coefficients sequence starts with the section $(A_j)_{j=0}^n$.

Interpolation problems have a rich history. Important results for the scalar case were already obtained in the first half of the 20th century. In the early 1950's a new period started, where interpolation problems for matrix-valued functions were considered. These investigations culminated in a series of monographs (see, e.g., \cite{BGR}, \cite{DFK}--\cite{FFGK}, \cite{Sak}).
An essential common feature of these monographs is that the considerations mainly concentrated on the 
so-called nondegenerate case which is connected to positive Hermitian block Pick matrices built 
from the given data in the interpolation problem. 

The study of the degenerate case (where the associated block Pick matrix is nonnegative Hermitian and singular) began with the pioneering work \cite{Dub1}--\cite{Dub6} of V.K. Dubovoj in the framework 
of the matricial Schur problem. In the sequel, quite different approaches to handle degenerate cases 
of matrix interpolation were used (see, e.g., \cite{BH}--\cite{DGK}, \cite[Chapter 7]{Dy}, and \cite[Chapter 5]{Sak}). 

This paper makes extensive use of the 
authors' recent investigations \cite{FKR1}--\cite{FKRS} on $J$-Potapov functions,
$J$-Potapov sequences, and their interrelations. For the basic strategy for treating the interpolation problem we draw upon the technique which was introduced in \cite{FKL1} and \cite{FKL}
to handle simultaneously the degenerate and nondegenerate cases in the matrix versions of the the 
interpolation problems named after Carath\'eodory and Schur (see \cite{Sch}), respectively. 
Our method is essentially based on our investigations \cite{FKRS} on the $J$-central $J$-Potapov
functions associated with a finite $J$-Potapov sequence of complex matrices. In particular, we will
make frequently use of the matrix ball description of the elements of a $J$-Potapov sequence.
The main results of this paper (Theorems \ref{S1-T1}, \ref{DSS-T1}, \ref{DSS-T2}, and \ref{nd-T1})
contain descriptions of the solution set of the interpolation problem under consideration in terms
of a linear fractional transformation the generating matrix-valued function of which is a 
matrix polynomial. The canonical blocks of this matrix polynomial will be constructed with the
aid of those matrix polynomials which were used in \cite{FKRS} to derive right and left quotient representations of $J$-central $J$-Potapov functions (see Theorem \ref{LRQ-T1}).

This paper is organized as follows. In Section 1, we summarize some basic facts on $J$-Potapov functions in the open unit disk, $J$-Potapov sequences and their interrelations. These results originate from \cite{FKR1} and \cite{FKR2}.

In Section 2, we study $J$-central $J$-Potapov functions. We verify that an arbitrary function $f$ 
belonging to the class $\cP_{J,0}(\D)$ can be approximated by the sequence of associated 
$J$-central $J$-Potapov functions (see Corollary \ref{CPF-C0}). Moreover, we recall the left and right quotient representations of $J$-central $J$-Potapov functions which were obtained in
\cite{FKRS} (see Theorem \ref{LRQ-T1}). These representations play the key role in the sequel.

Section 3 is devoted to the description of the solution set $\cP_{J,0}[\D,(A_j)_{j=0}^n]$ of the interpolation problem studied in this paper. This problem will be solved in the most general case
(see Theorems \ref{DSS-T1}, \ref{DSS-T2}, and \ref{DSS-T3}). 

In Section 4, we state several characterizations of the unique solvability of the interpolation 
problem under consideration.

In Section 5, we treat the so-called nondegenerate case of the interpolation problem, 
namely the situation where some strict $J$-Potapov sequence $(A_j)_{j=0}^n$ is given. In this case the representation of the solution set $\cP_{J,0}[\D,(A_j)_{j=0}^n]$ contained in Theorem \ref{nd-T1}
generalizes the result which was obtained by Arov/Krein \cite{AK} (see also \cite[Section 3.10]{DFK}
and \cite{FK2a}). Moreover, we will compute the corresponding Weyl matrix balls (see Theorem \ref{nd-P1}).

The central theme of Section 6 is the investigation of the interrelations between the Weyl matrix balls which are associated with a strict $J$-Potapov sequence  $(A_j)_{j=0}^n$ on the one hand and with its $J$-Potapov-Ginzburg transform  $(B_j)_{j=0}^n$ on the other hand. The main result of Section 6 is Proposition \ref{nd-P1B} which contains explicite formulas which express several interrelations between the corresponding parameters of the two
Weyl matrix balls under consideration.

The final Section 7 deals with the study of the limit behaviour of the sequence of the Weyl matrix balls
 associated with a nondegenerate $J$-Potapov function $f$. Using the fact that the $J$-Potapov-Ginzburg transform $g$ of $f$ turns out to be a nondegenerate $m\times m$ Schur function and taking into account the formulas obtained in Section 6 the desired result on the ranks of the semi-radii of the Weyl matrix balls under consideration (see Proposition \ref{ndlim-P5}) will be derived from the corresponding result for nondegenerate $m\times m$ Schur functions (see \cite[Theorem 3.11.2]{DFK}).

\section{Some preliminaries on $J$-Potapov functions in the open unit disk and on $J$-Potapov sequences}
\label{S1}

Throughout this paper, let $m$ be a positive integer. 
We will use the notations $\N$, $\N_0$, and $\C$ for the set of all positive integers, the set of all nonnegative integers, and the set of all complex numbers, respectively. 
If $s\in\N_0$ and $\kappa\in\N_0\cup\{+\infty\}$ then $\N_{s,\kappa}$ denotes the set of all integers $n$ satisfying $s\le n\le\kappa$.
Further, let $\D:=\{z\in\C:|z|<1\}$ and $\T:=\{z\in\C:|z|=1\}$.

Let $p,q\in\N$. Then $\C^{p\times q}$ designates the set af all complex $p\times q$ matrices. 
The notation $0_{p\times q}$ stands for the null matrix which belongs to
$\C^{p\times q}$, and the identity matrix which belongs to $\C^{q\times q}$ will be designated by $I_q$.
In cases where the size of a null matrix or the size of an identity matrix is obvious, we will omit the indices.
If $A\in\C^{p\times q}$ then $A^+$ stands for the Moore-Penrose inverse of $A$.
Furthermore, for each $A\in\C^{p\times q}$, let $\cR(A)$ be the range of $A$, let $\cN(A)$ be the nullspace of $A$, 
and let $\|A\|$ denote the spectral norm of $A$. 
We will write $\C^{q\times q}_\ge$ 
(respectivly, $\C^{q\times q}_>$) to denote the set af all nonnegative (respectively, positive) Hermitian matrices belonging to $\C^{q\times q}$.
In the set $\C^{q\times q}_H$ of all Hermitian $q\times q$ matrices we will use the L\"owner semi-ordering, i.e., we will write $A\le B$ or $B\ge A$ to indicate that $A$ and $B$ are Hermitian matrices of the same size such that $B-A$ is nonnegative Hermitian.
Moreover, we will write $A<B$ or $B>A$ to indicate that $A$ and $B$ are Hermitian matrices of the same size such that $B-A$ is positive Hermitian.

Let $n$ and $p_1,\ldots,p_n$ be positive integers, and let $A_j\in\C^{p_j\times p_j}$ for each 
$j\in\N_{1,n}$. Then $\diag(A_1,\ldots,A_n)$ denotes the block diagonal matrix with diagonal blocks
$A_1,\ldots,A_n$.

If $f$ is an $m\times m$ matrix-valued function which is meromorphic in the open unit disk $\D$, then let $\mathbb H_f$ be the set of all points at which $f$ is holomorphic. Let $J$ be an $m\times m$ signature matrix and let $f$ be a $\C^{m\times m}$-valued function which is meromorphic in $\D$. Then $f$ is called a \emph{$J$-Potapov function in $\D$} (respectively, a \emph{strong $J$-Potapov function in $\D$}), if for each $w\in\mathbb H_f$ the matrix $f(w)$ is $J$-contractive (respectively, strictly $J$-contractive). Here a matrix $A\in\C^{m\times m}$ is called \emph{$J$-contractive} (respectively, \emph{strictly $J$-contractive}), if the matrix 
${J-A^*JA}$ is nonnegative Hermitian (respectively, positive Hermitian). For each $m\times m$ signature matrix $J$, we will use the notation $\cP_J(\D)$ (respectively, $\cP_J'(\D)$) to denote the set of all $J$-Potapov functions in $\D$ (respectively, strong $J$-Potapov functions in $\D$). We will turn particular attention to a distinguished subclass of $\cP_J(\D)$, namely the class
$$\cP_{J,0}(\D):=\{f\in\cP_J(\D):\ 0\in\mathbb H_f\}.$$
In the case $J=I_m$ the classes $\cP_J(\D)$ and $\cP_{J,0}(\D)$ coincide. Indeed, $\cP_{I_m}(\D)$ is exactly the set $\cS_{m\times m}(\D)$ of all $m\times m$ Schur functions in $\D$, i.e., the set of all matrix-valued functions $f:\D\rightarrow\C^{m\times m}$ which are holomorphic in $\D$ and the values of which are contractive complex matrices.\\
Observe that the well-known concept of Potapov-Ginzburg transformation yields an interrelation between the classes $\cP_J(\D)$ and $\cS_{m\times m}(\D)$ on the one-hand side and between the strong $J$-Potapov class $\cP_J'(\D)$ and the strong Schur class $\cS_{m\times m}'(\D)$ 
of all $f\in\cS_{m\times m}(\D)$ for which the matrix $f(w)$ is strictly contractive for each $w\in\D$ on the other-hand side (see \cite[Proposition 3.4]{FKR1}).

The sequences $(A_j)_{j=0}^\infty$ of Taylor coefficients of the matrix-valued functions which belong to the class $\cP_{J,0}(\D)$ can be characterized in a clear way. In order to recall this characterization we introduce some notations.
Observe that, for each $m\times m$ signature matrix $J$ and every nonnegative integer $n$, the complex $(n+1)m\times(n+1)m$ matrix 
\begin{equation}\label{JnDef}
J_{[n]}:=\diag(J,\ldots,J)
\end{equation}
is an $(n+1)m\times(n+1)m$ signature matrix. If $n\in\N_0$, then a sequence $(A_j)_{j=0}^n$ of complex $m\times m$ matrices is called a \emph{$J$-Potapov sequence} (respectively, a \emph{strict $J$-Potapov sequence}) if the block Toeplitz matrix
\begin{equation}\label{NrS1}
S_n :=\begin{pmatrix}
 A_0 & 0_{m\times m} & \ldots & 0_{m\times m} \\
          A_1 & A_0 & \ldots & 0_{m\times m} \\
          \vdots & \vdots & & \vdots \\
          A_n & A_{n-1} & \ldots & A_0
\end{pmatrix}.
\end{equation}
is $J_{[n]}$-contractive (respectively, strictly $J_{[n]}$-contractive). 
If necessary, we will write $S_n^{(A)}$ instead of $S_n$ to indicate the sequence $(A_j)_{j=0}^n$
from which the matrix is built.

For each $n\in\N_0$ we will use $\mathcal P_{J,n}^{\leq}$ (respectively, $\mathcal P_{J,n}^{<}$) to designate the set of all $J$-Potapov sequences (respectively, strict $J$-Potapov sequences) $(A_j)_{j=0}^n$. 
From \cite[Lemma 3.2 (respectively, Lemma 3.3)]{FKR2} it follows 
that if $(A_j)_{j=0}^n$ belongs to $\mathcal P_{J,n}^{\leq}$ (respectively, $\mathcal P_{J,n}^<$), then $(A_j)_{j=0}^k\in\mathcal P_{J,k}^{\leq}$ (respectively, $(A_j)_{j=0}^k\in\mathcal P_{J,k}^<$) for each $k\in\N_{0,n}$. A sequence $(A_j)_{j=0}^\infty$ of complex $m\times m$ matrices is said to be a \emph{$J$-Potapov sequence} (respectively, a \emph{strict $J$-Potapov sequence}) if for each $n\in\N_0$ the sequence $(A_j)_{j=0}^n$ is a $J$-Potapov sequence (respectively, a strict $J$-Potapov sequence). We will write $\mathcal P_{J,\infty}^{\leq}$ for the set of all $J$-Potapov sequences $(A_j)_{j=0}^\infty$ and $\mathcal P_{J,\infty}^{<}$ for the set of all strict $J$-Potapov sequences $(A_j)_{j=0}^\infty$.

Observe that the concept of $J$-Potapov sequences is a generalization of the well-known concept
of $m\times m$ Schur sequences. Here, for each $\kappa\in\N_0\cup\{+\infty\}$, a sequence
$(A_j)_{j=0}^\kappa$ of complex $m\times m$ matrices is said to be an \emph{$m\times m$ Schur sequence} (respectively, a \emph{strict $m\times m$ Schur sequence}) if for each $n\in\N_{0,\kappa}$ the 
matrix $S_n$ given by (\ref{NrS1}) is contractive (respectively, strictly contractive).

Now we can formulate the Taylor series characterization of the class $\cP_{J,0}(\D)$.
\begin{thm}\label{T62}
Let $J$ be an $m\times m$ signature matrix. Then:
\begin{enumerate}
 \item[(a)] If $f\in\cP_{J,0}(\D)$ and if
\begin{equation}\label{NrTS}
f(w)=\sum_{j=0}^\infty A_jw^j
\end{equation}
is the Taylor series representation of $f$ in some neighborhood of $0$, then $(A_j)_{j=0}^\infty$ is a $J$-Potapov sequence.
\item[(b)] If $(A_j)_{j=0}^\infty$ is a $J$-Potapov sequence, then there is a unique $f\in\cP_{J,0}(\D)$ such that (\ref{NrTS}) holds for all $w$ belonging to some neighborhood of $0$.
\end{enumerate}
\end{thm}
A proof of Theorem \ref{T62} is given in \cite[Theorem 6.2]{FKR1}.

Considering the special case $J=I_m$ one can see immediately that Theorem \ref{T62} is a generalization of a well-known characterization of the Taylor coefficients of matricial Schur functions defined on $\D$ (see, e.g., \cite[Theorem 5.1.1]{DFK}).

The main goal of this paper is to give a description of the solution set of the following interpolation problem for functions belonging to the class $\cP_{J,0}(\D)$:

\vspace{2ex}
\textbf{Interpolation problem for Potapov functions (P):} Let $J$ be an $m\times m$ signature matrix, let $n\in\N_0$, and let $(A_j)_{j=0}^n$ be a sequence of complex $m\times m$ matrices. Describe the set $\cP_{J,0}\left[\D,(A_j)_{j=0}^n\right]$ of all matrix-valued functions $f\in\cP_{J,0}(\D)$ such that
\begin{equation}\label{ICP}
\frac{f^{(j)}(0)}{j!}=A_j
\end{equation} 
for each $j\in\N_{0,n}$ where the notation $f^{(j)}$ stands for the $j$-th derivative of $f$. 
In particular, characterize the case that the set $\cP_{J,0}\left[\D,(A_j)_{j=0}^n\right]$ is nonempty.

\vspace{2ex}
The following theorem characterizes the situation that the problem (P) has a solution.
\begin{thm}\label{T72}
Let $J$ be an $m\times m$ signature matrix, let $n\in\N_0$, and let $(A_j)_{j=0}^n$ be a sequence of complex $m\times m$ matrices. Then the set $\cP_{J,0}\left[\D,(A_j)_{j=0}^n\right]$ is nonempty if and only if $(A_j)_{j=0}^n$ is a $J$-Potapov sequence.
\end{thm}
A proof of Theorem \ref{T72} is given in \cite[Theorem 7.2]{FKR1}.


We will now give some more notations that will be used throughout this paper.
For each $n\in\N_0$, let the matrix polynomials $e_{n,m}:\C\rightarrow\C^{m\times(n+1)m}$ 
and $\varepsilon_{n,m}:\C\rightarrow\C^{(n+1)m\times m}$ be defined by
\begin{equation}\label{Nrenm}
e_{n,m}(w):=(I_m,wI_m,\ldots,w^nI_m) \quad\mbox{and}\quad
\eps_{n,m}(w):=(\overline w^{\,n}I_m,\overline w^{\,n-1}I_m,\ldots,I_m)^*.
\end{equation}

In this paper, we will frequently use the notion of the reciprocal matrix polynomial.
Let $p,q\in\N$, and let $b$ be a $p\times q$ matrix
polynomial, i.e., there are an $n\in\N_0$ and a matrix $B\in\C^{(n+1)p\times q}$ such that 
$b(w)=e_{n,p}(w)B$ holds for each $w\in\C$. Then the \emph{reciprocal matrix polynomial
$\tilde b^{[n]}$ of $b$ with respect to the unit circle $\T$ and the formal degree $n$} is given
by $\tilde b^{[n]}(w):=B^*\eps_{n,p}(w)$ for each $w\in\C$.
If $\beta$ is the restriction of $b$ onto $\D$ then let
$\tilde\beta^{[n]}$ be the restriction of $\tilde b^{[n]}$ onto $\D$.

Let $J$ be an $m\times m$ signature matrix, and let $\kappa\in\N_0\cup\{+\infty\}$. 
Whenever a sequence $(A_j)_{j=0}^\kappa$ of complex $m\times m$ matrices
is given, then the following notations will be used throughout this paper.
For each $n\in\N_{0,\kappa}$, let $S_n$ be given by (\ref{NrS1}), and let
\begin{equation}\label{NrPQ}
P_{n,J}:=J_{[n]}-S_nJ_{[n]}S_n^*\qquad\mbox{ and }\qquad Q_{n,J}:=J_{[n]}-S_n^*J_{[n]}S_n.
\end{equation}
In the case $n\in\N_{1,\kappa}$ we will use the block matrices
\begin{equation}\label{Nrynzn}
y_n:=(A_1^*,A_2^*,\ldots,A_n^*)^* \qquad\mbox{ and }\qquad z_n:=(A_n,A_{n-1},\ldots,A_1).
\end{equation}
If necessary, we will write $y_n^{(A)}$ (respectively, $z_n^{(A)}$) instead of $y_n$ (respectively, $z_n$) to indicate the sequence $(A_j)_{j=0}^\kappa$ from which the matrix is built.
Moreover, for each $n\in\N_{0,\kappa}$, we will work with the matrices
\begin{equation}\label{NrMkDef}
M_{n+1,J}:=\left\{\begin{array}{cl} 0_{m\times m}, &\mbox{if}\ n=0 \\
                             -z_nJ_{[n-1]}S_{n-1}^*P_{n-1,J}^+y_n, &\mbox{if}\ n\in\N_{1,\kappa} ,
                \end{array}\right.
\end{equation}
\begin{equation}\label{NrLkDef}
L_{n+1,J}:=\left\{\begin{array}{cl} J-A_0JA_0^*, &\mbox{if}\ n=0 \\
                             J-A_0JA_0^*-z_nQ_{n-1,J}^+z_n^*, &\mbox{if}\ n\in\N_{1,\kappa} 
                \end{array}\right.
\end{equation}
and
\begin{equation}\label{NrRkDef}
R_{n+1,J}:=\left\{\begin{array}{cl} J-A_0^*JA_0, &\mbox{if}\ n=0 \\
                             J-A_0^*JA_0-y_n^*P_{n-1,J}^+y_n, &\mbox{if}\ n\in\N_{1,\kappa} 
                \end{array}\right.
\end{equation}

Observe that if $(A_j)_{j=0}^\kappa$ is a $J$-Potapov sequence, then for each $n\in\N_{0,\kappa}$ the matrices $L_{n+1,J}$ and $R_{n+1,J}$ are both nonnegative Hermitian (see \cite[Lemma 3.7]{FKR2}).

Let us now recall the notion of a matrix ball. Let $p,q\in\N$, and
denote by $\K_{p\times q}$ the set of all contractive 
matrices from $\C^{p\times q}$. Then,
for each $M\in\C^{p\times q}$, each $L\in\C^{p\times p}$ and each $R\in\C^{q\times q}$, the set 
$$ \gK(M;L,R):=\{M+LKR:K\in\K_{p\times q}\} $$
is called the \emph{matrix ball} with \emph{center} $M$, \emph{left semi-radius} $L$, and \emph{right semi-radius} $R$.
The theory of matrix and operator balls was developed by Yu.L. \v Smuljan \cite{Sm} 
(see also \cite[Section 1.5]{DFK}).

In \cite[Theorem 3.9]{FKR2} there has been shown the following result which enlightens the inner structure of $J$-Potapov sequences.

\begin{thm}\label{T}
Let $J$ be an $m\times m$ signature matrix, let $n\in\N_0$, and let $(A_j)_{j=0}^{n+1}$ be a sequence of complex $m\times m$ matrices. 
Then the following statements are equivalent:
\begin{enumerate}
 \item[(i)] $(A_j)_{j=0}^{n+1}$ is a $J$-Potapov sequence.
 \item[(ii)] $(A_j)_{j=0}^{n}$ is a $J$-Potapov sequence and $A_{n+1}$ belongs to the matrix ball
$$\gK\big(M_{n+1,J};\sqrt{L_{n+1,J}},\sqrt{R_{n+1,J}}\,\big).$$
\end{enumerate}
\end{thm}

Considering the special choice $J=I_m$ we see that Theorem \ref{T} is a generalization of a well-known result for $m\times m$ Schur sequences 
(see, e.g., \cite[Theorem 1]{FK2} or \cite[Theorem 3.5.1]{DFK}).

The main goal of the present paper is to describe the set of solutions of problem (P). In particular, we will
prove the following theorem, which is a generalization of \cite[Theorem 1.1]{FKL}.

\begin{thm}\label{S1-T1}
Let $J$ be an $m\times m$ signature matrix, let $n\in\N_0$, 
and let $(A_j)_{j=0}^n$ be a $J$-Potapov sequence. 
Let the $m\times m$ matrix polynomials 
$\pi_{n,J}$, $\rho_{n,J}$,
$\sigma_{n,J}$, and $\tau_{n,J}$ 
be given by
$$\pi_{n,J}(w):=\left\{\begin{array}{cl} A_0, &\mbox{if }n=0\\
               A_0+w e_{n-1,m}(w)(I_m+S_{n-1}Q_{n-1,J}^+S_{n-1}^*J_{[n-1]})y_n, 
                        &\mbox{if }n\in\N,\end{array}\right.$$
$$\rho_{n,J}(w):=\left\{\begin{array}{cl} I_m, &\mbox{if }n=0\\
                             I_m+w e_{n-1,m}(w)Q_{n-1,J}^+S_{n-1}^*J_{[n-1]}y_n, 
                               &\mbox{if }n\in\N,\end{array}\right.$$
$$\sigma_{n,J}(w):=\left\{\begin{array}{cl} A_0, &\mbox{if }n=0\\
          z_n(J_{[n-1]}S_{n-1}^*P_{n-1,J}^+S_{n-1}+I_m)w\eps_{n-1,m}(w)+A_0,
                  &\mbox{if }n\in\N, \end{array}\right.$$
and
$$\tau_{n,J}(w):=\left\{\begin{array}{cl} I_m, &\mbox{if }n=0\\
                     z_nJ_{[n-1]}S_{n-1}^*P_{n-1,J}^+w\eps_{n-1,m}(w)+I_m,
                       &\mbox{if }n\in\N \end{array}\right.$$
for each $w\in\C$. For every $S\in\cS_{m\times m}(\D)$ and each $w\in\D$ for which the matrix \\ 
$wJ\tilde\sigma_{n,J}^{[n]}(w)\sqrt{L_{n+1,J}}^+S(w)\sqrt{R_{n+1,J}}+\rho_{n,J}(w)$ is nonsingular,
let
\begin{align}\label{S1-T1-1}
f_S(w)&:=\big(wJ\tilde\tau_{n,J}^{[n]}(w)\sqrt{L_{n+1,J}}^+S(w)\sqrt{R_{n+1,J}}+\pi_{n,J}(w)\big) \nonumber\\
&\qquad\cdot    \big(wJ\tilde\sigma_{n,J}^{[n]}(w)\sqrt{L_{n+1,J}}^+S(w)\sqrt{R_{n+1,J}}+\rho_{n,J}(w)\big)^{-1}. 
\end{align}
Then, for each $S\in\cS_{m\times m}(\D)$, by (\ref{S1-T1-1}) a matrix-valued function $f_S$  meromorphic in $\D$ is given, and 
the set $\mathbb H_{f_S}$ of all $w\in\D$ at which $f_S$ is holomorphic fulfills
\begin{align*}
\mathbb H_{f_S}&=\big\{w\in\D:\det
\big(wJ\tilde\sigma_{n,J}^{[n]}(w)\sqrt{L_{n+1,J}}^+S(w)\sqrt{R_{n+1,J}}+\rho_{n,J}(w)\big)\ne0\big\}\\
&=\big\{w\in\D:\det
\big(w\sqrt{L_{n+1,J}}S(w)\sqrt{R_{n+1,J}}^+\tilde\pi_{n,J}^{[n]}(w)J+\tau_{n,J}(w)\big)\ne0\big\}.
\end{align*}
Further, for each $S\in\cS_{m\times m}(\D)$ and each $w\in\mathbb H_{f_S}$,
$f_S$ admits the representation
\begin{align*}
f_S(w)&=\big(w\sqrt{L_{n+1,J}}S(w)\sqrt{R_{n+1,J}}^+\tilde\pi_{n,J}^{[n]}(w)J+\tau_{n,J}(w)\big)^{-1} \\
&\qquad\cdot
   \big(w\sqrt{L_{n+1,J}}S(w)\sqrt{R_{n+1,J}}^+\tilde\rho_{n,J}^{[n]}(w)J+\sigma_{n,J}(w)\big).
\end{align*}
Moreover, 
$$ \cP_{J,0}[\D,(A_j)_{j=0}^n]=\{f_S:S\in\cS_{m\times m}(\D)\} $$
holds true.
\end{thm}

\section{$J$-central $J$-Potapov functions}

A crucial idea in our approach to the interpolation problem (P) consists in comparing possible candidates
for solutions with a distinguished solution, namely with the so-called $J$-central $J$-Potapov function corresponding to the given $J$-Potapov sequence $(A_j)_{j=0}^n$.
Let $n\in\N_0$, and let $(A_j)_{j=0}^n$ be a $J$-Potapov sequence. Then Theorem \ref{T} implies that
the sequence
$(A_j)_{j=0}^\infty$ defined recursively by $A_k:=M_{k,J}$ for each $k\in\N_{n+1,\infty}$ is
a $J$-Potapov sequence. The sequence $(A_j)_{j=0}^\infty$ is said to be 
\emph{the $J$-central $J$-Potapov sequence corresponding to $(A_j)_{j=0}^n$}.
In view of Theorem \ref{T62},
there is a unique $f_{c,n}\in\cP_{J,0}(\D)$ such that $f_{c,n}(w)=\sum_{j=0}^\infty A_jw^j$ holds for all $w$ belonging to some neighborhood of $0$.
The matrix-valued function $f_{c,n}$ is called \emph{the $J$-central $J$-Potapov function
corresponding to $(A_j)_{j=0}^n$}.
A more detailed study of $J$-central $J$-Potapov functions can be found in \cite{FKRS}.
The following result complements Theorem \ref{T72}.

\begin{prop}\label{CPF-P0}
Let $J$ be an $m\times m$ signature matrix, let $n\in\N_0$, and let $(A_j)_{j=0}^n$ be a $J$-Potapov
sequence. Then the $J$-central $J$-Potapov function $f_{c,n}$ corresponding to $(A_j)_{j=0}^n$
belongs to $\cP_{J,0}[\D,(A_j)_{j=0}^n]$.
\end{prop}

\begin{pf}
The assertion is an immediate consequence of the construction of $f_{c,n}$.
\qed\end{pf}

Observe that the concept of $J$-central $J$-Potapov sequences and $J$-central $J$-Potapov functions is a generalization of the well-known concept of central $m\times m$ Schur sequences
and central $m\times m$ Schur functions in $\D$.

\begin{prop}\label{CPF-P2}
Let $J$ be an $m\times m$ signature matrix, and let $f\in\cP_{J,0}(\D)$. Let (\ref{NrTS})
be the Taylor series representation of $f$ in some neighborhood of $0$. For each $n\in\N_0$,
let $f_n\in\cP_{J,0}[\D,(A_j)_{j=0}^n]$.
Then, for each compact subset $K$ of $\mathbb H_f$, there exists a nonnegative integer $n_0$ such that $f_n$ is holomorphic in $K$ for every integer $n$ with $n\ge n_0$. Moreover,
$$ \lim_{n\to\infty}f_n(w)=f(w) $$
holds for each $w\in\mathbb H_f$. This convergence is uniform in each compact subset of $\mathbb H_f$.
\end{prop}

\begin{pf}
Let $\mathbf P_J:=\frac12(I+J)$ and $\mathbf Q_J:=\frac12(I-J)$, and let $n\in\N_0$.
According to \cite[Proposition 3.4]{FKR1}, we have $\det(\mathbf Q_Jf(0)+\mathbf P_J)\ne0$ and $\det(\mathbf Q_Jf_n(0)+\mathbf P_J)\ne0$.
Let
$$g:=(\mathbf P_Jf+\mathbf Q_J)(\mathbf Q_Jf+\mathbf P_J)^{-1}\quad\mbox{and}\quad
g_n:=(\mathbf P_Jf_n+\mathbf Q_J)(\mathbf Q_Jf_n+\mathbf P_J)^{-1}.$$
Then an application of \cite[Proposition 3.4]{FKR1} yields that $g$ and $g_n$ both belong to 
$\cS_{m\times m}(\D)$.
Let $g(w)=\sum_{j=0}^\infty B_jw^j$ 
(respectively, $g_n(w)=\sum_{j=0}^\infty B_j^{(n)}w^j$)
be the Taylor series representation of $g$ (respectively, of $g_n$) in $\D$. 
Then $(B_j)_{j=0}^\infty$ and $(B_j^{(n)})_{j=0}^\infty$
are both $m\times m$ Schur sequences 
(see, e.g., \cite[Theorem 5.1.1.]{DFK} or Theorem \ref{T62} with $J=I_m$). 
In particular, for each $j\in\N_0$ we have $\|B_j\|\le1$ and $\|B_j^{(n)}\|\le1$.
Moreover, from \cite[Remark 6.1]{FKR1} we get $B_j=B_j^{(n)}$ for each $j\in\N_{0,n}$.
Consequently, for each $w\in\D$ we get
$$ \|g(w)-g_n(w)\|\le\sum_{j=n+1}^\infty\|B_j-B_j^{(n)}\||w|^j\le\frac{2|w|^{n+1}}{1-|w|} \ .$$
Thus, we obtain that
$\lim_{n\to\infty}g_n(w)=g(w)$
holds for each $w\in\D$ and that this convergence is uniform in each compact subset of $\D$.
From \cite[Proposition 3.4]{FKR1} we get 
\begin{equation}\label{CPF-P2-5}
\mathbb H_f=\{w\!\in\!\D\!:\det(\mathbf Q_Jg(w)+\mathbf P_J)\!\ne\!0\}, \ \ 
\mathbb H_{f_n}=\{w\!\in\!\D\!:\det(\mathbf Q_Jg_n(w)+\mathbf P_J)\!\ne\!0\},
\end{equation}
as well as
\begin{equation}\label{CPF-P2-6}
f=(\mathbf P_Jg+\mathbf Q_J)(\mathbf Q_Jg+\mathbf P_J)^{-1} \quad\mbox{and}\quad
f_n=(\mathbf P_Jg_n+\mathbf Q_J)(\mathbf Q_Jg_n+\mathbf P_J)^{-1}
\end{equation}
for each $n\in\N_0$.
Let $K$ be a compact subset of $\mathbb H_f$. Since $g$ is continuous, 
the image $g(K)$ of $K$ under $g$ is a compact subset 
of $\C^{m\times m}$. 
Let $\cQ_{(\mathbf Q_J,\mathbf P_J)}:=\{X\in\C^{m\times m}:\det(\mathbf Q_JX+\mathbf P_J)\ne0\}$.
In view of (\ref{CPF-P2-5}) we have $g(K)\sq\cQ_{(\mathbf Q_J,\mathbf P_J)}$.
Because $\cQ_{(\mathbf Q_J,\mathbf P_J)}$ is an open subset of $\C^{m\times m}$,
there exists a positive real number $\eta$ such that
\begin{equation}\label{CPF-P2-7}
g(K)\sq g(K)+\mathbf B_\eta\sq\cQ_{(\mathbf Q_J,\mathbf P_J)}
\end{equation}
holds where $\mathbf B_\eta:=\{X\in\C^{m\times m}:\|X\|\le\eta\}$. Moreover,
$g(K)+\mathbf B_\eta$ is a compact subset of $\C^{m\times m}$. 
Since the sequence $(g_n)_{n\in\N_0}$ converges uniformly to $g$ in $K$,
there exists an $n_0\in\N$ such that for each $n\in\N_{n_0,\infty}$ and each $w\in K$
the relation $g_n(w)\in g(w)+\mathbf B_\eta$ holds. 
Hence 
\begin{equation}\label{CPF-P2-8}
g_n(K)\sq g(K)+\mathbf B_\eta
\end{equation}
is valid for each integer $n$ with $n\ge n_0$.
Taking into account (\ref{CPF-P2-5}), (\ref{CPF-P2-7}), and (\ref{CPF-P2-8}), we get that
$f_n$ is holomorphic in $K$ for every integer $n$ with $n\ge n_0$. 
Let the map $F:\cQ_{(\mathbf Q_J,\mathbf P_J)}\rightarrow\C^{m\times m}$
be defined by $X\mapsto(\mathbf P_JX+\mathbf Q_J)(\mathbf Q_JX+\mathbf P_J)^{-1}$.
Then, in view of  (\ref{CPF-P2-8}), (\ref{CPF-P2-7}), and (\ref{CPF-P2-6})
we obtain
\begin{equation}\label{CPF-P2-9}
f(w)=F(g(w)) \qquad\mbox{and}\qquad f_n(w)=F(g_n(w))
\end{equation}
for each $w\in K$ and each integer $n$ with $n\ge n_0$.
Obviously, $F$ is continuous in $\cQ_{(\mathbf Q_J,\mathbf P_J)}$
and, therefore, uniformly continuous in the compact set $g(K)+\mathbf B_\eta$. 
Consequently, because of (\ref{CPF-P2-7}), (\ref{CPF-P2-8}), and (\ref{CPF-P2-9}), 
the uniform convergence $g_n\to g$ in $K$ implies that the sequence
$(f_n)_{n\ge n_0}$ converges uniformly to $f$ in $K$. 
\qed\end{pf}

In particular we see that each $f\in\cP_{J,0}(\D)$ can be approximated by its associated 
sequence of $J$-central $J$-Potapov functions.
This is a generalization of the corresponding well-known result for $m\times m$ Schur functions in $\D$ 
(see, e.g., \cite[Theorem 3.5.3]{DFK}).

\begin{cor}\label{CPF-C0}
Let $J$ be an $m\times m$ signature matrix, and let $f\in\cP_{J,0}(\D)$. Let (\ref{NrTS})
be the Taylor series representation of $f$ in some neighborhood of $0$. For each $n\in\N_0$,
let $f_{c,n}$ be the $J$-central $J$-Potapov function corresponding to $(A_j)_{j=0}^n$.
Then, for each compact subset $K$ of $\mathbb H_f$, there exists a nonnegative integer $n_0$ such that $f_{c,n}$ is holomorphic in $K$ for every integer $n$ with $n\ge n_0$. Moreover,
$$ \lim_{n\to\infty}f_{c,n}(w)=f(w) $$
holds for each $w\in\mathbb H_f$. This convergence is uniform in each compact subset of $\mathbb H_f$.
\end{cor}

\begin{pf}
Combine Proposition \ref{CPF-P0} and Proposition \ref{CPF-P2}.
\qed\end{pf}

Let us now introduce some notations. Let $J$ be an $m\times m$ signature matrix, let $n\in\N$,
and let $(A_j)_{j=0}^n$ be a $J$-Potapov sequence. Then let 
\begin{equation}\label{NrYnJ}
\cY_{n,J}:=\{V\in\C^{nm\times m}:Q_{n-1,J}V=S_{n-1}^*J_{[n-1]}y_n\}
\end{equation}
and
\begin{equation}\label{NrZnJ}
\cZ_{n,J}:=\{W\in\C^{m\times nm}:WP_{n-1,J}=z_nJ_{[n-1]}S_{n-1}^*\}.
\end{equation}

\begin{rem}\label{LRQ-R2}
Let $J$ be an $m\times m$ signature matrix, let $n\in\N$,
and let $(A_j)_{j=0}^n$ be a $J$-Potapov sequence. 
Furthermore, let
\begin{equation}\label{NrVnWn}
V_n^\Box:=Q_{n-1,J}^+S_{n-1}^*J_{[n-1]}y_n
\quad\mbox{and}\quad
W_n^\Box:=z_nJ_{[n-1]}S_{n-1}^*P_{n-1,J}^+.
\end{equation}
Then $V_n^\Box\in\cY_{n,J}$ and $W_n^\Box\in\cZ_{n,J}$ hold (see \cite[Remark 2.5]{FKRS}).
\end{rem}

Henceforth, whenever an $m\times m$ signature matrix $J$, an $n\in\N_0$, and some matrix polynomials 
$\pi_{n,J}$, $\rho_{n,J}$, $\sigma_{n,J}$, $\tau_{n,J}$ (defined on $\C$) are given, then 
$\pi_{n,J,\D}$ (respectively, $\rho_{n,J,\D}$, $\sigma_{n,J,\D}$, $\tau_{n,J,\D}$) always stands for the restriction of $\pi_{n,J}$ (respectively, $\rho_{n,J}$, $\sigma_{n,J}$, $\tau_{n,J}$) onto $\D$. 

\begin{thm}\label{LRQ-T1}
Let $J$ be an $m\times m$ signature matrix, let $n\in\N_0$, and let $(A_j)_{j=0}^n$ be a $J$-Potapov
sequence. Denote by $f_{c,n}$ the $J$-central $J$-Potapov function corresponding to $(A_j)_{j=0}^n$.
If $n\in\N$, then let $V_n\in\cY_{n,J}$ and $W_n\in\cZ_{n,J}$. 
Let the $m\times m$ matrix polynomials 
$\pi_{n,J}$, $\rho_{n,J}$,
$\sigma_{n,J}$, and $\tau_{n,J}$
be given by
\begin{equation}\label{LRQ-T1-1}
\pi_{n,J}(w):=\left\{\begin{array}{cl} A_0, &\mbox{if }n=0\\
               A_0+w e_{n-1,m}(w)(y_n+S_{n-1}V_n), &\mbox{if }n\in\N,\end{array}\right.
\end{equation}
\begin{equation}\label{LRQ-T1-2}
\rho_{n,J}(w):=\left\{\begin{array}{cl} I_m, &\mbox{if }n=0\\
                             I_m+w e_{n-1,m}(w)V_n, &\mbox{if }n\in\N,\end{array}\right.
\end{equation}
\begin{equation}\label{LRQ-T2-1}
\sigma_{n,J}(w):=\left\{\begin{array}{cl} A_0, &\mbox{if }n=0\\
          (W_nS_{n-1}+z_n)w\eps_{n-1,m}(w)+A_0,&\mbox{if }n\in\N, \end{array}\right.
\end{equation}
and
\begin{equation}\label{LRQ-T2-2}
\tau_{n,J}(w):=\left\{\begin{array}{cl} I_m, &\mbox{if }n=0\\
                     W_nw\eps_{n-1,m}(w)+I_m,&\mbox{if }n\in\N \end{array}\right.
\end{equation}
for each $w\in\C$. Then $f_{c,n}$ admits the representations
\begin{equation}\label{LRQ-T1-5}
f_{c,n}=\pi_{n,J,\D}\,\rho_{n,J,\D}^{-1}\quad\mbox{and}\quad f_{c,n}=\tau_{n,J,\D}^{-1}\,\sigma_{n,J,\D}.
\end{equation}
\end{thm}

A proof of Theorem \ref{LRQ-T1} is given in \cite[Theorems 2.7 and 2.8]{FKRS}.

%

In the sequel, we will use the following notation. If $J$ is an $m\times m$ signature matrix, 
if $n\in\N_0$, and if $(A_j)_{j=0}^{n+1}$ is a $J$-Potapov sequence, then we will work with the sets
\begin{equation}\label{RF-E1}
\cL_{n+1,J}:=\{t\in\C^{m\times m}:L_{n+1,J}t=A_{n+1}-M_{n+1,J}\}
\end{equation}
and
\begin{equation}\label{RF-E2}
\cR_{n+1,J}:=\{u\in\C^{m\times m}:uR_{n+1,J}=A_{n+1}-M_{n+1,J}\}.
\end{equation}

\begin{rem}\label{RF-R1}
Let $J$ be an $m\times m$ signature matrix, let $n\in\N_0$, and let $(A_j)_{j=0}^{n+1}$ be a $J$-Potapov
sequence. Let 
\begin{equation}\label{RF-R1-1}
t_{n+1}:=L_{n+1,J}^+(A_{n+1}-M_{n+1,J}) \quad\mbox{and}\quad u_{n+1}:=(A_{n+1}-M_{n+1,J})R_{n+1,J}^+.
\end{equation}
In view of Theorem \ref{T}, there is a contractive $m\times m$ matrix $K$ such that the identity
$A_{n+1}-M_{n+1,J}=\sqrt{L_{n+1,J}}K\sqrt{R_{n+1,J}}$ holds. Consequently, $t_{n+1}\in\cL_{n+1,J}$
and $u_{n+1}\in\cR_{n+1,J}$.
\end{rem}


\begin{prop}\label{CPF-P1}
Let $J$ be an $m\times m$ signature matrix, let $n\in\N_0$, let $k\in\N$, 
and let $(A_j)_{j=0}^{n+k}$ be a $J$-Potapov sequence. If $n\ge1$, then let $V_n\in\cY_{n,J}$ 
and $W_n\in\cZ_{n,J}$. Furthermore, let $\pi_{n,J}$, $\rho_{n,J}$, $\sigma_{n,J}$, and $\tau_{n,J}$
be the matrix polynomials given by (\ref{LRQ-T1-1}), (\ref{LRQ-T1-2}), (\ref{LRQ-T2-1}), and (\ref{LRQ-T2-2}). For $s\in\N_{0,k-1}$ let $\pi_{n+s+1,J}$, $\rho_{n+s+1,J}$, $\sigma_{n+s+1,J}$, and $\tau_{n+s+1,J}$ be the matrix polynomials which are recursively defined by
\begin{equation}\label{CPF-P1-1}
\pi_{n+s+1,J}(w):=\pi_{n+s,J}(w)+wJ\tilde\tau_{n+s,J}^{[n+s]}(w)t_{n+s+1},
\end{equation}
\begin{equation}\label{CPF-P1-2}
\rho_{n+s+1,J}(w):=\rho_{n+s,J}(w)+wJ\tilde\sigma_{n+s,J}^{[n+s]}(w)t_{n+s+1},
\end{equation}
\begin{equation}\label{CPF-P1-3}
\sigma_{n+s+1,J}(w):=\sigma_{n+s,J}(w)+u_{n+s+1}w\tilde\rho_{n+s,J}^{[n+s]}(w)J,
\end{equation}
and
\begin{equation}\label{CPF-P1-4}
\tau_{n+s+1,J}(w):=\tau_{n+s,J}(w)+u_{n+s+1}w\tilde\pi_{n+s,J}^{[n+s]}(w)J
\end{equation}
for each $w\!\in\!\C$, where $t_{n+s+1}\!:=\!L_{n+s+1,J}^+(A_{n+s+1}-M_{n+s+1,J})$ and $u_{n+s+1}\!:=\!(A_{n+s+1}-M_{n+s+1,J})R_{n+s+1,J}^+$.
Then, for each $s\in\N_{0,k-1}$, 
both $\det\rho_{n+s+1,J}$ and $\det\tau_{n+s+1,J}$ do not vanish identically in $\D$, and
the $J$-central $J$-Potapov function $f_{c,n+s+1}$
corresponding to $(A_j)_{j=0}^{n+s+1}$ admits the representations
$$ f_{c,n+s+1}=\pi_{n+s+1,J,\D}\,\rho_{n+s+1,J,\D}^{-1}\quad\mbox{and}\quad f_{c,n+s+1}=\tau_{n+s+1,J,\D}^{-1}\,\sigma_{n+s+1,J,\D}.$$
\end{prop}

\begin{pf}
Apply Remark \ref{RF-R1}, \cite[Proposition 3.4, Remark 3.2, and Lemma 3.3]{FKRS}, 
and Theorem \ref{LRQ-T1}.
\qed\end{pf}

In the following, for each $k\in\N_0$ let the map $\cE_k:\D\rightarrow\D$ be defined by $w\mapsto w^k$.

\begin{cor}\label{CPF-C1}
Let $J$ be an $m\times m$ signature matrix, let $n\in\N_0$, 
and let $(A_j)_{j=0}^n$ be a $J$-Potapov sequence. If $n\ge1$, then let $V_n\in\cY_{n,J}$ 
and $W_n\in\cZ_{n,J}$. Furthermore, let $\pi_{n,J}$, $\rho_{n,J}$, $\sigma_{n,J}$, 
and $\tau_{n,J}$ be given by (\ref{LRQ-T1-1}), (\ref{LRQ-T1-2}), (\ref{LRQ-T2-1}), and (\ref{LRQ-T2-2}).
Let $K$ be a contractive matrix from $\C^{m\times m}$ and let
$$ A_{n+1}:=M_{n+1,J}+\sqrt{L_{n+1,J}}K\sqrt{R_{n+1,J}}.$$
Then $(A_j)_{j=0}^{n+1}$ is a $J$-Potapov sequence.
Furthermore, each of the functions\linebreak[0] $\det\!\big(\!\cE_1\!J\tilde\sigma_{n,J,\D}^{[n]}\sqrt{\!L_{n+1,J}}^+\!K\!\sqrt{\!R_{n+1,J}} \linebreak[0]
 + \linebreak[0] 
 \rho_{n,J,\D}\!\big)\!$
and \linebreak[0]
$\det\!\big(\!\cE_1\!\sqrt{L_{n+1,J}}K\!\sqrt{R_{n+1,J}}^+\!\tilde\pi_{n,J,\D}^{[n]}J \linebreak[0]
 + \linebreak[0]
 \tau_{n,J,\D}\!\big)\!$
does not vanish identically in $\D$. 
Moreover, the $J$-central
$J$-Potapov function $f_{c,n+1}$ corresponding to $(A_j)_{j=0}^{n+1}$
admits the representations
\begin{align*}
f_{c,n+1}&=\big(\cE_1J\tilde\tau_{n,J,\D}^{[n]}\sqrt{L_{n+1,J}}^+K\sqrt{R_{n+1,J}}+\pi_{n,J,\D}\big)\\
& \qquad \cdot
    \big(\cE_1J\tilde\sigma_{n,J,\D}^{[n]}\sqrt{L_{n+1,J}}^+K\sqrt{R_{n+1,J}}+\rho_{n,J,\D}\big)^{-1} 
\end{align*}
and
\begin{align*}
f_{c,n+1}&=\big(\cE_1\sqrt{L_{n+1,J}}K\sqrt{R_{n+1,J}}^+\tilde\pi_{n,J,\D}^{[n]}J+\tau_{n,J,\D}\big)^{-1} \\
&\qquad\cdot
  \big(\cE_1\sqrt{L_{n+1,J}}K\sqrt{R_{n+1,J}}^+\tilde\rho_{n,J,\D}^{[n]}J+\sigma_{n,J,\D}\big).
\end{align*}
\end{cor}

\begin{pf}
According to Theorem \ref{T}, $(A_j)_{j=0}^{n+1}$ is a $J$-Potapov sequence. 
Furthermore, from $\rho_{n,J}(0)=I_m$ and $\tau_{n,J}(0)=I_m$ we see that
each of the functions \linebreak[0] $\det\!\big(\!\cE_1J\tilde\sigma_{n,J,\D}^{[n]}\sqrt{L_{n+1,J}}^+\!K\!\sqrt{R_{n+1,J}} \linebreak[0]
 + \linebreak[0]
 \rho_{n,J,\D}\!\big)\!$
and  \linebreak[0]
$\det\!\big(\!\cE_1\!\sqrt{L_{n+1,J}}K\!\sqrt{R_{n+1,J}}^+\!\tilde\pi_{n,J,\D}^{[n]}J \linebreak[0]
 + \linebreak[0]
 \tau_{n,J,\D}\!\big)$
does not vanish identically in $\D$.
In view of $L_{n+1,J}^+\sqrt{L_{n+1,J}}=\sqrt{L_{n+1,J}}^+$ and $\sqrt{R_{n+1,J}}R_{n+1,J}^+=\sqrt{R_{n+1,J}}^+$, 
the matrices $t_{n+1}$ and $u_{n+1}$ given by (\ref{RF-R1-1}) satisfy
$t_{n+1}=\sqrt{L_{n+1,J}}^+K\sqrt{R_{n+1,J}}$ and $u_{n+1}=$ $\sqrt{L_{n+1,J}}K\sqrt{R_{n+1,J}}^+$.
Thus, an application of Proposition \ref{CPF-P1} (with $s=0$) completes the proof.
\qed\end{pf}
\section{Description of the solution set $\cP_{J,0}\left[\D,(A_j)_{j=0}^n\right]$}

The main goal of this section is to prove Theorem \ref{S1-T1}. In fact, the results obtained in
this section are even more general than the description of the
solution set of problem (P) which is given in Theorem \ref{S1-T1}.
This section generalizes the corresponding results for $m\times m$ Schur functions, which were obtained
in \cite{FKL}, to the case of $J$-Potapov functions.

In this section, we will continue to use the notation $\cE_k$ introduced above, i.e., for each 
$k\in\N_0$ let $\cE_k:\D\rightarrow\D$ be defined by $w\mapsto w^k$.

\begin{rem}\label{DSS-R1}
Let $J$ be an $m\times m$ signature matrix, let $n\in\N_0$, 
and let $(A_j)_{j=0}^n$ be a $J$-Potapov sequence. If $n\ge1$, then let $V_n\in\cY_{n,J}$ 
and $W_n\in\cZ_{n,J}$. Furthermore, let $\pi_{n,J}$, $\rho_{n,J}$, $\sigma_{n,J}$, and $\tau_{n,J}$
be given by (\ref{LRQ-T1-1}), (\ref{LRQ-T1-2}), (\ref{LRQ-T2-1}), and (\ref{LRQ-T2-2}). 
Let $S\in\cS_{m\times m}(\D)$. 
Then, because of $\rho_{n,J}(0)=I_m$ and $\tau_{n,J}(0)=I_m$, each of the functions $\det\big(\cE_1J\tilde\sigma_{n,J,\D}^{[n]}\sqrt{L_{n+1,J}}^+S\sqrt{R_{n+1,J}}+\rho_{n,J,\D}\big)$
and 
$\det\big(\cE_1\sqrt{\!L_{n+1,J}}S\sqrt{\!R_{n+1,J}}^+\tilde\pi_{n,J,\D}^{[n]}J+\tau_{n,J,\D}\big)$ does not vanish identically in $\D$. 
Let $\cA_1$ (resp., $\cA_2$) be the set of all zeros of 
$\det\!\big(\cE_1J\tilde\sigma_{n,J,\D}^{[n]}\sqrt{\!L_{n+1,J}}^+\!S\sqrt{\!R_{n+1,J}}+\rho_{n,J,\D}\big)$
(resp.,
$\det\big(\cE_1\sqrt{L_{n+1,J}}S\sqrt{R_{n+1,J}}^+\tilde\pi_{n,J,\D}^{[n]}J+\tau_{n,J,\D}\big)$).
Setting
\begin{align}\label{DSS-R1-1}
f&:=\big(\cE_1J\tilde\tau_{n,J,\D}^{[n]}\sqrt{L_{n+1,J}}^+S\sqrt{R_{n+1,J}}+\pi_{n,J,\D}\big) \nonumber\\
&\qquad\cdot
\big(\cE_1J\tilde\sigma_{n,J,\D}^{[n]}\sqrt{L_{n+1,J}}^+S\sqrt{R_{n+1,J}}+\rho_{n,J,\D}\big)^{-1} 
\end{align}
we see from Corollary \ref{CPF-C1} that
\begin{align*}
f(w)&=\big(w\sqrt{L_{n+1,J}}S(w)\sqrt{R_{n+1,J}}^+\tilde\pi_{n,J}^{[n]}(w)J+\tau_{n,J}(w)\big)^{-1} \\
&\qquad\cdot
   \big(w\sqrt{L_{n+1,J}}S(w)\sqrt{R_{n+1,J}}^+\tilde\rho_{n,J}^{[n]}(w)J+\sigma_{n,J}(w)\big)
\end{align*}
holds for each $w\in\D\setminus(\cA_1\cup\cA_2)$. Since $\cA_1\cup\cA_2$ is a discrete subset of $\D$,
this implies
\begin{align}\label{DSS-R1-2}
f&=\big(\cE_1\sqrt{L_{n+1,J}}S\sqrt{R_{n+1,J}}^+\tilde\pi_{n,J,\D}^{[n]}J+\tau_{n,J,\D}\big)^{-1}\nonumber\\
&\qquad\cdot
   \big(\cE_1\sqrt{L_{n+1,J}}S\sqrt{R_{n+1,J}}^+\tilde\rho_{n,J,\D}^{[n]}J+\sigma_{n,J,\D}\big).
\end{align}
If $n\ge1$, and if $\fpi_{n,J}$, $\frho_{n,J}$, $\fsigma_{n,J}$, and $\ftau_{n,J}$ are further
matrix polynomials which can be represented, for each $w\in\C$, via
\begin{equation}\label{DSS-R1-3}
\fpi_{n,J}(w)=A_0+w e_{n-1,m}(w)(y_n+S_{n-1}\mathbf V_n),\ \: 
 \frho_{n,J}(w)=I_m+w e_{n-1,m}(w)\mathbf V_n, 
\end{equation}
\begin{equation}\label{DSS-R1-4}
\fsigma_{n,J}(w)=(\mathbf W_nS_{n-1}+z_n)w\eps_{n-1,m}(w)+A_0,\ \: 
\ftau_{n,J}(w)\!=\!\mathbf W_nw\eps_{n-1,m}(w)+I_m 
\end{equation}
with some $\mathbf V_n\in\cY_{n,J}$ and $\mathbf W_n\in\cZ_{n,J}$, then Corollary \ref{CPF-C1}
and a similar consideration as above provide us
\begin{align*}
f&=\big(\cE_1J\tilde\ftau_{n,J,\D}^{[n]}\sqrt{L_{n+1,J}}^+S\sqrt{R_{n+1,J}}+\fpi_{n,J,\D}\big) \\
&\qquad\cdot
\big(\cE_1J\tilde\fsigma_{n,J,\D}^{[n]}\sqrt{L_{n+1,J}}^+S\sqrt{R_{n+1,J}}+\frho_{n,J,\D}\big)^{-1}
\end{align*}
where $\fpi_{n,J,\D}$ (respectively, $\frho_{n,J,\D}$, $\fsigma_{n,J,\D}$, $\ftau_{n,J,\D}$) denotes the restriction of $\fpi_{n,J}$ (respectively, $\frho_{n,J}$, $\fsigma_{n,J}$, $\ftau_{n,J}$) onto $\D$. 
\end{rem}

Now we will show that for every $S\in\cS_{m\times m}(\D)$ the matrix-valued function $f$ constructed in
Remark \ref{DSS-R1} is a solution of the interpolation problem (P) stated in section \ref{S1}.

\begin{thm}\label{DSS-T1}
Let $J$ be an $m\times m$ signature matrix, let $n\in\N_0$, 
and let $(A_j)_{j=0}^n$ be a $J$-Potapov sequence. If $n\ge1$, then let $V_n\in\cY_{n,J}$ 
and $W_n\in\cZ_{n,J}$. Furthermore, let $\pi_{n,J}$, $\rho_{n,J}$, $\sigma_{n,J}$, and $\tau_{n,J}$
be given by (\ref{LRQ-T1-1}), (\ref{LRQ-T1-2}), (\ref{LRQ-T2-1}), and (\ref{LRQ-T2-2}). 
Let $S\in\cS_{m\times m}(\D)$, and let the matrix-valued function $f$ be given by (\ref{DSS-R1-1}).
Then $f$ belongs to $\cP_{J,0}\left[\D,(A_j)_{j=0}^n\right]$ and admits the representation
(\ref{DSS-R1-2}).
\end{thm}

\begin{pf}
In view of Remark \ref{DSS-R1} we see that $f$ is a well-defined matrix-valued function which is meromorphic in $\D$. 
Moreover, because of $\rho_{n,J}(0)=I_m$ we have $0\in\mathbb H_f$.
Let $\cA_1$ (respectively, $\cA_2$) be the set of all zeros of the function
$\det\big(\cE_1J\tilde\sigma_{n,J,\D}^{[n]}\sqrt{L_{n+1,J}}^+S\sqrt{R_{n+1,J}}+\rho_{n,J,\D}\big)$ (respectively, of the function $\det\big(\cE_1\sqrt{L_{n+1,J}}S\sqrt{R_{n+1,J}}^+\tilde\pi_{n,J,\D}^{[n]}J+\tau_{n,J,\D}\big)$).
Then $\cA_1\cup\cA_2$ is a discrete subset of $\D$. 
Now let $w_0\in\D\setminus(\cA_1\cup\cA_2)$.
Then $K:=S(w_0)$ is a contractive $m\times m$ matrix.
Consequently, from Corollary \ref{CPF-C1} we obtain that the matrix $f(w_0)$ is $J$-contractive.  
By continuity we get that $f(w)$ is $J$-contractive for all $w\in\mathbb H_f$, i.e., we have
$f\in\cP_{J,0}(\D)$. 
Now we prove that $f$ fulfills condition (\ref{ICP}) for each $j\in\N_{0,n}$.
Let $F_n:=\sqrt{L_{n+1,J}}^+S\sqrt{R_{n+1,J}}$.
Because of $\rho_{n,J}(0)=I_m$ and $\tau_{n,J}(0)=I_m$ we see that 
$$ h_n:=\tau_{n,J,\D}^{-1}L_{n+1,J}F_n(\cE_1J\tilde\sigma_{n,J,\D}^{[n]}F_n+\rho_{n,J,\D})^{-1} $$
is a well-defined matrix-valued function which is meromorphic in $\D$ and holomorphic at $0$.
According to Theorem \ref{LRQ-T1}, the $J$-central $J$-Potapov function $f_{c,n}$ corresponding to
$(A_j)_{j=0}^n$ admits the representations
\begin{equation}\label{DSS-T1-7}
f_{c,n}=\pi_{n,J,\D}\rho_{n,J,\D}^{-1} \qquad\mbox{and}\qquad
 f_{c,n}=\tau_{n,J,\D}^{-1}\sigma_{n,J,\D}. 
\end{equation}
Moreover, \cite[Proposition 2.11]{FKRS} yields
\begin{equation}\label{DSS-T1-8}
\tau_{n,J,\D}J\tilde\tau_{n,J,\D}^{[n]}-\sigma_{n,J,\D}J\tilde\sigma_{n,J,\D}^{[n]}=\cE_nL_{n+1,J}.
\end{equation}
Using (\ref{DSS-R1-1}), (\ref{DSS-T1-7}), and (\ref{DSS-T1-8}) we get 
\begin{align}\label{DSS-T1-9}
f-f_{c,n}
&= \tau_{n,J,\D}^{-1}[\tau_{n,J,\D}(\cE_1J\tilde\tau_{n,J,\D}^{[n]}F_n+\pi_{n,J,\D})
  -\sigma_{n,J,\D}(\cE_1J\tilde\sigma_{n,J,\D}^{[n]}F_n+\rho_{n,J,\D})] \nonumber\\
&\qquad\cdot    (\cE_1J\tilde\sigma_{n,J,\D}^{[n]}F_n+\rho_{n,J,\D})^{-1} \nonumber
\displaybreak[0]\\
&= \tau_{n,J,\D}^{-1}(\cE_1\tau_{n,J,\D}J\tilde\tau_{n,J,\D}^{[n]}F_n
  -\cE_1\sigma_{n,J,\D}J\tilde\sigma_{n,J,\D}^{[n]}F_n)
   (\cE_1J\tilde\sigma_{n,J,\D}^{[n]}F_n+\rho_{n,J,\D})^{-1} \nonumber
\displaybreak[0]\\
&= \tau_{n,J,\D}^{-1}\cE_{n+1}L_{n+1,J}F_n(\cE_1J\tilde\sigma_{n,J,\D}^{[n]}F_n+\rho_{n,J,\D})^{-1}
= \cE_{n+1}h_n.
\end{align}
Let $h_n(w)=\sum_{j=0}^\infty H_jw^j$ be the Taylor series representation of $h_n$ in some
neighborhood of $0$. Then (\ref{DSS-T1-9}) implies
$ f(w)-f_{c,n}(w)=\sum_{j=n+1}^\infty H_{j-n-1}w^j $
for each $w$ belonging to some neighborhood of $0$. Consequently, because of 
$f_{c,n}\in\cP_{J,0}[\D,(A_j)_{j=0}^n]$ and $f\in\cP_{J,0}(\D)$ we obtain $f\in\cP_{J,0}[\D,(A_j)_{j=0}^n]$.
Finally, having in mind Remark \ref{DSS-R1} we see that $f$ can be represented via (\ref{DSS-R1-2}).
\qed\end{pf}

The following considerations show that the parameter function $S\in\cS_{m\times m}(\D)$
in the linear fractional transformations stated in (\ref{DSS-R1-1}) and (\ref{DSS-R1-2})
is generally not uniquely determined by $f$.

\begin{rem}\label{DSS-R2}
Let $J$ be an $m\times m$ signature matrix, let $n\in\N_0$, 
and let $(A_j)_{j=0}^n$ be a $J$-Potapov sequence. If $n\ge1$, then let $V_n\in\cY_{n,J}$ 
and $W_n\in\cZ_{n,J}$. Furthermore, let $\pi_{n,J}$, $\rho_{n,J}$, $\sigma_{n,J}$, and $\tau_{n,J}$
be given by (\ref{LRQ-T1-1}), (\ref{LRQ-T1-2}), (\ref{LRQ-T2-1}), and (\ref{LRQ-T2-2}). 
Let $S\in\cS_{m\times m}(\D)$, and let the matrix-valued function $f$ be given by
(\ref{DSS-R1-1}). 
Since $L_{n+1,J}L_{n+1,J}^+$ and $R_{n+1,J}^+R_{n+1,J}$ are orthoprojection matrices,
the matrix-valued function $S^\sharp:=L_{n+1,J}L_{n+1,J}^+SR_{n+1,J}^+R_{n+1,J}$
belongs to $\cS_{m\times m}(\D)$ and satisfies the identity 
$L_{n+1,J}L_{n+1,J}^+S^\sharp R_{n+1,J}^+R_{n+1,J}=S^\sharp$.
Moreover, it is readily checked that $f$ admits the representation
\begin{align*}
f&=
\big(\cE_1J\tilde\tau_{n,J,\D}^{[n]}\sqrt{L_{n+1,J}}^+S^\sharp\sqrt{R_{n+1,J}}+\pi_{n,J,\D}\big)\\
&\qquad\cdot    \big(\cE_1J\tilde\sigma_{n,J,\D}^{[n]}\sqrt{L_{n+1,J}}^+S^\sharp\sqrt{R_{n+1,J}}+\rho_{n,J,\D}\big)^{-1}.
\end{align*}
\end{rem}

\begin{prop}\label{DSS-P1}
Let $J$ be an $m\times m$ signature matrix, let $n\in\N_0$, 
and let $(A_j)_{j=0}^n$ be a $J$-Potapov sequence. If $n\ge1$, then let $V_n\in\cY_{n,J}$ 
and $W_n\in\cZ_{n,J}$. Furthermore, let $\pi_{n,J}$, $\rho_{n,J}$, $\sigma_{n,J}$, and $\tau_{n,J}$
be given by (\ref{LRQ-T1-1}), (\ref{LRQ-T1-2}), (\ref{LRQ-T2-1}), and (\ref{LRQ-T2-2}). 
For $j\in\{1,2\}$, let $S_j\in\cS_{m\times m}(\D)$ and let (in view of Remark \ref{DSS-R1})
the matrix-valued meromorphic functions $f_j$ and $\tilde f_j$ be given by
\begin{align*}
f_j&:=
\big(\cE_1J\tilde\tau_{n,J,\D}^{[n]}\sqrt{L_{n+1,J}}^+S_j\sqrt{R_{n+1,J}}+\pi_{n,J,\D}\big)\\
&\qquad\cdot    \big(\cE_1J\tilde\sigma_{n,J,\D}^{[n]}\sqrt{L_{n+1,J}}^+S_j\sqrt{R_{n+1,J}}+\rho_{n,J,\D}\big)^{-1}
\end{align*}
and
\begin{align*}
\tilde f_j
&:=\big(\cE_1\sqrt{L_{n+1,J}}S_j\sqrt{R_{n+1,J}}^+\tilde\pi_{n,J,\D}^{[n]}J+\tau_{n,J,\D}\big)^{-1}\\
&\qquad\cdot    \big(\cE_1\sqrt{L_{n+1,J}}S_j\sqrt{R_{n+1,J}}^+\tilde\rho_{n,J,\D}^{[n]}J+\sigma_{n,J,\D}\big).
\end{align*}
\begin{enumerate}
\item[(a)] The following statements are equivalent:
\begin{enumerate}
\item[(i)] $\sqrt{L_{n+1,J}}S_1\sqrt{R_{n+1,J}}=\sqrt{L_{n+1,J}}S_2\sqrt{R_{n+1,J}}$.
\item[(ii)] $f_1=f_2$.
\item[(iii)] $\tilde f_1=\tilde f_2$.
\end{enumerate}
\item[(b)] Let $f_1=f_2$ or $\tilde f_1=\tilde f_2$. If $S_j=L_{n+1,J}L_{n+1,J}^+S_jR_{n+1,J}^+R_{n+1,J}$
holds for each $j\in\{1,2\}$, then $S_1=S_2$.
\end{enumerate}
\end{prop}

\begin{pf}
The equivalence of (ii) and (iii) is an immediate consequence of Remark \ref{DSS-R1}.
In view of $\sqrt{L_{n+1,J}}^+=L_{n+1,J}^+\sqrt{L_{n+1,J}}$, the implication (i)$\Rightarrow$(ii)
is obvious. Now suppose that (ii) is satisfied. 
From Theorem \ref{LRQ-T1} we obtain
\begin{equation}
\tau_{n,J,\D}\pi_{n,J,\D}=\sigma_{n,J,\D}\rho_{n,J,\D}
\quad\mbox{and}\quad
\tilde\pi_{n,J,\D}^{[n]}\tilde\tau_{n,J,\D}^{[n]}=\tilde\rho_{n,J,\D}^{[n]}\tilde\sigma_{n,J,\D}^{[n]}.
\end{equation}
Moreover, \cite[Proposition 2.11]{FKRS} provides us
(\ref{DSS-T1-8})
and
\begin{equation}\label{DSS-P1-2}
\tilde\rho_{n,J,\D}^{[n]}J\rho_{n,J,\D}-\tilde\pi_{n,J,\D}^{[n]}J\pi_{n,J,\D}=\cE_nR_{n+1,J}.
\end{equation}
By succinct setting
$X_1\!:=\!\cE_1J\tilde\sigma_{n,J,\D}^{[n]}\sqrt{L_{n+1,J}}^+\!S_1\!\sqrt{R_{n+1,J}}+\rho_{n,J,\D}$ and
$X_2\!:=\!\cE_1\!\sqrt{L_{n+1,J}}S_2\!\sqrt{R_{n+1,J}}^+\tilde\pi_{n,J,\D}^{[n]}J + \tau_{n,J,\D}$,
from (ii), Remark \ref{DSS-R1}, (\ref{DSS-T1-8}), and (\ref{DSS-P1-2}) we get
\begin{align*}
0&=f_1-f_2=f_1-\tilde f_2 \\
&= X_2^{-1}\Big[X_2
    \Big(\cE_1J\tilde\tau_{n,J,\D}^{[n]}\sqrt{L_{n+1,J}}^+S_1\sqrt{R_{n+1,J}}+\pi_{n,J,\D}\Big) \\
& \qquad  -\Big(\cE_1\sqrt{L_{n+1,J}}S_2\sqrt{R_{n+1,J}}^+\tilde\rho_{n,J,\D}^{[n]}J+\sigma_{n,J,\D}\Big)X_1\Big]
 X_1^{-1}
\displaybreak[0] \\
&= X_2^{-1}\Big[
 \cE_2\sqrt{L_{n+1,J}}S_2\sqrt{R_{n+1,J}}^+\!
 (\tilde\pi_{n,J,\D}^{[n]}\tilde\tau_{n,J,\D}^{[n]}-\tilde\rho_{n,J,\D}^{[n]}\tilde\sigma_{n,J,\D}^{[n]})
  \sqrt{L_{n+1,J}}^+S_1\sqrt{R_{n+1,J}} \\
& \qquad 
 +\cE_1(\tau_{n,J,\D}J\tilde\tau_{n,J,\D}^{[n]}-\sigma_{n,J,\D}J\tilde\sigma_{n,J,\D}^{[n]})
  \sqrt{L_{n+1,J}}^+S_1\sqrt{R_{n+1,J}} \\
& \qquad
 +\cE_1\sqrt{L_{n+1,J}}S_2\sqrt{R_{n+1,J}}^+
  (\tilde\pi_{n,J,\D}^{[n]}J\pi_{n,J,\D}-\tilde\rho_{n,J,\D}^{[n]}J\rho_{n,J,\D}) \\
& \qquad
 +\tau_{n,J,\D}\pi_{n,J,\D}-\sigma_{n,J,\D}\rho_{n,J,\D} \Big]X_1^{-1}
\displaybreak[0] \\
&= 
  \cE_{n+1}X_2^{-1}\Big(\sqrt{L_{n+1,J}}S_1\sqrt{R_{n+1,J}}-\sqrt{L_{n+1,J}}S_2\sqrt{R_{n+1,J}}\Big)
  X_1^{-1}. 
\end{align*}
Consequently, (i) is fulfilled. Hence part (a) is verified.
Now suppose that $f_1=f_2$ or $\tilde f_1=\tilde f_2$ is satisfied. In view of (a) then (i) holds.
Thus, using the identities $L_{n+1,J}L_{n+1,J}^+=\sqrt{L_{n+1,J}}^+\sqrt{L_{n+1,J}}$
and $R_{n+1,J}^+R_{n+1,J}=\sqrt{R_{n+1,J}}\sqrt{R_{n+1,J}}^+$ we obtain (b).
\qed\end{pf}

Now we are going to prove an inverse statement to Theorem \ref{DSS-T1}, i.e., we will show that
any function $f$ belonging to $\cP_{J,0}[\D,(A_j)_{j=0}^n]$ can be represented via
(\ref{DSS-R1-1}) and (\ref{DSS-R1-2}) with some $S\in\cS_{m\times m}(\D)$.

\begin{thm}\label{DSS-T2}
Let $J$ be an $m\times m$ signature matrix, let $n\in\N_0$, 
and let $(A_j)_{j=0}^n$ be a $J$-Potapov sequence. If $n\ge1$, then let $V_n\in\cY_{n,J}$ 
and $W_n\in\cZ_{n,J}$. Furthermore, let $\pi_{n,J}$, $\rho_{n,J}$, $\sigma_{n,J}$, and $\tau_{n,J}$
be given by (\ref{LRQ-T1-1}), (\ref{LRQ-T1-2}), (\ref{LRQ-T2-1}), and (\ref{LRQ-T2-2}). 
If $f\in\cP_{J,0}\left[\D,(A_j)_{j=0}^n\right]$, then there exists an $S\in\cS_{m\times m}(\D)$
such that $f$ admits the representations
\begin{align}\label{DSS-T2-1}
f&=\big(\cE_1J\tilde\tau_{n,J,\D}^{[n]}\sqrt{L_{n+1,J}}^+S\sqrt{R_{n+1,J}}+\pi_{n,J,\D}\big)\nonumber\\
&\qquad\cdot    \big(\cE_1J\tilde\sigma_{n,J,\D}^{[n]}\sqrt{L_{n+1,J}}^+S\sqrt{R_{n+1,J}}+\rho_{n,J,\D}\big)^{-1}
\end{align}
and (\ref{DSS-R1-2}).
\end{thm}

\begin{pf}
In view of $f\in\cP_{J,0}\left[\D,(A_j)_{j=0}^n\right]$, let $(A_j)_{j=n+1}^\infty$ be the sequence
of complex $m\times m$ matrices such that (\ref{NrTS}) is satisfied for each $w$ belonging to 
some neighborhood of $0$. According to Theorem \ref{T62}, for every $k\in\N_0$, the sequence
$(A_j)_{j=0}^k$ is a $J$-Potapov sequence. For each $k\in\N_0$, let $f_{c,k}$ be the
$J$-central $J$-Potapov function corresponding to $(A_j)_{j=0}^k$.
Then Theorem \ref{LRQ-T1} yields that $f_{c,n}$ admits the representation
$f_{c,n}=\pi_{n,J,\D}\rho_{n,J,\D}^{-1}$. For each $s\in\N_0$, let the matrix polynomials
$\pi_{n+s+1,J}$, $\rho_{n+s+1,J}$, $\sigma_{n+s+1,J}$, and $\tau_{n+s+1,J}$
be recursively defined by (\ref{CPF-P1-1}), (\ref{CPF-P1-2}), (\ref{CPF-P1-3}), and (\ref{CPF-P1-4})
where $t_{n+s+1}:=L_{n+s+1,J}^+(A_{n+s+1}-M_{n+s+1,J})$ and $u_{n+s+1}:=(A_{n+s+1}-M_{n+s+1,J})R_{n+s+1,J}^+$. 
From (\ref{LRQ-T1-2}), (\ref{LRQ-T2-2}),
(\ref{CPF-P1-2}), and (\ref{CPF-P1-4}) we see that, for every $j\in\N_0$,
the identities $\rho_{n+j,J}(0)=I_m$ and $\tau_{n+j,J}(0)=I_m$ hold.
Thus, for each $j\in\N_0$ and each $S\in\cS_{m\times m}(\D)$, the function
$$ \det\big(\cE_1J\tilde\sigma_{n+j,J,\D}^{[n+j]}\sqrt{L_{n+j+1,J}}^+S\sqrt{R_{n+j+1,J}}
   +\rho_{n+j,J,\D}\big) $$
does not vanish identically in $\D$. We are now going to prove that, for each $k\in\N_0$
and each $j\in\N_0$, there is an $S_{jk}\in\cS_{m\times m}(\D)$ such that
\begin{equation}\label{DSS-T2-3}
f_{c,n+j+k}=(\cE_1J\tilde\tau_{n+j,J,\D}^{[n+j]}F_{jk}+\pi_{n+j,J,\D})
    (\cE_1J\tilde\sigma_{n+j,J,\D}^{[n+j]}F_{jk}+\rho_{n+j,J,\D}\big)^{-1}
\end{equation}
where $F_{jk}:=\sqrt{L_{n+j+1,J}}^+S_{jk}\sqrt{R_{n+j+1,J}}$.
In the case $k=0$ we choose the constant $m\times m$ Schur function $S_{j0}$ defined on $\D$ with value 
$0_{m\times m}$ for all $j\in\N_0$. Then Theorem \ref{LRQ-T1} and Proposition \ref{CPF-P1} yield 
(\ref{DSS-T2-3}). Hence, there exists a $\kappa\in\N_0$ such that, for each $k\in\N_{0,\kappa}$,
there is a sequence $(S_{\ell k})_{\ell=0}^\infty$ of $m\times m$ Schur functions defined on $\D$
satisfying
\begin{equation}\label{DSS-T2-4}
f_{c,n+\ell+k}=(\cE_1J\tilde\tau_{n+\ell,J,\D}^{[n+\ell]}F_{\ell k}+\pi_{n+\ell,J,\D})
    (\cE_1J\tilde\sigma_{n+\ell,J,\D}^{[n+\ell]}F_{\ell k}+\rho_{n+\ell,J,\D})^{-1}
\end{equation}
for all $\ell\in\N_0$ where $F_{\ell k}:=\sqrt{L_{n+\ell+1,J}}^+S_{\ell k}\sqrt{R_{n+\ell+1,J}}$.
Let $j\in\N_0$. According to \cite[Proposition 4.1]{FKR2}, the matrix 
\begin{equation}\label{DSS-T2-4A}
K_{n+j+1,J}:=\sqrt{L_{n+j+1,J}}^+(A_{n+j+1}-M_{n+j+1,J})\sqrt{R_{n+j+1,J}}^+
\end{equation}
is contractive and
fulfills $A_{n+j+1}-M_{n+j+1,J}$ $=\sqrt{L_{n+j+1,J}}K_{n+j+1,J}\sqrt{R_{n+j+1,J}}$.
Consequently,
\begin{equation}\label{DSS-T2-5}
\sqrt{L_{n+j+1,J}}^+K_{n+j+1,J}\sqrt{R_{n+j+1,J}}=L_{n+j+1,J}^+(A_{n+j+1}-M_{n+j+1,J})
  =t_{n+j+1}
\end{equation}
and
\begin{equation}\label{DSS-T2-6}
\sqrt{L_{n+j+1,J}}K_{n+j+1,J}\sqrt{R_{n+j+1,J}}^+=(A_{n+j+1}-M_{n+j+1,J})R_{n+j+1,J}^+
  =u_{n+j+1}
\end{equation}
hold. Furthermore, the matrix-valued functions
\begin{equation}\label{DSS-T2-7}
\Theta_{j,\kappa+1}
 :=\sqrt{L_{n+j+1,J}}\sqrt{L_{n+j+2,J}}^+S_{j+1,\kappa}\sqrt{R_{n+j+2,J}}\sqrt{R_{n+j+1,J}}^+,
\end{equation}
\begin{equation}\label{DSS-T2-8}
\Phi_{j,\kappa+1}:=K_{n+j+1,J}+\cE_1\Theta_{j,\kappa+1}, \quad\mbox{and}\quad
\Psi_{j,\kappa+1}:=I_m+\cE_1K_{n+j+1,J}^*\Theta_{j,\kappa+1}
\end{equation}
are holomorphic in $\D$. 
Obviously, $\det\Psi_{j,\kappa+1}$ does not vanish identically in $\D$.
Thus, $S_{j,\kappa+1}:=\Phi_{j,\kappa+1}\Psi_{j,\kappa+1}^{-1}$ is a well-defined matrix-valued
function which is meromorphic in $\D$. 
Further,
\begin{align}\label{DSS-T2-9}
&  [\Psi_{j,\kappa+1}(z)]^*\Psi_{j,\kappa+1}(z)-[\Phi_{j,\kappa+1}(z)]^*\Phi_{j,\kappa+1}(z) \nonumber\\
& = I_m-K_{n+j+1,J}^*K_{n+j+1,J}-|z|^2[\Theta_{j,\kappa+1}(z)]^*(I_m-K_{n+j+1,J}K_{n+j+1,J}^*)
   \Theta_{j,\kappa+1}(z)
\end{align}
holds for each $z\in\D$.
According to \cite[Proposition 4.1]{FKR2}, we have 
\begin{equation}\label{DSS-T2-10}
\sqrt{L_{n+j+1,J}}(I_m-K_{n+j+1,J}K_{n+j+1,J}^*)\sqrt{L_{n+j+1,J}}=L_{n+j+2,J}
\end{equation}
and
\begin{equation}\label{DSS-T2-11}
\sqrt{R_{n+j+1,J}}(I_m-K_{n+j+1,J}^*K_{n+j+1,J})\sqrt{R_{n+j+1,J}}=R_{n+j+2,J}.
\end{equation}
In view of (\ref{DSS-T2-7}) and $\sqrt{L_{n+j+2,J}}^+L_{n+j+2,J}\sqrt{L_{n+j+2,J}}^+=L_{n+j+2,J}L_{n+j+2,J}^+$, equation (\ref{DSS-T2-10}) implies 
\begin{align}\label{DSS-T2-12}
& [\Theta_{j,\kappa+1}(z)]^*(I_m-K_{n+j+1,J}K_{n+j+1,J}^*)\Theta_{j,\kappa+1}(z) \nonumber\\
&=\! \sqrt{\!R_{n+j+1,J}}^+\!\!\sqrt{\!R_{n+j+2,J}}[S_{j+1,\kappa}(z)]^*L_{n+j+2,J}\!
      L_{n+j+2,J}^+S_{j+1,\kappa}(z)\sqrt{\!R_{n+j+2,J}}\sqrt{\!R_{n+j+1,J}}^+
\end{align}
for each $z\in\D$. Moreover, from (\ref{DSS-T2-4A}) and (\ref{DSS-T2-11}) we get
\begin{align}\label{DSS-T2-13}
& I_m-K_{n+j+1,J}^*K_{n+j+1,J} \nonumber\\
&= I_m-\sqrt{R_{n+j+1,J}}^+\sqrt{R_{n+j+1,J}} \nonumber\\
&\quad  +\sqrt{R_{n+j+1,J}}^+\sqrt{R_{n+j+1,J}}
 (I_m-K_{n+j+1,J}^*K_{n+j+1,J})\sqrt{R_{n+j+1,J}}\sqrt{R_{n+j+1,J}}^+ \nonumber\\
& = I_m-\sqrt{R_{n+j+1,J}}^+\sqrt{R_{n+j+1,J}}+\sqrt{R_{n+j+1,J}}^+R_{n+j+2,J}\sqrt{R_{n+j+1,J}}^+.
\end{align}
Since $S_{j+1,\kappa}$ belongs to $\cS_{m\times m}(\D)$ and
$L_{n+j+2,J}L_{n+j+2,J}^+$ is an orthoprojection matrix, we see that the matrix
$$ X(z):=I_m-|z|^2[S_{j+1,\kappa}(z)]^*L_{n+j+2,J}L_{n+j+2,J}^+S_{j+1,\kappa}(z) $$
is nonnegative Hermitian for each $z\in\D$ .
Taking into account (\ref{DSS-T2-9}), (\ref{DSS-T2-12}), and (\ref{DSS-T2-13}) we obtain
\begin{align}\label{DSS-T2-14}
& [\Psi_{j,\kappa+1}(z)]^*\Psi_{j,\kappa+1}(z)-[\Phi_{j,\kappa+1}(z)]^*\Phi_{j,\kappa+1}(z) \nonumber\\
&= \!I_m\!-\!\sqrt{\!R_{n+j+1,J}}^+\!\sqrt{\!R_{n+j+1,J}}\!+\!\sqrt{\!R_{n+j+1,J}}^+\!\sqrt{\!R_{n+j+2,J}}
  X\!(z)\sqrt{\!R_{n+j+2,J}}\sqrt{\!R_{n+j+1,J}}^+
\end{align}
for each $z\in\D$. The right-hand side of (\ref{DSS-T2-14}) is nonnegative Hermitian for every
$z\in\D$. Let $\cA$ be the set of all zeros of $\det\Psi_{j,\kappa+1}$.
Then
\begin{align}\label{DSS-T2-15}
& I_m-[S_{j,\kappa+1}(z)]^*S_{j,\kappa+1}(z) \nonumber\\
& = [\Psi_{j,\kappa+1}(z)]^{-*}
 \Big([\Psi_{j,\kappa+1}(z)]^*\Psi_{j,\kappa+1}(z)-[\Phi_{j,\kappa+1}(z)]^*\Phi_{j,\kappa+1}(z)\Big)
 [\Psi_{j,\kappa+1}(z)]^{-1} 
\end{align}
holds for all $z\in\D\setminus\cA$. Consequently, $S_{j,\kappa+1}$ is both holomorphic and contractive
in $\D\setminus\cA$. Since $\cA$ is a discrete subset of $\D$, Riemann's theorem on removable 
singularities of bounded holomorphic functions yields $S_{j,\kappa+1}\in\cS_{m\times m}(\D)$.
Furthermore, from (\ref{DSS-T2-10}) we obtain in particular \linebreak[0]
$\cR\left(\sqrt{L_{n+j+2,J}}^+\right)\linebreak[0]
 = \linebreak[0]
 \cR \linebreak[0]
 (L_{n+j+2,J})\linebreak[0]
 \sq\cR\left(\sqrt{L_{n+j+1,J}}\right)$
and, consequently, \linebreak[0]
$\sqrt{L_{n+j+1,J}}^+\sqrt{L_{n+j+1,J}}\sqrt{L_{n+j+2,J}}^+ \linebreak[0] 
 =\linebreak[0]
 \sqrt{L_{n+j+2,J}}^+$.
Similarly, (\ref{DSS-T2-11}) implies \linebreak[0]
$\sqrt{R_{n+j+2,J}}\sqrt{R_{n+j+1,J}}^+\sqrt{R_{n+j+1,J}} \linebreak[0]
 =\linebreak[0]
 \sqrt{R_{n+j+2,J}}$.
Hence, by using the setting \linebreak[0]
$F_{j+1,\kappa} \linebreak[0]
 := \linebreak[0]
 \sqrt{L_{n+j+2,J}}^+S_{j+1,\kappa}\sqrt{R_{n+j+2,J}}$ 
we get
\begin{equation}\label{DSS-T2-17}
\sqrt{L_{n+j+1,J}}^+\sqrt{L_{n+j+1,J}}F_{j+1,\kappa}\sqrt{R_{n+j+1,J}}^+\sqrt{R_{n+j+1,J}}
 = F_{j+1,\kappa}.
\end{equation}
From (\ref{DSS-T2-8}) and (\ref{DSS-T2-7}) 
we see that
$$ \det\left(I_m+zK_{n+j+1,J}^*\sqrt{L_{n+j+1,J}}F_{j+1,\kappa}(z)\sqrt{R_{n+j+1,J}}^+\right)\ne0$$
holds for each $z\in\D\setminus\cA$.
Taking into account (\ref{DSS-T2-4A}) and (\ref{DSS-T2-17}), 
from \cite[Remark 3.5]{FKL} we get therefore
$$ \det\left(I_m+z\sqrt{R_{n+j+1,J}}^+K_{n+j+1,J}^*\sqrt{L_{n+j+1,J}}F_{j+1,\kappa}(z)\right)\ne0$$
and
\begin{align}\label{DSS-T2-16}
& \left(I_m+zK_{n+j+1,J}^*\sqrt{L_{n+j+1,J}}F_{j+1,\kappa}(z)\sqrt{R_{n+j+1,J}}^+\right)^{-1}\sqrt{R_{n+j+1,J}} \nonumber\\
&= \sqrt{R_{n+j+1,J}}
\left(I_m+z\sqrt{R_{n+j+1,J}}^+K_{n+j+1,J}^*\sqrt{L_{n+j+1,J}}F_{j+1,\kappa}(z)\right)^{-1}
\end{align}
for each $z\in\D\setminus\cA$.
Let $F_{j,\kappa+1}:=\sqrt{L_{n+j+1,J}}^+S_{j,\kappa+1}\sqrt{R_{n+j+1,J}}$.
Taking into account (\ref{DSS-T2-8}), (\ref{DSS-T2-7}), (\ref{DSS-T2-16}), (\ref{DSS-T2-17}), (\ref{DSS-T2-5}), and (\ref{DSS-T2-6}) we obtain for each $z\in\D\setminus\cA$ then
\begin{align*}
F_{j,\kappa+1}(z)
&=  \sqrt{L_{n+j+1,J}}^+\Big(K_{n+j+1,J}+z\sqrt{L_{n+j+1,J}}F_{j+1,\kappa}(z)\sqrt{R_{n+j+1,J}}^+\Big)
  \nonumber\\
&\qquad\cdot 
 \Big(I_m+zK_{n+j+1,J}^*\sqrt{L_{n+j+1,J}}F_{j+1,\kappa}(z)\sqrt{R_{n+j+1,J}}^+\Big)^{\!-1}
 \!\sqrt{R_{n+j+1,J}} \nonumber
\displaybreak[0]\\
&=  \sqrt{L_{n+j+1,J}}^+\!\Big(\!\!K_{n+j+1,J}\!+\!z\sqrt{L_{n+j+1,J}}F_{j+1,\kappa}(z)\sqrt{R_{n+j+1,J}}^+\!\Big)\!
  \sqrt{R_{n+j+1,J}}   \nonumber\\
&\qquad\cdot 
\Big(I_m+z\sqrt{R_{n+j+1,J}}^+K_{n+j+1,J}^*\sqrt{L_{n+j+1,J}}F_{j+1,\kappa}(z)\Big)^{-1} \nonumber
\displaybreak[0]\\
&=  (t_{n+j+1}+zF_{j+1,\kappa}(z))(I_m+zu_{n+j+1}^*F_{j+1,\kappa}(z))^{-1}.
\end{align*}
Since $\cA$ is a discrete subset of $\D$, this implies
\begin{equation}\label{DSS-T2-18}
F_{j,\kappa+1}=(t_{n+j+1}+\cE_1F_{j+1,\kappa})(I_m+\cE_1u_{n+j+1}^*F_{j+1,\kappa})^{-1}.
\end{equation}
Moreover, from (\ref{CPF-P1-4}) and (\ref{CPF-P1-3}) the equations
$$ \tilde\tau_{n+j+1,J,\D}^{[n+j+1]}=\cE_1\tilde\tau_{n+j,J,\D}^{[n+j]}+J\pi_{n+j,J,\D}u_{n+j+1}^* $$
and
$$ \tilde\sigma_{n+j+1,J,\D}^{[n+j+1]}=\cE_1\tilde\sigma_{n+j,J,\D}^{[n+j]}+J\rho_{n+j,J,\D}u_{n+j+1}^*$$
follow. 
Hence, by succinct setting $Y:=I_m+\cE_1u_{n+j+1}^*F_{j+1,\kappa}$ from (\ref{DSS-T2-18}),
(\ref{CPF-P1-1}), and (\ref{CPF-P1-2}) we obtain
\begin{align*}
&\cE_1J\tilde\tau_{n+j,J,\D}^{[n+j]}F_{j,\kappa+1}+\pi_{n+j,J,\D}\\
&=[\cE_1J\tilde\tau_{n+j,J,\D}^{[n+j]}(t_{n+j+1}+\cE_1F_{j+1,\kappa})+\pi_{n+j,J,\D}(I_m+\cE_1u_{n+j+1}^*F_{j+1,\kappa})]Y^{-1} \displaybreak[0]\\
&=(\cE_1J\tilde\tau_{n+j+1,J,\D}^{[n+j+1]}F_{j+1,\kappa}+\pi_{n+j+1,J,\D})Y^{-1}
\end{align*}
and
\begin{align*}
&\cE_1J\tilde\sigma_{n+j,J,\D}^{[n+j]}F_{j,\kappa+1}+\rho_{n+j,J,\D}\\
&=[\cE_1J\tilde\sigma_{n+j,J,\D}^{[n+j]}(t_{n+j+1}+\cE_1F_{j+1,\kappa})+\rho_{n+j,J,\D}(I_m+\cE_1u_{n+j+1}^*F_{j+1,\kappa})]Y^{-1}
\displaybreak[0] \\
&= (\cE_1J\tilde\sigma_{n+j+1,J,\D}^{[n+j+1]}F_{j+1,\kappa}+\rho_{n+j+1,J,\D})Y^{-1}.
\end{align*}
Consequently, an application of (\ref{DSS-T2-4}) yields
\begin{align*}
&f_{c,n+j+(\kappa+1)}=f_{c,n+(j+1)+\kappa} \\
&=(\cE_1J\tilde\tau_{n+j+1,J,\D}^{[n+j+1]}F_{j+1,k}+\pi_{n+j+1,J,\D})
    (\cE_1J\tilde\sigma_{n+j+1,J,\D}^{[n+j+1]}F_{j+1,k}+\rho_{n+j+1,J,\D})^{-1}
\displaybreak[0] \\
&=(\cE_1J\tilde\tau_{n+j,J,\D}^{[n+j]}F_{j,\kappa+1}+\pi_{n+j,J,\D})
     (\cE_1J\tilde\sigma_{n+j,J,\D}^{[n+j]}F_{j,\kappa+1}+\rho_{n+j,J,\D})^{-1}.
\end{align*}
Thus, we have shown by induction that for all $j,k\in\N_0$ there is an $S_{jk}\in\cS_{m\times m}(\D)$
such that (\ref{DSS-T2-3}) is satisfied where $F_{jk}:=\sqrt{L_{n+j+1,J}}^+S_{jk}\sqrt{R_{n+j+1,J}}$.
The matricial version of Montel's theorem yields that there are an $S\!\in\!\cS_{m\times m}\!(\D)$
and a subsequence $(S_{0k_r})_{r=0}^\infty$ of $(S_{0k})_{k=0}^\infty$ such that
\begin{equation}\label{DSS-T2-19}
\lim\limits_{r\to\infty}S_{0k_r}(w)=S(w)
\end{equation}
holds for each $w\in\D$. 
Denote by $\cA_0$ the set of all zeros of the function 
$\det(\cE_1J\tilde\sigma_{n,J,\D}^{[n]}\sqrt{L_{n+1,J}}^+S\sqrt{R_{n+1,J}}+\rho_{n,J,\D})$. Then 
(\ref{DSS-T2-19}) and a continuity argument yield that, for each $w\in\mathbb H_f\setminus\cA_0$,
there is an $r_0\in\N$ such that
$$\det\big(wJ\tilde\sigma_{n,J}^{[n]}(w)\sqrt{L_{n+1,J}}^+S_{0k_r}(w)\sqrt{R_{n+1,J}}+\rho_{n,J}(w)\big)\ne0 $$
holds for each $r\in\N_{r_0,\infty}$. Thus, using (\ref{DSS-T2-3}), (\ref{DSS-T2-19}), 
and Corollary \ref{CPF-C0} we obtain 
\begin{align*}
f(w)&=\big(wJ\tilde\tau_{n,J}^{[n]}(w)\sqrt{L_{n+1,J}}^+S(w)\sqrt{R_{n+1,J}}+\pi_{n,J}(w)\big) \\
&\qquad\cdot    \big(wJ\tilde\sigma_{n,J}^{[n]}(w)\sqrt{L_{n+1,J}}^+S(w)\sqrt{R_{n+1,J}}+\rho_{n,J}(w)\big)^{-1}
\end{align*}
for each $w\in\mathbb H_f\setminus\cA_0$ and hence (\ref{DSS-T2-1}).
Finally, in view of (\ref{DSS-T2-1}) and Remark \ref{DSS-R1}, $f$ can be represented via
(\ref{DSS-R1-2}).
\qed\end{pf}

\begin{rem}\label{DSS-R3}
Let $J$ be an $m\times m$ signature matrix, let $n\in\N_0$, 
and let $(A_j)_{j=0}^n$ be a $J$-Potapov sequence. Denote by $\cS_{m\times m}^\sharp(\D)$ the set of
all $S\in\cS_{m\times m}(\D)$ satisfying $L_{n+1,J}L_{n+1,J}^+SR_{n+1,J}^+R_{n+1,J}=S$.
Then a combination of Theorem \ref{DSS-T1}, Remark \ref{DSS-R2}, part (b) of Proposition \ref{DSS-P1},
and Theorem \ref{DSS-T2} shows that there is a bijective correspondence between the sets 
$\cS_{m\times m}^\sharp(\D)$ and $\cP_{J,0}[\D,(A_j)_{j=0}^n]$.
\end{rem}

Theorems \ref{DSS-T1} and \ref{DSS-T2} give a complete description of the solution set of problem (P) via some kind of linear fractional transformation which is given by (\ref{DSS-R1-1}) (respectively, (\ref{DSS-R1-2})).
A closer look at these two theorems 
shows that there is still some freedom in the construction of the matrix polynomials
from which the linear fractional transformation is built.
This freedom consists in the special choice of $V_n\in\cY_{n,J}$ and $W_n\in\cZ_{n,J}$.
Now the question arises whether we can find suitable matrix polynomials 
$\pi_{n,J}$, $\rho_{n,J}$, $\sigma_{n,J}$, and $\tau_{n,J}$ such that, 
for each $S\in\cS_{m\times m}(\D)$, the zeros of the determinant of the 
"denominator function" in the representation (\ref{DSS-R1-1}) (respectively, (\ref{DSS-R1-2}))
of the $J$-Potapov function $f$ are exactly the poles of $f$.
The following considerations give some answer to this question.

Henceforth we will use the following notation.
Let $J$ be an $m\times m$ signature matrix, let $n\in\N$, 
and let $(A_j)_{j=0}^n$ be a $J$-Potapov sequence. 
Then we will write $\tilde\cY_{n,J}$ to denote the set of all $V_n\in\cY_{n,J}$ for which 
the matrix polynomials $\pi_{n,J}$ and $\rho_{n,J}$ defined by (\ref{LRQ-T1-1}) and (\ref{LRQ-T1-2})
satisfy the condition 
\begin{equation}\label{DSS-E1}
 \cN(\rho_{n,J}(w))\cap\cN(\pi_{n,J}(w))=\{0_{m\times1}\} 
\end{equation}
for all $w\in\D$. Similarly, $\tilde\cZ_{n,J}$ denotes the set of all $W_n\in\cZ_{n,J}$ for which 
the matrix polynomials $\sigma_{n,J}$ and $\tau_{n,J}$ defined by (\ref{LRQ-T2-1}) and (\ref{LRQ-T2-2})
fulfill
\begin{equation}\label{DSS-E2}
 \cN([\tau_{n,J}(w)]^*)\cap\cN([\sigma_{n,J}(w)]^*)=\{0_{m\times1}\} 
\end{equation}
for all $w\in\D$. 

\begin{rem}\label{DSS-R4}
Let $J$ be an $m\times m$ signature matrix, let $n\in\N$, 
and let $(A_j)_{j=0}^n$ be a $J$-Potapov sequence. 
Let $V_n^\Box$ and $W_n^\Box$ be given by (\ref{NrVnWn}). Then according to 
\cite[Lemmata 2.16 and 2.17]{FKRS} we have $V_n^\Box\in\tilde\cY_{n,J}$ and $W_n^\Box\in\tilde\cZ_{n,J}$.
\end{rem}

\begin{prop}\label{DSS-P2}
Let $J$ be an $m\times m$ signature matrix, let $n\in\N_0$, 
and let $(A_j)_{j=0}^n$ be a $J$-Potapov sequence. If $n\ge1$, then let $V_n\in\tilde\cY_{n,J}$ and
$W_n\in\tilde\cZ_{n,J}$. Let the matrix polynomials
$\pi_{n,J}$, $\rho_{n,J}$, $\sigma_{n,J}$, and $\tau_{n,J}$
be given by (\ref{LRQ-T1-1}), (\ref{LRQ-T1-2}), (\ref{LRQ-T2-1}), and (\ref{LRQ-T2-2}). 
Let $S\in\cS_{m\times m}(\D)$, and let the matrix-valued function $f$ be given by (\ref{DSS-R1-1}).
Then
\begin{equation}\label{DSS-P2-1}
\mathbb H_f=\big\{w\in\D:\det
\big(wJ\tilde\sigma_{n,J}^{[n]}(w)\sqrt{L_{n+1,J}}^+S(w)\sqrt{R_{n+1,J}}+\rho_{n,J}(w)\big)\ne0\big\} 
\end{equation}
and 
\begin{equation}\label{DSS-P2-2}
\mathbb H_f=\big\{w\in\D:\det
\big(w\sqrt{L_{n+1,J}}S(w)\sqrt{R_{n+1,J}}^+\tilde\pi_{n,J}^{[n]}(w)J+\tau_{n,J}(w)\big)\ne0\big\}
\end{equation}
hold true, and $f$ admits for each $w\in\mathbb H_f$ the representations
\begin{align}\label{DSS-P2-3}
f(w)&=\big(wJ\tilde\tau_{n,J}^{[n]}(w)\sqrt{L_{n+1,J}}^+S(w)\sqrt{R_{n+1,J}}+\pi_{n,J}(w)\big) \nonumber\\
&\qquad\cdot    \big(wJ\tilde\sigma_{n,J}^{[n]}(w)\sqrt{L_{n+1,J}}^+S(w)\sqrt{R_{n+1,J}}+\rho_{n,J}(w)\big)^{-1}
\end{align}
and
\begin{align}\label{DSS-P2-4}
f(w)&=\big(w\sqrt{L_{n+1,J}}S(w)\sqrt{R_{n+1,J}}^+\tilde\pi_{n,J}^{[n]}(w)J+\tau_{n,J}(w)\big)^{-1}
\nonumber\\
&\qquad\cdot    \big(w\sqrt{L_{n+1,J}}S(w)\sqrt{R_{n+1,J}}^+\tilde\rho_{n,J}^{[n]}(w)J+\sigma_{n,J}(w)\big).
\end{align}
\end{prop}

\begin{pf}
In view of (\ref{DSS-R1-1}) the inclusion
$$\big\{w\in\D:\det
\big(wJ\tilde\sigma_{n,J}^{[n]}(w)\sqrt{L_{n+1,J}}^+S(w)\sqrt{R_{n+1,J}}+\rho_{n,J}(w)\big)\ne0\big\} 
\sq\mathbb H_f$$
is obvious.
Now let $w_0\in\mathbb H_f$. 
From (\ref{DSS-R1-1}) then
\begin{align}\label{DSS-P2-5}
& w_0J\tilde\tau_{n,J}^{[n]}(w_0)\sqrt{L_{n+1,J}}^+S(w_0)\sqrt{R_{n+1,J}}+\pi_{n,J}(w_0) 
 \nonumber\\
&= f(w_0) \Big(w_0J\tilde\sigma_{n,J}^{[n]}(w_0)\sqrt{L_{n+1,J}}^+S(w_0)\sqrt{R_{n+1,J}}
             +\rho_{n,J}(w_0)\Big)
\end{align}
follows.
Further, $K:=S(w_0)$ is a contractive $m\times m$ matrix. Let 
$A_{n+1}:=M_{n+1,J}+\sqrt{L_{n+1,J}}K\sqrt{R_{n+1,J}}$. According to Theorem \ref{T}, $(A_j)_{j=0}^{n+1}$
is a $J$-Potapov sequence. Let $t_{n+1}$ be given by (\ref{RF-R1-1}),
and for each $w\in\C$ let the matrix polynomials $\pi_{n+1,J}$ and $\rho_{n+1,J}$
be defined by (\ref{CPF-P1-1}) and (\ref{CPF-P1-2}) (with $s=0$).
Because of
$t_{n+1}=L_{n+1,J}^+\sqrt{L_{n+1,J}}K\sqrt{R_{n+1,J}}=$ $\sqrt{L_{n+1,J}}^+S(w_0)\sqrt{R_{n+1,J}}$
from (\ref{DSS-P2-5}) we get $\pi_{n+1,J}(w_0)=f(w_0)\rho_{n+1,J}(w_0)$
and consequently 
\begin{equation}\label{DSS-P2-6}
\cN(\rho_{n+1,J}(w_0))\sq\cN(\pi_{n+1,J}(w_0)).
\end{equation}
In the case $n\ge1$, because of $V_n\in\tilde\cY_{n,J}$ condition (\ref{DSS-E1}) is satisfied
for each $w\in\D$. If $n=0$ then (\ref{DSS-E1}) obviously holds true for all $w\in\D$. 
According to \cite[Proposition 3.6]{FKRS},
this implies 
$\cN(\rho_{n+1,J}(w))\cap\cN(\pi_{n+1,J}(w))=\{0_{m\times1}\}$
for each $w\in\D$. 
In view of (\ref{DSS-P2-6}) we get therefore
$$\det
\Big(\!w_0J\tilde\sigma_{n,J}^{[n]}(w_0)\sqrt{L_{n+1,J}}^+S(w_0)\sqrt{R_{n+1,J}}+\rho_{n,J}(w_0)\!\Big)
=\det(\rho_{n+1,J}(w_0))\ne0. $$
Hence (\ref{DSS-P2-1}) is verified. Having in mind Remark \ref{DSS-R1} we see that $f$ admits the 
representation (\ref{DSS-R1-2}). Thus, equation (\ref{DSS-P2-2}) can be checked analogously to 
(\ref{DSS-P2-1}).
Furthermore, from (\ref{DSS-R1-1}), (\ref{DSS-R1-2}), (\ref{DSS-P2-1}), 
and (\ref{DSS-P2-2}) it follows immediately that $f$ can 
be represented via (\ref{DSS-P2-3}) and (\ref{DSS-P2-4}) for each $w\in\mathbb H_f$.
\qed\end{pf}

Now we are able to prove Theorem \ref{S1-T1}.\\[2ex]
{\bf Proof of Theorem \ref{S1-T1}.}
Use Theorems \ref{DSS-T1} and \ref{DSS-T2} as well as Remark \ref{DSS-R4} and Proposition \ref{DSS-P2}.
\qed

\vspace{2ex}
Theorems \ref{DSS-T1} and \ref{DSS-T2} as well as Proposition \ref{DSS-P2} should be regarded in 
connection to \cite[Section 3]{FKRS}.
As it was mentioned above, there is some freedom in the construction of the matrix polynomials
$\pi_{n,J}$, $\rho_{n,J}$, $\sigma_{n,J}$, and $\tau_{n,J}$ occuring in the representations
(\ref{DSS-R1-1}) and (\ref{DSS-R1-2}) of the functions belonging to $\cP_{J,0}[\D,(A_j)_{j=0}^n]$.
Section 3 in \cite{FKRS} points out a possibility for a recursive construction of such matrix polynomials that fulfill the requirements of Theorems \ref{DSS-T1} and \ref{DSS-T2} (and even of Proposition \ref{DSS-P2}).
More precisely, we can state the following.

\begin{rem}\label{DSS-R5}
Let $J$ be an $m\times m$ signature matrix, let $n\in\N_0$, 
and let $(A_j)_{j=0}^n$ be a $J$-Potapov sequence. Further, let
$$\fpi_{0,J}(w):=A_0,\quad \frho_{0,J}:=I_m,\quad \fsigma_{0,J}(w):=A_0,\quad\mbox{and}\quad 
\ftau_{0,J}(w):=I_m $$
for all $w\in\C$. If $n\ge1$, then let, for each $k\in\N_{0,n-1}$ and each $w\in\C$, the matrix polynomials 
$\fpi_{k+1,J}$, $\frho_{k+1,J}$, $\fsigma_{k+1,J}$, and $\ftau_{k+1,J}$ 
be defined recursively by
$$\fpi_{k+1,J}(w)\!:=\fpi_{k,J}(w)+wJ\tilde\ftau_{k,J}^{[k]}(w)t_{k+1},
\ \ 
\frho_{k+1,J}(w)\!:=\frho_{k,J}(w)+wJ\tilde\fsigma_{k,J}^{[k]}(w)t_{k+1},$$
$$\fsigma_{k+1,J}(w)\!:=\fsigma_{k,J}(w)+u_{k+1}w\tilde\frho_{k,J}^{[k]}(w)J,
\ \ 
\ftau_{k+1,J}(w)\!:=\ftau_{k,J}(w)+u_{k+1}w\tilde\fpi_{k,J}^{[k]}(w)J$$
where $t_{k+1}$ and $u_{k+1}$ are given by (\ref{RF-R1-1}). Then, in the case $n\ge1$, 
there exist some $\mathbf V_n\in\tilde\cY_{n,J}$ and $\mathbf W_n\in\tilde\cZ_{n,J}$
such that (\ref{DSS-R1-3}) and (\ref{DSS-R1-4}) hold for each $w\in\C$
(see \cite[Remark 4.2 and Proposition 3.6]{FKRS}).
However, observe that in general 
the matrix polynomials $\fpi_{n,J}$, $\frho_{n,J}$, $\fsigma_{n,J}$, and $\ftau_{n,J}$ 
constructed in this way do not coincide with $\pi_{n,J}$, $\rho_{n,J}$, $\sigma_{n,J}$, and $\tau_{n,J}$, respectively, given in Theorem \ref{S1-T1}.
\end{rem}

\begin{thm}\label{DSS-T3}
Let $J$ be an $m\times m$ signature matrix, let $n\in\N_0$, 
and let $(A_j)_{j=0}^n$ be a $J$-Potapov sequence. Further, let 
the matrix polynomials $\fpi_{n,J}$, $\frho_{n,J}$, $\fsigma_{n,J}$, and $\ftau_{n,J}$
be recursively defined as in Remark \ref{DSS-R5}.
Then, for each $S\in\cS_{m\times m}(\D)$, 
\begin{align*}
f_S&:=\big(\cE_1J\tilde\ftau_{n,J,\D}^{[n]}\sqrt{L_{n+1,J}}^+S\sqrt{R_{n+1,J}}+\fpi_{n,J,\D}\big) \nonumber\\
&\qquad\cdot
\big(\cE_1J\tilde\fsigma_{n,J,\D}^{[n]}\sqrt{L_{n+1,J}}^+S\sqrt{R_{n+1,J}}+\frho_{n,J,\D}\big)^{-1} 
\end{align*}
is a (well-defined) matrix-valued function meromorphic in $\D$, and the set $\mathbb H_{f_S}$ of all $w\in\D$ at which $f_S$ is holomorphic fulfills
\begin{align}\label{DSS-T3-1}
\mathbb H_{f_S}&=\big\{w\in\D:\det
\big(wJ\tilde\fsigma_{n,J}^{[n]}(w)\sqrt{L_{n+1,J}}^+S(w)\sqrt{R_{n+1,J}}+\frho_{n,J}(w)\big)\ne0\big\}\nonumber\\
&=\big\{w\in\D:\det
\big(w\sqrt{L_{n+1,J}}S(w)\sqrt{R_{n+1,J}}^+\tilde\fpi_{n,J}^{[n]}(w)J+\ftau_{n,J}(w)\big)\ne0\big\}.
\end{align}
Further, for each $S\in\cS_{m\times m}(\D)$,
$f_S$ admits the representations
\begin{align}\label{DSS-T3-2}
f_S&=\big(\cE_1\sqrt{L_{n+1,J}}S\sqrt{R_{n+1,J}}^+\tilde\fpi_{n,J,\D}^{[n]}J+\ftau_{n,J,\D}\big)^{-1} \nonumber\\
&\qquad\cdot
   \big(\cE_1\sqrt{L_{n+1,J}}S\sqrt{R_{n+1,J}}^+\tilde\frho_{n,J,\D}^{[n]}J+\fsigma_{n,J,\D}\big)
\end{align}
and
\begin{align}\label{DSS-T3-3}
f_S&=\big(\cE_1J\tilde\tau_{n,J,\D}^{[n]}\sqrt{L_{n+1,J}}^+S\sqrt{R_{n+1,J}}+\pi_{n,J,\D}\big) \nonumber\\
&\qquad\cdot
\big(\cE_1J\tilde\sigma_{n,J,\D}^{[n]}\sqrt{L_{n+1,J}}^+S\sqrt{R_{n+1,J}}+\rho_{n,J,\D}\big)^{-1} 
\end{align}
where $\pi_{n,J}$, $\rho_{n,J}$, $\sigma_{n,J}$, and $\tau_{n,J}$
are the matrix polynomials given in Theorem \ref{S1-T1}.
Moreover, 
\begin{equation}\label{DSS-T3-4}
 \cP_{J,0}[\D,(A_j)_{j=0}^n]=\{f_S:S\in\cS_{m\times m}(\D)\} 
\end{equation}
holds true.
\end{thm}

\begin{pf}
In the case $n=0$ the assertion is an immediate consequence of Theorem \ref{S1-T1}. 
Now suppose $n\ge1$. 
Then Remark \ref{DSS-R5} shows that
there exist some $\mathbf V_n\in\tilde\cY_{n,J}$ and $\mathbf W_n\in\tilde\cZ_{n,J}$
such that (\ref{DSS-R1-3}) and (\ref{DSS-R1-4}) hold for each $w\in\C$. Thus, applying Theorems \ref{DSS-T1} and \ref{DSS-T2} and Proposition \ref{DSS-P2} we obtain 
(\ref{DSS-T3-1}), (\ref{DSS-T3-2}), and (\ref{DSS-T3-4}). 
Finally, Remarks \ref{LRQ-R2} and \ref{DSS-R1} yield (\ref{DSS-T3-3}). 
\qed\end{pf}

\section{The case of a unique solution}
\label{s4}

In this section we study the case of a given $J$-Potapov sequence $(A_j)_{j=0}^n$ for which the $J$-central $J$-Potapov function corresponding to $(A_j)_{j=0}^n$ is the unique 
$J$-Potapov function fulfilling (\ref{ICP}) for each $j\in\N_{0,n}$.
This section generalizes the corresponding results obtained in \cite[Section 7]{FKL} where the case of  matricial Schur functions is treated.

\begin{lem}\label{eind-L1}
Let $J$ be an $m\times m$ signature matrix, let $n\in\N_0$, 
and let $(A_j)_{j=0}^n$ be a $J$-Potapov sequence. Further, let $(A_j)_{j=0}^\infty$ be the $J$-central $J$-Potapov sequence corresponding to $(A_j)_{j=0}^n$. Then the following statements are equivalent:
\begin{enumerate}\item[(i)] There is a unique $f\in\cP_{J,0}(\D)$ such that 
(\ref{ICP}) is fulfilled for each $j\in\N_{0,n}$ (namely the $J$-central $J$-Potapov function $f=f_{c,n}$ corresponding to $(A_j)_{j=0}^n$).
\item[(ii)] For each $k\in\N_{n+1,\infty}$, the identities $L_{k,J}=0_{m\times m}$ and 
$R_{k,J}=0_{m\times m}$ are satisfied.
\item[(iii)] $L_{n+1,J}=0_{m\times m}$ or $R_{n+1,J}=0_{m\times m}$.\end{enumerate}
\end{lem}

\begin{pf}
The equivalence of (i) and (ii) follows easily from Theorem \ref{T72} and
\cite[Proposition 3.12]{FKR2}, whereas the equivalence of (ii) and (iii) is an immediate consequence of \cite[Proposition 3.12 and Proposition 5.5]{FKR2}.
\qed\end{pf}

\begin{prop}\label{eind-P1}
Let $J$ be an $m\times m$ signature matrix, let $n\in\N_0$, 
and let $(A_j)_{j=0}^n$ be a $J$-Potapov sequence. If $n\ge1$, then let $V_n\in\cY_{n,J}$ 
and $W_n\in\cZ_{n,J}$. Furthermore, let the matrix polynomials $\pi_{n,J}$, $\rho_{n,J}$, $\sigma_{n,J}$, and $\tau_{n,J}$
be given by (\ref{LRQ-T1-1}), (\ref{LRQ-T1-2}), (\ref{LRQ-T2-1}), and (\ref{LRQ-T2-2}). 
Then the following statements are equivalent:
\begin{enumerate}\item[(i)] There is a unique $f\in\cP_{J,0}(\D)$ such that 
(\ref{ICP}) is fulfilled for each $j\in\N_{0,n}$.
\item[(ii)] The identities $\tilde\rho_{n,J}^{[n]}J\rho_{n,J}=\tilde\pi_{n,J}^{[n]}J\pi_{n,J}$
and $\tau_{n,J}J\tilde\tau_{n,J}^{[n]}=\sigma_{n,J}J\tilde\sigma_{n,J}^{[n]}$ hold.
\item[(iii)] For each $z\in\T$, the identities $[\rho_{n,J}(z)]^*J\rho_{n,J}(z)=[\pi_{n,J}(z)]^*J\pi_{n,J}(z)$ and 
$\tau_{n,J}(z)J[\tau_{n,J}(z)]^*=\sigma_{n,J}(z)J[\sigma_{n,J}(z)]^*$ are satisfied.
\item[(iv)] There is some $z\in\T$ such that $[\rho_{n,J}(z)]^*J\rho_{n,J}(z)=[\pi_{n,J}(z)]^*J\pi_{n,J}(z)$ or 
$\tau_{n,J}(z)J[\tau_{n,J}(z)]^*=\sigma_{n,J}(z)J[\sigma_{n,J}(z)]^*$ is fulfilled.
\end{enumerate}
\end{prop}

\begin{pf}
The assertion is an immediate consequence of Lemma \ref{eind-L1} and \cite[Proposition 2.11]{FKRS}.
\qed\end{pf}

In the sequel, if $\beta$ is some matrix polynomial of quadratic size, then we will use the 
notation $\cN_\beta:=\{w\in\C:\det\beta(w)=0\}$.

\begin{rem}\label{eind-R1}
Let $J$ be an $m\times m$ signature matrix, let $n\in\N_0$, 
and let $(A_j)_{j=0}^n$ be a $J$-Potapov sequence. If $n\ge1$, then let $V_n\in\cY_{n,J}$ 
and $W_n\in\cZ_{n,J}$. Furthermore, let the matrix polynomials $\rho_{n,J}$ and $\tau_{n,J}$
be given by (\ref{LRQ-T1-2}) and (\ref{LRQ-T2-2}), respectively. 
Then one can easily see that each of the sets $\cN_{\rho_{n,J}}$ and $\cN_{\tau_{n,J}}$ 
consists of at most $nm$ elements. Moreover, taking into account \cite[Lemma 1.2.2]{DFK},
it follows that each of the sets $\cN_{\tilde\rho_{n,J}^{[n]}}$ and $\cN_{\tilde\tau_{n,J}^{[n]}}$
also consists of at most $nm$ elements. 
\end{rem}

In view of Theorem \ref{LRQ-T1} and Remark \ref{eind-R1} there exists a unique rational extension
$f_{c,n}^\diamondsuit$ of the $J$-central $J$-Potapov function $f_{c,n}$ to
$\C\setminus\{w_1,w_2,\ldots,w_r\}$ for some complex numbers $w_1,w_2,\ldots,w_r$ with $r\le nm$.

Recall that a matrix $A\in\C^{m\times m}$ is said to be $J$-unitary if $A^*JA=J$, where $J$ is some $m\times m$ signature matrix. Obviously, $A\in\C^{m\times m}$ is $J$-unitary if and only if $AJA^*=J$.

\begin{cor}\label{eind-C1}
Let $J$ be an $m\times m$ signature matrix, let $n\in\N_0$, 
and let $(A_j)_{j=0}^n$ be a $J$-Potapov sequence. Further, let $f_{c,n}$ be the $J$-central $J$-Potapov function corresponding to $(A_j)_{j=0}^n$ and let $f_{c,n}^\diamondsuit$ 
be the unique rational extension of $f_{c,n}$. 
Then the following statements are equivalent:
\begin{enumerate}\item[(i)] There is a unique $f\in\cP_{J,0}(\D)$ such that 
(\ref{ICP}) is fulfilled for each $j\in\N_{0,n}$.
\item[(ii)] There exists a finite subset $\mathfrak F$ of $\T$ such that 
for each $z\in\T\setminus\mathfrak F$ the matrix $f_{c,n}^\diamondsuit(z)$ is $J$-unitary.
\item[(iii)] There is some $z\in\T\setminus(\cN_{\rho_{n,J}}\cap\cN_{\tau_{n,J}})$ such that 
the matrix $f_{c,n}^\diamondsuit(z)$ is $J$-unitary, where $\rho_{n,J}$ and $\tau_{n,J}$
are given by (\ref{LRQ-T1-2}) and (\ref{LRQ-T2-2}) with some $V_n\in\cY_{n,J}$ 
and $W_n\in\cZ_{n,J}$ if $n\ge1$.
\end{enumerate}
\end{cor}

\begin{pf}
Use Proposition \ref{eind-P1}, Remark \ref{eind-R1}, and Theorem \ref{LRQ-T1}.
\qed\end{pf}

\begin{cor}\label{eind-C2}
Let $J$ be an $m\times m$ signature matrix, let $n\in\N_0$, 
and let $(A_j)_{j=0}^n$ be a $J$-Potapov sequence. Further, let $f_{c,n}$ be the $J$-central $J$-Potapov function corresponding to $(A_j)_{j=0}^n$ and let $f_{c,n}^\diamondsuit$ 
be the unique rational extension of $f_{c,n}$. 
If $n\ge1$, then let $V_n\in\cY_{n,J}$ 
and $W_n\in\cZ_{n,J}$. 
Let the matrix polynomials $\pi_{n,J}$, $\rho_{n,J}$, $\sigma_{n,J}$, and $\tau_{n,J}$
be given by (\ref{LRQ-T1-1}), (\ref{LRQ-T1-2}), (\ref{LRQ-T2-1}), and (\ref{LRQ-T2-2}).
Then the following statements are equivalent:
\begin{enumerate}\item[(i)] There is a unique $f\in\cP_{J,0}(\D)$ such that 
(\ref{ICP}) is fulfilled for each $j\in\N_{0,n}$.
\item[(ii)] The complex-valued functions $\det\tilde\pi_{n,J}^{[n]}$ and 
$\det\tilde\sigma_{n,J}^{[n]}$ do not vanish identically, and 
$f_{c,n}^\diamondsuit\!=\!J(\tilde\pi_{n,J}^{[n]})^{\!-1}\tilde\rho_{n,J}^{[n]}J \linebreak[0]
 =J\tilde\tau_{n,J}^{[n]}(\tilde\sigma_{n,J}^{[n]})^{-1}J$ holds.
\item[(iii)] Each of the sets $\cN_{\tilde\pi_{n,J}^{[n]}}$ and $\cN_{\tilde\sigma_{n,J}^{[n]}}$
consists of at most $nm$ elements, and for each $z\in\T\setminus\cN_{\tilde\pi_{n,J}^{[n]}}$ 
the matrix $(\tilde\pi_{n,J}^{[n]}(z))^{-1}\tilde\rho_{n,J}^{[n]}(z)$ is $J$-unitary,
and for each $z\in\T\setminus\cN_{\tilde\sigma_{n,J}^{[n]}}$ 
the matrix $\tilde\tau_{n,J}^{[n]}(z)(\tilde\sigma_{n,J}^{[n]}(z))^{-1}$ is $J$-unitary.
\item[(iv)] There exists some $z\in\T\setminus\cN_{\tilde\pi_{n,J}^{[n]}}$ such that
$(\tilde\pi_{n,J}^{[n]}(z))^{-1}\tilde\rho_{n,J}^{[n]}(z)$ is $J$-contrac\-tive,
or there exists some $z\in\T\setminus\cN_{\tilde\sigma_{n,J}^{[n]}}$ 
such that $\tilde\tau_{n,J}^{[n]}(z)(\tilde\sigma_{n,J}^{[n]}(z))^{-1}$ is $J$-contractive.
\end{enumerate}
\end{cor}

\begin{pf}
From Proposition \ref{eind-P1} and Remark \ref{eind-R1} we infer that, if (i) is fulfilled, then
each of the sets $\cN_{\tilde\pi_{n,J}^{[n]}}$ and $\cN_{\tilde\sigma_{n,J}^{[n]}}$
consists of at most $nm$ elements. Furthermore, Theorem \ref{LRQ-T1} implies 
$f_{c,n}^\diamondsuit=\pi_{n,J}\rho_{n,J}^{-1}=\tau_{n,J}^{-1}\sigma_{n,J}$. Hence the equivalence
of (i) and (ii) follows from Proposition \ref{eind-P1}.
Moreover, \cite[Proposition 2.11]{FKRS} yields
$$ [\rho_{n,J}(z)]^*J\rho_{n,J}(z)-[\pi_{n,J}(z)]^*J\pi_{n,J}(z)=R_{n+1,J} $$
and
$$ \tau_{n,J}(z)J[\tau_{n,J}(z)]^*-\sigma_{n,J}(z)J[\sigma_{n,J}(z)]^*=L_{n+1,J} $$
for each $z\in\T$. Thus, taking into account \cite[Lemma 1.2.2]{DFK}, we obtain
\begin{align}
& \Big([\tilde\pi_{n,J}^{[n]}(z)]^{-1}\tilde\rho_{n,J}^{[n]}(z)\Big)J
  \Big([\tilde\pi_{n,J}^{[n]}(z)]^{-1}\tilde\rho_{n,J}^{[n]}(z)\Big)^*-J \nonumber\\
& = [\tilde\pi_{n,J}^{[n]}(z)]^{-1} \Big(\tilde\rho_{n,J}^{[n]}(z)J[\tilde\rho_{n,J}^{[n]}(z)]^*
                               -\tilde\pi_{n,J}^{[n]}(z)J[\tilde\pi_{n,J}^{[n]}(z)]^*\Big)
      [\tilde\pi_{n,J}^{[n]}(z)]^{-*} \nonumber\\
& = [\tilde\pi_{n,J}^{[n]}(z)]^{-1} R_{n+1,J} [\tilde\pi_{n,J}^{[n]}(z)]^{-*} 
 \,\ge\, 0 \label{eind-C2-1}
\end{align}
for each $z\in\T\setminus\cN_{\tilde\pi_{n,J}^{[n]}}$ and, analogously,
\begin{equation}\label{eind-C2-2}
\Big(\!\tilde\tau_{n,J}^{[n]}(z)[\tilde\sigma_{n,J}^{[n]}(z)]^{-1}\!\Big)^{\!*}J
  \Big(\!\tilde\tau_{n,J}^{[n]}(z)[\tilde\sigma_{n,J}^{[n]}(z)]^{-1}\!\Big)-J 
= [\tilde\sigma_{n,J}^{[n]}(z)]^{-*} L_{n+1,J} [\tilde\sigma_{n,J}^{[n]}(z)]^{-1} 
 \,\ge\, 0 
\end{equation}
for each $z\in\T\setminus\cN_{\tilde\sigma_{n,J}^{[n]}}$.
Thus, the implications (i)$\Rightarrow$(iii) and (iv)$\Rightarrow$(i) are immediate consequences
of (\ref{eind-C2-1}), (\ref{eind-C2-2}), Lemma \ref{eind-L1}, and \cite[Theorem 1.3.3]{DFK}. Finally, the implication
(iii)$\Rightarrow$(iv) is obvious.
\qed\end{pf}

In his landmark paper \cite{Po1}, V.P. Potapov obtained  a multiplicative decomposition of a function $f\in\cP_J(\D)$ into special factors belonging to $\cP_J(\D)$. Perhaps the simplest functions belonging to the class $\cP_J(\D)$ are particular rational $m\times m$ matrix-valued functions which have exactly one pole of order $1$ in the extended complex plane. These functions will be discussed in the following two examples.\\
For $\alpha\in\D$, we denote by $b_\alpha$ the normalized elementary Blaschke factor associated with $\alpha$, i. e., for $w\in\C\backslash\{\frac{1}{\,\overline{\alpha}\,}\}$ we have
$$b_\alpha(w):=\left\{\begin{matrix}
	w & \text{, if } \alpha=0 \\
	\frac{|\alpha|}{\alpha}\frac{\alpha-w}{1-\overline{\alpha}w} & \text{, if } \alpha\neq 0.
                      \end{matrix}\right.
$$
\begin{ex}\label{eind-E1}
Let $J$ be an $m\times m$ signature matrix and let $\alpha\in\D$.
\begin{itemize}
 \item[(a)] Let $P$ be a complex $m\times m$ matrix satisfying $P\neq 0_{m\times m}$, $JP\geq 0$, and $P^2=P$. Let the function $B_{\alpha, P}:\C\backslash\{\frac{1}{\,\overline{\alpha}\,}\}\rightarrow\C^{m\times m}$ be defined by
$$B_{\alpha,P}(w):=I_m+\left[b_\alpha(w)-1\right]P.$$
For $w\in\C\backslash\{\frac{1}{\,\overline{\alpha}\,}\}$, then the identity
$$J-\left[B_{\alpha,P}(w)\right]^*J\left[B_{\alpha,P}(w)\right]=\left(1-|b_\alpha(w)|^2\right)JP$$
holds. 
Thus, the restriction of the function $B_{\alpha,P}$ onto $\D$ belongs to $\cP_{J,0}(\D)$. The function $B_{\alpha,P}$ is called \emph{the Blaschke-Potapov $J$-elementary factor of first kind associated with $\alpha$ and $P$}.
\item[(b)] Let $Q$ be a complex $m\times m$ matrix satisfying $Q\neq 0_{m\times m}$, $-JQ\geq 0$ and $Q^2=Q$. Let the function $C_{\alpha, Q}:\C\backslash\{\alpha\}\rightarrow\C^{m\times m}$ be defined by
$$C_{\alpha, Q}(w):=I_m+\left[\frac{1}{b_\alpha(w)}-1\right]Q.$$
For $w\in\C\backslash\{\alpha\}$, then the identity
$$J-\left[C_{\alpha, Q}(w)\right]^*J\left[C_{\alpha, Q}(w)\right]=\frac{1-|b_\alpha(w)|^2}{|b_\alpha(w)|^2}(-JQ)$$
holds. 
Thus, the restriction of the function $C_{\alpha, Q}$ onto $\D$ belongs to $\cP_J(\D)$. In the case $\alpha\neq 0$ this restriction even belongs to $\cP_{J,0}(\D)$. The function $C_{\alpha, Q}$ is called \emph{the Blaschke-Potapov $J$-elementary factor of second kind associated with $\alpha$ and $Q$}.
\end{itemize}
\end{ex}
Observe that if $J=-I_m$ (respectively, $J=I_m$), then there is no matrix $P\in\C^{m\times m}\backslash\{0_{m\times m}\}$ (respectively, $Q\in\C^{m\times m}\backslash\{0_{m\times m}\}$) such that $JP\geq 0$ and $P^2=P$ (respectively, $-JQ\geq 0$ and $Q^2=Q$). Consequently, if $J=-I_m$ (respectively, $J=I_m$), then there is no Blaschke-Potapov $J$-elementary factor of first (respectively, second) kind.

\begin{ex}\label{eind-E2}
Let $J$ be an $m\times m$ signature matrix and let $u\in\T$. Furthermore, let $R$ be a complex $m\times m$ matrix satisfying $R\neq 0_{m\times m}$, $JR\geq 0$, and $R^2=0_{m\times m}$. Let the function $D_{u,R}:\C\backslash\{u\}\rightarrow\C^{m\times m}$ be defined by
$$D_{u,R}(w):=I_m-\frac{u+w}{u-w}R.$$
For $w\in\C\backslash\{u\}$, then the identity
$$J-\left[D_{u,R}(w)\right]^*J\left[D_{u,R}(w)\right]=\frac{2(1-|w|^2)}{|u-w|^2}JR$$
holds. 
Thus, the restriction of the function $D_{u,R}$ onto $\D$ belongs to $\cP_{J,0}(\D)$. The function $D_{u,R}$ is called \emph{the Blaschke-Potapov $J$-elementary factor of third kind associated with $u$ and $R$}. 
\end{ex}

Observe that in the cases $J=I_m$ and $J=-I_m$ there does not exist a complex $m\times m$ matrix $R\neq 0_{m\times m}$ satisfying $JR\geq 0$ and $R^2=0_{m\times m}$, i.e., there is no 
Blaschke-Potapov $J$-elementary factor of third kind in these cases.

\vspace{2ex}
An $m\times m$ matrix function $B$ meromorphic in $\C$ is said to be a \emph{finite Blaschke-Potapov
product with respect to $J$} if 
either $B$ is a constant function with $J$-unitary value or
$B$ admits a representation $B=B_1\cdot\ldots\cdot B_n$ with some $n\in\N$
where each of the factors $B_j$ has the shape $B_j=U_j\tilde B_jV_j$ where 
$U_j$ and $V_j$ are constant $J$-unitary matrices and $\tilde B_j$ is a Blaschke-Potapov $J$-elementary factor of first or second or third kind.

\vspace{2ex}
Let us now introduce some notations. Let $\C_0:=\C\cup\{\infty\}$ 
be the extended complex plane, and $\E:=\C_0\setminus(\D\cup\T)$.
Suppose that $G$ is a simply connected domain of $\C_0$. Then let $\cNM(G)$ be the \emph{Nevanlinna class} of all functions which are meromorphic in $G$ and which can be represented as quotient of two
bounded holomorphic functions in $G$. If $f\in\cNM(\D)$ (resp., $f\in\cNM(\E)$), then a well-known
theorem due to Fatou implies that $f$ has radial boundary values $\ul\lambda$-almost everywhere on $\T$, i.e., there exists a Borel measurable function $\ul f:\T\rightarrow\C$ 
such that
$$ \lim\limits_{r\to1-0}f(rz)=\ul f(z) \quad\mbox{(resp.,}
 \lim\limits_{r\to1+0}f(rz)=\ul f(z) \mbox{)} $$
for $\ul\lambda$-a.e. $z\in\T$, where $\underline{\lambda}$ stands for the linear Lebesgue measure on $\T$.

Recall that each entry function of a given matrix-valued function $f\in\cP_J(\D)$ belongs to $\cNM(\D)$ 
(see, e.g., \cite[Corollary 3.6]{FKR1}). In particular, every $f\in\cP_J(\D)$
has radial boundary values $\ul\lambda$-a.e. on $\T$. This observation leads to the following notion. If $J$ is an $m\times m$ signature matrix and if $f\in\cP_J(\D)$, then the function $f$ is called \emph{$J$-inner}, if $f$ has $J$-unitary radial boundary values $\underline{\lambda}$-a.e. on $\T$.
In the special case $J=I_m$, the class of all $J$-inner functions in $\D$ coincides with the class
of all \emph{inner $m\times m$ Schur functions} in $\D$, i.e. the class of all functions $f\in\cS_{m\times m}(\D)$ which have unitary radial boundary values $\underline{\lambda}$-a.e. on $\T$.

\begin{rem}\label{eind-R2}
Let $J$ be an $m\times m$ signature matrix, and let $B$ be a finite Blaschke-Potapov product with respect to $J$. Then Examples \ref{eind-E1} and \ref{eind-E2} and \cite[Lemma 1.3.13]{DFK} show that $B$ is a rational $m\times m$ matrix-valued function and the restriction of $B$ onto $\D\cap\mathbb H_B$ is a $J$-inner function.
\end{rem}

For our purposes it will be important that the converse statement given below is also true. This 
well-known fact is a particular consequence of V.P. Potapov's famous factorization theorem for $J$-Potapov functions \cite{Po1}.

\begin{prop}\label{eind-P2}
Let $J$ be an $m\times m$ signature matrix, and let $B$ be a rational $m\times m$ matrix-valued function such that the restriction of $B$ onto $\D\cap\mathbb H_B$ 
is $J$-inner. Then $B$ is a finite Blaschke-Potapov product with respect to $J$.
\end{prop}

Now we can give a further characterization of the case that the interpolation problem (P) has a unique solution.

\begin{thm}\label{eind-T1}
Let $J$ be an $m\times m$ signature matrix, let $n\in\N_0$, 
and let $(A_j)_{j=0}^n$ be a $J$-Potapov sequence. Further, let $f_{c,n}$ be the $J$-central $J$-Potapov function corresponding to $(A_j)_{j=0}^n$ and let $f_{c,n}^\diamondsuit$ 
be the unique rational extension of $f_{c,n}$. 
Then the following statements are equivalent:
\begin{enumerate}\item[(i)] There is a unique $f\in\cP_{J,0}(\D)$ such that 
(\ref{ICP}) is fulfilled for each $j\in\N_{0,n}$ (namely $f=f_{c,n}$).
\item[(ii)] $f_{c,n}^\diamondsuit$ is a finite Blaschke-Potapov product with respect to $J$.
\item[(iii)] $f_{c,n}$ is $J$-inner.
\item[(iv)] There is some $z\in\T\setminus(\cN_{\rho_{n,J}}\cap\cN_{\tau_{n,J}})$ such that 
the matrix $f_{c,n}^\diamondsuit(z)$ is $J$-unitary, where $\rho_{n,J}$ and $\tau_{n,J}$
are given by (\ref{LRQ-T1-2}) and (\ref{LRQ-T2-2}) with some $V_n\in\cY_{n,J}$ 
and $W_n\in\cZ_{n,J}$ if $n\ge1$.
\end{enumerate}
\end{thm}

\begin{pf}
The equivalence of (i), (iii), and (iv) is an immediate consequence of Corollary \ref{eind-C1}.
Furthermore, Remark \ref{eind-R2} and Proposition \ref{eind-P2} yield the equivalence of (ii) and (iii), 
since $f_{c,n}$ is the restriction of a rational 
$m\times m$ matrix-valued function onto $\D\cap\mathbb H_{f_{c,n}}$.
\qed\end{pf}

We are now going to derive some kind of converse statement. More precisely, we will show that
for each finite Blaschke-Potapov product $B$ (with respect to $J$) which is holomorphic at $0$, there exist an $n\in\N_0$
and a sequence $(A_j)_{j=0}^n$ such that the restriction of $B$ onto $\D\cap\mathbb H_B$ 
is the unique $f\in\cP_{J,0}(\D)$ such that (\ref{ICP}) is fulfilled for each $j\in\N_{0,n}$.
In our proof we will refer to the corresponding result for matricial Schur functions obtained in 
\cite{FKL}. The following observation, which is based on the well-known concept of Potapov-Ginzburg transformation (Potapov \cite{Po3}, Ginzburg \cite{Gi}, see also, e.g., \cite{FKR1}), will enable us to do so.

\begin{rem}\label{eind-R3}
Let $J$ be an $m\times m$ signature matrix, and let
\begin{equation}\label{eind-R3-1}
\mbox{$\mathbf P_J:=\frac12(I+J)$ \quad  and \quad  $\mathbf Q_J:=\frac12(I-J)$.}
\end{equation}
Let $B$ be a finite Blaschke-Potapov product with respect to $J$. 
Then we see from
Remark \ref{eind-R2} and \cite[Proposition 3.8]{FKR1} that $\det(\mathbf Q_JB(w)+\mathbf P_J)\ne0$ for each $w\in\D\cap\mathbb H_B$, that
the matrix-valued function
\begin{equation}\label{eind-R3-2}
 S:=(\mathbf P_JB+\mathbf Q_J)(\mathbf Q_JB+\mathbf P_J)^{-1} 
\end{equation}
is holomorphic in $\D$, and that the restriction of $S$ onto $\D$ is an inner $m\times m$ Schur function. Since $S$ is a rational matrix-valued function, Proposition \ref{eind-P2}
yields that $S$ is a finite Blaschke-Potapov product with respect to the signature matrix $I_m$. 

Conversely, if $S$ is a finite Blaschke-Potapov product with respect to the signature matrix $I_m$ such that $\det(\mathbf Q_JS+\mathbf P_J)$ does not vanish identically, then a combination 
of Remark \ref{eind-R2} and \cite[Proposition 3.8]{FKR1} shows that the 
restriction of the matrix-valued function
\begin{equation}\label{eind-R3-3}
 B:=(\mathbf P_JS+\mathbf Q_J)(\mathbf Q_JS+\mathbf P_J)^{-1} 
\end{equation}
onto $\D\cap\mathbb H_B$ is $J$-inner. 
Applying again Proposition \ref{eind-P2}, we infer that $B$ is a finite Blasch\-ke-Potapov product with respect to $J$.
\end{rem}

\begin{prop}\label{eind-P3}
Let $J$ be an $m\times m$ signature matrix, and let $B$ be a finite Blaschke-Potapov product with respect to $J$ such that $0\in\mathbb H_B$. Furthermore, let 
\begin{equation}\label{eind-P3-1}
B(w)=\sum_{k=0}^\infty A_kw^k
\end{equation}
be the Taylor series representation of $B$ for each $w$ belonging to some neighborhood of $0$.
Then there exists an $n\in\N_0$ such that 
$\cP_{J,0}[\D,(A_j)_{j=0}^n]=\{f\}$ where $f$ denotes
the restriction of $B$ onto $\D\cap\mathbb H_B$.
\end{prop}

\begin{pf}
Let the matrix-valued function $S$ be defined by (\ref{eind-R3-2}), where $\mathbf P_J$ and $\mathbf Q_J$ are given
by (\ref{eind-R3-1}).
Remark \ref{eind-R3}
yields that $S$ is a finite Blaschke-Potapov product with respect to the signature matrix $I_m$. Let $S(w)=\sum_{k=0}^\infty C_kw^k$ be the Taylor series representation of
$S$ for each $w\in\D$. Then \cite[Proposition 7.7]{FKL} implies that there exists an 
$n\in\N_0$ such that the restriction of $S$ onto $\D$ is the unique $m\times m$ Schur function $g$
satisfying $g^{(j)}(0)=\frac{C_j}{j!}$ 
for each $j\in\N_{0,n}$. Remark \ref{eind-R2} provides us that 
$f$ belongs to $\cP_{J,0}[\D,(A_j)_{j=0}^n]$. 
Now suppose $\tilde f\in\cP_{J,0}[\D,(A_j)_{j=0}^n]$. 
From \cite[Proposition 3.4]{FKR1} it follows that $\det(\mathbf Q_J\tilde f(w)+\mathbf P_J)\ne0$ holds for each $w\in\mathbb H_{\tilde f}$ and that
the matrix-valued function
$$ \tilde g:=(\mathbf P_J\tilde f+\mathbf Q_J)(\mathbf Q_J\tilde f+\mathbf P_J)^{-1} $$
is an $m\times m$ Schur function. Let $\tilde g(w)=\sum_{k=0}^\infty \tilde C_kw^k$ be the Taylor series representation of $\tilde g$ for each $w\in\D$. 
From \cite[Remark 6.1]{FKR1}
we obtain $\tilde C_j=C_j$ for each $j\in\N_{0,n}$
and therefore $\tilde g=g$. Finally, an application of \cite[Proposition 3.4]{FKR1} yields
$\det(\mathbf Q_J\tilde g(0)+\mathbf P_J)\ne0$, $\det(\mathbf Q_J g(0)+\mathbf P_J)\ne0$, and
$$ \tilde f=(\mathbf P_J\tilde g+\mathbf Q_J)(\mathbf Q_J\tilde g+\mathbf P_J)^{-1}=(\mathbf P_Jg+\mathbf Q_J)(\mathbf Q_Jg+\mathbf P_J)^{-1}=f.$$
This completes the proof.
\qed\end{pf}

\begin{cor}\label{eind-C3}
Let $J$ be an $m\times m$ signature matrix, let $\mathbf P_J$ and $\mathbf Q_J$ be given by (\ref{eind-R3-1}),
and let $B$ be a rational $m\times m$ matrix-valued function. Then the following statements are equivalent:
\begin{enumerate}
\item[(i)] $B$ is a finite Blaschke-Potapov product with respect to $J$.
\item[(ii)] There are some $m\times m$ matrix polynomials $\pi$ and $\rho$ such that 
$\det\rho$ does not vanish identically and
the following three conditions are satisfied:
\begin{enumerate}
\item[(I)] $B=\pi\rho^{-1}$.
\item[(II)] $\det[\mathbf Q_J\pi(w)+\mathbf P_J\rho(w)]\ne0$ for each $w\in\D$.
\item[(III)] $[\rho(z)]^*J\rho(z)=[\pi(z)]^*J\pi(z)$ for each $z\in\T$.
\end{enumerate}
\item[(iii)] There are some $m\times m$ matrix polynomials $\sigma$ and $\tau$ such that $\det\tau$ does not vanish identically and the following three conditions are satisfied:
\begin{enumerate}
\item[(IV)] $B=\tau^{-1}\sigma$.
\item[(V)] $\det[\sigma(w)\mathbf Q_J-\tau(w)\mathbf P_J]\ne0$ for each $w\in\D$.
\item[(VI)] $\tau(z)J[\tau(z)]^*=\sigma(z)J[\sigma(z)]^*$ for each $z\in\T$.
\end{enumerate}
\end{enumerate}
\end{cor}

\begin{pf}
(i)$\Rightarrow$(ii): Let the matrix-valued function $S$ be defined by (\ref{eind-R3-2}). 
Then from Remark \ref{eind-R3} we know that $S$ is a finite Blaschke-Potapov product 
with respect to the signature matrix $I_m$. Thus, \cite[Corollary 7.8]{FKL} yields the existence
of some $m\times m$ matrix polynomials $\pi_S$ and $\rho_S$ such that the following three
conditions are satisfied:
\begin{enumerate}
\item[(I')] $S=\pi_S\rho_S^{-1}$.
\item[(II')] $\det\rho_S(w)\ne0$ for each $w\in\D$.
\item[(III')] $[\rho_S(z)]^*\rho_S(z)=[\pi_S(z)]^*\pi_S(z)$ for each $z\in\T$.
\end{enumerate}
Let $\pi:=\mathbf P_J\pi_S+\mathbf Q_J\rho_S$ and $\rho:=\mathbf Q_J\pi_S+\mathbf P_J\rho_S$. Then we have 
$\rho=(\mathbf Q_JS+\mathbf P_J)\rho_S$. From Remark \ref{eind-R2} and \cite[Proposition 3.4]{FKR1} it follows that $\det(\mathbf Q_JS+\mathbf P_J)$ does not vanish identically in $\D$, and therefore $\det\rho$ does not vanish identically. 
Taking into account (I'), (II'), and (III'), straightforward calculations show that the matrix polynomials $\pi$ and $\rho$ fulfill conditions (I), (II), and (III).

(ii)$\Rightarrow$(i): From (I) and (II) we see that $\det(\mathbf Q_JB+\mathbf P_J)$ does not vanish identically 
in $\D$, i.e., the matrix-valued function
$$ S:=(\mathbf P_JB+\mathbf Q_J)(\mathbf Q_JB+\mathbf P_J)^{-1} $$
is well-defined. Using (I), (II), and (III), it is readily checked
that the matrix polynomials $\pi_S:=\mathbf P_J\pi+\mathbf Q_J\rho$ and $\rho_S:=\mathbf Q_J\pi+\mathbf P_J\rho$
fulfill conditions (I'), (II'), and (III'). In particular, $S$ is a rational matrix-valued function
holomorphic in $\D$.
Moreover, from (I') and (III') we infer that $\T\setminus\cN_{\rho_S}\sq\mathbb H_S$ and that $S(z)$ is unitary for 
each $z\in\T\setminus\cN_{\rho_S}$. Since $\cN_{\rho_S}$ is a finite set, it follows that $\T\sq\mathbb H_S$ and that $S(z)$ is unitary for each $z\in\T$. Thus, a matrix version of the maximum modulus principle yields that the restriction of $S$ onto $\D$ is an inner $m\times m$ Schur function. Proposition \ref{eind-P2} ensures therefore that $S$ is a finite Blaschke-Potapov product with respect to the signature matrix $I_m$. 
Using well-known properties of the Potapov-Ginzburg transformation 
(see, e.g., \cite[formula (2.6)]{FKR1}) we obtain that $\det(\mathbf Q_JS+\mathbf P_J)$ does not vanish identically 
in $\D$ and that $B=(\mathbf P_JS+\mathbf Q_J)(\mathbf Q_JS+\mathbf P_J)^{-1}$ holds. 
Consequently, an application of Remark \ref{eind-R3} yields (i).
Thus, the equivalence of (i) and (ii) is proved. 

Taking into account additionally \cite[Remark 2.1]{FKR1},
the equivalence of (i) and (iii) can be shown similarly.
\qed\end{pf}
\section{The nondegenerate case}

The main goal of this section is to specify some of the preceding results 
for the interpolation problem (P) in the nondegenerate case.
Moreover, we will state a complete description of the Weyl matrix balls 
associated with a strict $J$-Potapov sequence $(A_j)_{j=0}^n$.
The results in this section extend the corresponding well-known results
for the nondegenerate matricial Schur problem (see, e.g., \cite{DFK}). 

Let $J$ be an $m\times m$ signature matrix and let $\kappa\in\N_0\cup\{+\infty\}$.
Throughout this section, we consider a strict $J$-Potapov sequence $(A_j)_{j=0}^\kappa$.
Observe that in this case the matrices $P_{n,J}$ and $Q_{n,J}$ are nonsingular for each $n\in\N_{0,\kappa}$.
Moreover, from \cite[Lemmata 3.3 and 3.7]{FKR2} we know that, for each $n\in\N_{0,\kappa}$,
the matrices $L_{n+1,J}$ and $R_{n+1,J}$ are nonsingular as well.
Obviously, for each $n\in\N_{1,\kappa}$, the sets $\cY_{n,J}$ and $\cZ_{n,J}$ defined by (\ref{NrYnJ}) and (\ref{NrZnJ}), respectively, fulfill
$\cY_{n,J}=\{V_n^\Box\}$ and $\cZ_{n,J}=\{W_n^\Box\}$
where
\begin{equation}\label{nd-E0}
V_n^\Box:=Q_{n-1,J}^{-1}S_{n-1}^*J_{[n-1]}y_n \quad\mbox{and}\quad
 W_n^\Box:=z_nJ_{[n-1]}S_{n-1}^*P_{n-1,J}^{-1}. 
\end{equation}
In the sequel, 
whenever a strict $J$-Potapov sequence $(A_j)_{j=0}^\kappa$ is given, then let
for every $n\in\N_{0,\kappa}$ the matrix polynomials 
$\pi_{n,J}$, $\rho_{n,J}$, $\sigma_{n,J}$, and $\tau_{n,J}$ be defined by
\begin{equation}\label{nd-E1}
\pi_{n,J}(w):=\left\{\begin{array}{cl} A_0, &\mbox{if }n=0\\
               A_0+w e_{n-1,m}(w)J_{[n-1]}P_{n-1,J}^{-1}y_n, &\mbox{if }n\in\N,
               \end{array}\right.
\end{equation}
\begin{equation}\label{nd-E2}
\rho_{n,J}(w):=\left\{\begin{array}{cl} I_m, &\mbox{if }n=0\\
                             I_m+w e_{n-1,m}(w)J_{[n-1]}S_{n-1}^*P_{n-1,J}^{-1}y_n, 
                                 &\mbox{if }n\in\N,\end{array}\right.
\end{equation}
\begin{equation}\label{nd-E3}
\sigma_{n,J}(w):=\left\{\begin{array}{cl} A_0, &\mbox{if }n=0\\
          z_nQ_{n-1,J}^{-1}J_{[n-1]}w\eps_{n-1,m}(w)+A_0,&\mbox{if }n\in\N, \end{array}\right.
\end{equation}
and
\begin{equation}\label{nd-E4}
\tau_{n,J}(w):=\left\{\begin{array}{cl} I_m, &\mbox{if }n=0\\
                     z_nQ_{n-1,J}^{-1}S_{n-1}^*J_{[n-1]}w\eps_{n-1,m}(w)+I_m,&\mbox{if }n\in\N \end{array}\right.
\end{equation}
for each $w\in\C$.
In view of \cite[Remark 5.2]{FKRS},
for each $n\in\N_{1,\kappa}$, the matrix polynomials 
$\pi_{n,J}$, $\rho_{n,J}$, $\sigma_{n,J}$, and $\tau_{n,J}$ satisfy
\begin{equation}\label{nd-E5}
\pi_{n,J}(w)=A_0+w e_{n-1,m}(w)(y_n+S_{n-1}V_n^\Box), 
\end{equation}
\begin{equation}\label{nd-E6}
\rho_{n,J}(w)=I_m+w e_{n-1,m}(w)V_n^\Box, 
\end{equation}
\begin{equation}\label{nd-E7}
\sigma_{n,J}(w)=(W_n^\Box S_{n-1}+z_n)w\eps_{n-1,m}(w)+A_0, 
\end{equation}
and
\begin{equation}\label{nd-E8}
\tau_{n,J}(w)=W_n^\Box w\eps_{n-1,m}(w)+I_m 
\end{equation}
for all $w\in\C$.
Furthermore, whenever a strict $J$-Potapov sequence $(A_j)_{j=0}^\kappa$ is given, then let
the matrix polynomials $\gC_{n,J}$ and
$\gD_{n,J}$ for every $n\in\N_{0,\kappa}$ be defined by
\begin{equation}\label{nd-E9}
\gC_{n,J}(w):=\begin{pmatrix} 
        w J\tilde\tau_{n,J}^{[n]}(w) & \ \pi_{n,J}(w) \\
        w J\tilde\sigma_{n,J}^{[n]}(w) &\  \rho_{n,J}(w) \end{pmatrix}
        \begin{pmatrix} \sqrt{L_{n+1,J}}^{-1} & 0_{m\times m} \\
                         0_{m\times m} & \sqrt{R_{n+1,J}}^{-1} \end{pmatrix}
\end{equation}
and
\begin{equation}\label{nd-E10}
\gD_{n,J}(w):=\begin{pmatrix} \sqrt{R_{n+1,J}}^{-1} & 0_{m\times m} \\
                         0_{m\times m} & \sqrt{L_{n+1,J}}^{-1} \end{pmatrix}
       \begin{pmatrix} 
        w \tilde\rho_{n,J}^{[n]}(w)J &
                 \  w \tilde\pi_{n,J}^{[n]}(w)J \\
       \sigma_{n,J}(w)  & \tau_{n,J}(w) \end{pmatrix}
\end{equation}
for each $w\in\C$.
In the case $\kappa\ge1$ we observe that,
for each $n\in\N_{1,\kappa}$, the matrix polynomials $\gC_{n,J}$ and $\gD_{n,J}$ can be factorized via
\begin{equation}\label{nd-E11}
\gC_{n,J}=\gC_{0,J}G_{1,J}\cdot\ldots\cdot G_{n,J} \qquad\mbox{and}\qquad
\gD_{n,J}=H_{n,J}\cdot\ldots\cdot H_{1,J}\gD_{0,J}
\end{equation}
where for each $k\in\N_{1,n}$, $G_{k,J}:\C\rightarrow\C^{2m\times2m}$ and 
$H_{k,J}:\C\rightarrow\C^{2m\times2m}$ are defined by
\begin{equation}\label{nd-E12}
 G_{k,J}(w)
 :=\begin{pmatrix} I_m & K_{k,J} \\ K_{k,J}^* & I_m \end{pmatrix}
  \begin{pmatrix} w\sqrt{L_{k,J}}\sqrt{L_{k+1,J}}^{-1} & 0_{m\times m} \\
          0_{m\times m} & \sqrt{R_{k,J}}\sqrt{R_{k+1,J}}^{-1}  \end{pmatrix} 
\end{equation}
and
\begin{equation}\label{nd-E13}
 H_{k,J}(w)
 :=  \begin{pmatrix} w\sqrt{R_{k+1,J}}^{-1}\sqrt{R_{k,J}} & 0_{m\times m} \\
          0_{m\times m} & \sqrt{L_{k+1,J}}^{-1}\sqrt{L_{k,J}}  \end{pmatrix} 
   \begin{pmatrix} I_m & K_{k,J}^* \\ K_{k,J} & I_m \end{pmatrix} 
\end{equation}
with 
\begin{equation}\label{nd-DefSP}
K_{k,J}:=\sqrt{L_{k,J}}^{-1}(A_k-M_{k,J})\sqrt{R_{k,J}}^{-1}
\end{equation}
(see \cite[Proposition 5.6]{FKRS}).
For a more detailed discussion of the matrix polynomials introduced in (\ref{nd-E9})--(\ref{nd-E13})
we refer the reader to \cite[Section 5]{FKRS}. The special case $J=I_m$ was already treated in 
\cite[Section 3]{FK3}.

Let us now recall the notion of linear fractional transformations of matrices
(see Potapov \cite{Po2} or \cite[Section 1.6]{DFK}). Let $p,q\in\N$, and let
$A$ and $B$ be  complex $(p+q)\times(p+q)$ matrices and let 
$$ A = \begin{pmatrix} a & \ \ b \\ c & \ \ d \end{pmatrix} 
\qquad\mbox{and}\qquad
  B = \begin{pmatrix} \alpha & \ \ \gamma \\ \beta & \ \ \delta \end{pmatrix} $$
be the block representations of $A$ and $B$ with $p\times p$ blocks $a$ and $\delta$. If the set
$$ \cQ_A^{(p,q)}
    := \{x\in\C^{p\times q} : \det(cx + d) \ne 0\} $$
is nonempty, then let $\cS_A^{(p,q)}: \cQ_A^{(p,q)}\rightarrow\C^{p\times q}$ be defined by
$$ \cS_A^{(p,q)}(x):=(ax + b)(cx + d)^{-1} .$$
If the set
$$   \cR_B^{(q,p)}
  := \{x\in\C^{p\times q} : \det(x\gamma +\delta) \ne 0\} $$
is nonempty, then let $\cT_B^{(q,p)}: \cR_B^{(q,p)} \rightarrow\C^{p\times q}$ be defined by
$$ \cT_B^{(q,p)}(x):=(x\gamma+\delta)^{-1}(x\alpha+\beta) .$$
Observe that $\cQ_A^{(p,q)}\ne\emptyset$ if and only if $\rank(c, d) = q$. 
Moreover, $\cR_B^{(q,p)}\ne\emptyset$ if and only if $\rank{\gamma\choose\delta}=p$.
%

We will now specify Theorem \ref{S1-T1} in the nondegenerate case.

\begin{thm}\label{nd-T1}
Let $J$ be an $m\times m$ signature matrix, let $n\in\N_0$, 
and let $(A_j)_{j=0}^n$ be a strict $J$-Potapov sequence. 
For every $S\in\cS_{m\times m}(\D)$ and each $w\in\D$ for which 
$S(w)\in\cQ_{\gC_{n,J}(w)}^{(m,m)}$ is satisfied,
let
\begin{equation}\label{nd-T1-1}
f_S(w):=\cS_{\gC_{n,J}(w)}^{(m,m)}(S(w)).
\end{equation}
Then, for each $S\in\cS_{m\times m}(\D)$, by (\ref{nd-T1-1}) a matrix-valued function $f_S$  meromorphic in $\D$ is given, and 
the set $\mathbb H_{f_S}$ of all $w\in\D$ at which $f_S$ is holomorphic fulfills
$$\mathbb H_{f_S}=\{w\in\D:S(w)\in\cQ_{\gC_{n,J}(w)}^{(m,m)}\}
   =\{w\in\D:S(w)\in\cR_{\gD_{n,J}(w)}^{(m,m)}\}.$$
Further, for each $S\in\cS_{m\times m}(\D)$ and each $w\in\mathbb H_{f_S}$,
$f_S$ admits the representation
$f_S(w)=\cT_{\gD_{n,J}(w)}^{(m,m)}(S(w))$.
Moreover, 
$$\cP_{J,0}[\D,(A_j)_{j=0}^n]=\{f_S:S\in\cS_{m\times m}(\D)\}$$
holds true.
\end{thm}

\begin{pf}
If $n\ge1$, then $\cY_{n,J}=\{V_n^\Box\}$, $\cZ_{n,J}=\{W_n^\Box\}$, and
(\ref{nd-E5})--(\ref{nd-E8}) hold true. Thus, due to
$\det\sqrt{L_{n+1,J}}\ne0$ and $\det\sqrt{R_{n+1,J}}\ne0$, all
the assertions follow immediately from Theorem \ref{S1-T1}.
\qed\end{pf}

In the special case $J=I_m$, Theorem \ref{nd-T1} goes back to Arov/Krein \cite{AK}, who obtained
this result as a consequence of the Adamjan/Arov/Krein paper \cite{AAK} on the matricial Nehari problem. For an alternate approach to Theorem \ref{nd-T1} in the case $J=I_m$ we refer to 
\cite[Section 3.10]{DFK} and \cite{FK2a}.

\begin{lem}\label{nd-L1A}
Let $J$ be an $m\times m$ signature matrix, let $n\in\N_0$, 
and let $(A_j)_{j=0}^n$ be a strict $J$-Potapov sequence. 
Then 
$$\det\gC_{n,J}(w)=\det\gD_{n,J}(w)=w^{(n+1)m}$$ 
holds for each $w\in\C$.
\end{lem}

\begin{pf}
From \cite[Lemma 3.11]{FKR2} we know that
\begin{equation}\label{eind-L1A-1}
\det L_{k+1,J}=\det R_{k+1,J}
\end{equation}
for each $k\in\N_{0,n}$. 
In view of 
$$ \gC_{0,J}(w)=\begin{pmatrix} 
        w J & \ A_0 \\
        w JA_0^* &\  I \end{pmatrix}
        \begin{pmatrix} \sqrt{L_{1,J}}^{-1} & 0_{m\times m} \\
                         0_{m\times m} & \sqrt{R_{1,J}}^{-1} \end{pmatrix}, $$
an application of \cite[Lemma 1.1.7]{DFK} yields therefore for each $w\in\C$
\begin{equation}\label{eind-L1A-1A}
\det\gC_{0,J}(w)
 =w^m\det(J-A_0JA_0^*)\sqrt{\det L_{1,J}}^{\,-1}\sqrt{\det R_{1,J}}^{\,-1}=w^m.
\end{equation}
Now suppose $n\ge1$.
For each $k\in\N_{1,n}$, let $K_{k,J}$ be defined by (\ref{nd-DefSP}),
and let the matrix polynomial $G_{k,J}$ for each $w\in\C$ be given by (\ref{nd-E12}).
Then \cite[Proposition 4.1]{FKR2} provides us 
\begin{equation}\label{eind-L1A-2}
L_{k+1,J}=\sqrt{L_{k,J}}(I-K_{k,J}K_{k,J}^*)\sqrt{L_{k,J}}
\end{equation}
for each $k\in\N_{1,n}$. Taking into account (\ref{nd-E12}), \cite[Lemma 1.1.7]{DFK},
(\ref{eind-L1A-1}), and (\ref{eind-L1A-2}), we obtain for each $k\in\N_{1,n}$ and each $w\in\C$
\begin{align}\label{eind-L1A-3}
\det G_{k,J}(w)
&=w^m\det(I-K_{k,J}K_{k,J}^*)\det\!\Big(\!\!\sqrt{\!L_{k,J}}\sqrt{\!L_{k+1,J}}^{-1}\Big)
                        \det\!\Big(\!\!\sqrt{\!R_{k,J}}\sqrt{\!R_{k+1,J}}^{-1}\Big) \nonumber\\
&=w^m\det(I-K_{k,J}K_{k,J}^*)\cdot\det L_{k,J}\cdot(\det L_{k+1,J})^{-1}
=w^m.
\end{align}
Similarly, it can be checked that for each $k\in\N_{1,n}$ and each $w\in\C$ the equations
\begin{equation}\label{eind-L1A-4}
\det H_{k,J}(w)=w^m \quad\mbox{and}\quad \det\gD_{0,J}(w)=w^m
\end{equation}
hold true where $H_{k,J}:\C\rightarrow\C^{2m\times 2m}$ is defined by (\ref{nd-E13}).
Hence the assertion follows from (\ref{eind-L1A-1A}), (\ref{eind-L1A-3}), (\ref{eind-L1A-4}), 
and (\ref{nd-E11}).
\qed\end{pf}

We are now interested in describing particular subsets of $\cP_{J,0}[\D,(A_j)_{j=0}^n]$
for a given strict $J$-Potapov sequence $(A_j)_{j=0}^n$.
For this purpose, let us recall the notion of pseudocontinuable functions.
Let $\cNM(\D)$ and $\cNM(\E)$ be the 
Nevanlinna classes of 
meromorphic functions introduced in Section \ref{s4}.
Then a (scalar) meromorphic function $f$ in $\D$ is said to be \emph{pseudocontinuable (into $\E$)} if it 
belongs to $\cNM(\D)$ and if there is an $f^\#\in\cNM(\E)$ such that the radial boundary values 
$\ul f$ and $\ul{f^\#}$ of $f$ and $f^\#$, respectively, coincide $\ul\lambda$-a.e. on $\T$. This function $f^\#$ is then called
the \emph{pseudocontinuation of $f$ (into $\E$)}. Clearly, a function $f\in\cNM(\D)$
admits at most one pseudocontinuation $f^\#$.
A matrix-valued meromorphic function $f$ in $\D$ is said to be \emph{pseudocontinuable} if 
each entry function of $f$ is pseudocontinuable.
It is easy to see that every rational matrix-valued function in $\D$ and every $J$-inner function
is pseudocontinuable.

Now we are going to describe the sets $\cP_{J,0}^{pc}[\D,(A_j)_{j=0}^n]$ (resp., $\cP_{J,0}^r[\D,(A_j)_{j=0}^n]$, $\ul\cP_{J,0}[\D,(A_j)_{j=0}^n]$)
of all matrix-valued functions $f\in\cP_{J,0}[\D,(A_j)_{j=0}^n]$ which are 
pseudocontinuable (resp., rational, $J$-inner). 
In the following, we denote by $\cS_{m\times m}^{pc}(\D)$ (resp., $\cS_{m\times m}^r(\D)$,
$\ul\cS_{m\times m}(\D)$) the set of all pseudocontinuable (resp., rational, inner)
$m\times m$ Schur functions in $\D$.

\begin{prop}\label{TK-P1}
Let $J$ be an $m\times m$ signature matrix, let $n\in\N_0$, 
and let $(A_j)_{j=0}^n$ be a strict $J$-Potapov sequence. 
For each $S\in\cS_{m\times m}(\D)$, let the matrix-valued function $f_S$ be defined as in Theorem \ref{nd-T1}.
Then:
\begin{enumerate}
\item[(a)] $\cP_{J,0}^{pc}[\D,(A_j)_{j=0}^n]=\{f_S:S\in\cS_{m\times m}^{pc}(\D)\}$.
\item[(b)] $\cP_{J,0}^r[\D,(A_j)_{j=0}^n]=\{f_S:S\in\cS_{m\times m}^r(\D)\}$.
\item[(c)] $\ul\cP_{J,0}[\D,(A_j)_{j=0}^n]=\{f_S:S\in\ul\cS_{m\times m}(\D)\}$.
\end{enumerate}
\end{prop}

\begin{pf}
in view of Theorem \ref{nd-T1}, Lemma \ref{nd-L1A}, and \cite[Proposition 1.6.2]{DFK}, it follows
\begin{equation}\label{TK-P1-1}
f_S(w)=\cS_{\gC_{n,J}(w)}^{(m,m)}(S(w)) \quad\mbox{and}\quad
S(w)=\cS_{\gC_{n,J}^{-1}(w)}^{(m,m)}(f_S(w))
\end{equation}
for each $S\in\cS_{m\times m}(\D)$ and each $w\in\mathbb H_{f_S}\setminus\{0\}$.
Since the matrix-valued functions $\gC_{n,J}$ and $\gC_{n,J}^{-1}$ are both rational, 
(a) and (b) follow from (\ref{TK-P1-1}), Theorem \ref{nd-T1}, and the arithmetic of pseudocontinuable 
and of rational matrix-valued functions.
Now we will prove part (c). Let 
$\gC_{n,J}=\begin{pmatrix} \gC_{n,J}^{[11]} & \gC_{n,J}^{[12]} \\ \gC_{n,J}^{[21]} & \gC_{n,J}^{[22]} 
           \end{pmatrix}$
be the block partition of $\gC_{n,J}$ into $m\times m$ blocks, and let $S\in\cS_{m\times m}(\D)$.
If $\ul S$ (resp., $\ul{f_S}$) stands for a radial boundary function of $S$ (resp., $f_S$),
then we get $\ul S(z)\in\cQ_{\gC_{n,J}(z)}^{(m,m)}$
and
$\ul{f_S}(z)=\cS_{\gC_{n,J}(z)}^{(m,m)}(\ul S(z))$
for $\ul{\lambda}$-a.e. $z\in\T$. 
Further, from \cite[Proposition 5.7]{FKRS} we know that 
$\gC_{n,J}^*(z)J^\Box\gC_{n,J}(z)=\diag(I_m,-I_m)$ holds for every $z\in\T$, where $J^\Box:=\diag(J,-J)$.
Hence a straihgtforward calculation yields
\begin{align*}
& J-\ul{f_S}^*J\ul{f_S} \\
&= (\gC_{n,J}^{[21]}\ul S+\gC_{n,J}^{[22]})^{-*} \\
& \qquad\cdot
 \big[(\gC_{n,J}^{[21]}\ul S+\gC_{n,J}^{[22]})^*J(\gC_{n,J}^{[21]}\ul S+\gC_{n,J}^{[22]})
   -(\gC_{n,J}^{[11]}\ul S+\gC_{n,J}^{[12]})^*J(\gC_{n,J}^{[11]}\ul S+\gC_{n,J}^{[12]})  \big] \\
& \qquad\cdot
  (\gC_{n,J}^{[21]}\ul S+\gC_{n,J}^{[22]})^{-1}
\displaybreak[0] \\
&= (\gC_{n,J}^{[21]}\ul S+\gC_{n,J}^{[22]})^{-*} 
  \big[-(\ul S^*,I)\cdot\gC_{n,J}^*J^\Box\gC_{n,J}\cdot(\ul S^*,I)^*\big]
    (\gC_{n,J}^{[21]}\ul S+\gC_{n,J}^{[22]})^{-1}
\displaybreak[0] \\
&  = (\gC_{n,J}^{[21]}\ul S+\gC_{n,J}^{[22]})^{-*}\big(I-\ul S^*\ul S\big)
   (\gC_{n,J}^{[21]}\ul S+\gC_{n,J}^{[22]})^{-1} 
\end{align*}
$\ul\lambda$-a.e. on $\T$. Thus, taking into account Theorem \ref{nd-T1}, we obtain (c).
\qed\end{pf}


We are now going to describe the Weyl matrix balls corresponding to a given strict $J$-Potapov
sequence $(A_j)_{j=0}^n$. For this purpose, we need some preparation.
For $p,q\in\N$ we will work with the special $(p+q)\times(p+q)$ signature matrix
$$j_{pq}:=\diag(I_p,-I_q).$$ 
Further, let $\K_{p\times q}$ be the set of all contractive 
matrices from $\C^{p\times q}$.

A closer look at the proof of \cite[Theorem 1.6.3]{DFK} shows that the following result holds,
which is slightly more general than \cite[Theorem 1.6.3]{DFK}.

\begin{lem}\label{nd-L1}
Let $p,q\in\N$, let $A$ be a nonsingular matrix from $\C^{(p+q)\times(p+q)}$, and let
\begin{equation}\label{nd-L1-0}
 W:=A^{-*}j_{pq}A^{-1} \qquad\mbox{and}\qquad V:=Aj_{pq}A^*. 
\end{equation}
Let
\begin{equation}\label{nd-L1-1}
A=\begin{pmatrix}A_{11}&\ A_{12}\\A_{21}&\ A_{22}\end{pmatrix}, \quad
W=\begin{pmatrix}W_{11}&\ W_{12}\\W_{21}&\ W_{22}\end{pmatrix}, \quad\mbox{and}\quad
V=\begin{pmatrix}V_{11}&\ V_{12}\\V_{21}&\ V_{22}\end{pmatrix} 
\end{equation}
be the block partitions of $A$, $W$, and $V$, respectively,
where $A_{11}$, $W_{11}$, and $V_{11}$ are $p\times p$ blocks.
Suppose that $A_{22}$ is nonsingular and that the matrix $A_{22}^{-1}A_{21}$ is strictly contractive.
Then the matrices $W_{11}$ and $V_{22}$ are nonsingular, the inclusion 
$\K_{p\times q}\sq\cQ_A^{(p,q)}$ holds, and the identities
$$\cS_A^{(p,q)}(\K_{p\times q})=\{Y\in\C^{p\times q}:-(Y^*,I)W(Y^*,I)^*\ge0\}
 =\gK(M;\sqrt L,\sqrt R) $$
are satisfied where
$$M:=-W_{11}^{-1}W_{12}=V_{12}V_{22}^{-1},$$
\begin{equation}\label{nd-L1-2}
L:=W_{11}^{-1}=V_{11}-V_{12}V_{22}^{-1}V_{21}>0,
\end{equation}
and
\begin{equation}\label{nd-L1-3}
R:=W_{21}W_{11}^{-1}W_{12}-W_{22}=-V_{22}^{-1}>0.
\end{equation}
\end{lem}

\begin{rem}\label{nd-RL1det}
Under the assumptions of Lemma \ref{nd-L1}, the identities
$$\frac{\det L}{\det R}=\frac{(-1)^q}{\det W}=|\det A|^2$$
hold true, where $L$ and $R$ are given by (\ref{nd-L1-2}) and (\ref{nd-L1-3}). 
This can be easily seen from the well-known Schur block decomposition of the matrix $W$
(see, e.g., \cite[Lemma 1.1.7]{DFK}).
\end{rem}

\begin{lem}\label{nd-L2}
Let $J$ be an $m\times m$ signature matrix, let $n\in\N_0$, 
and let $(A_j)_{j=0}^n$ be a strict $J$-Potapov sequence. 
Let $\cE_1:\C\rightarrow\C$ be defined by $w\mapsto w$.
In view of $\det\rho_{n,J}(0)\ne0$, let 
\begin{equation}\label{nd-L2-1}
\chi_{n,J}:=\cE_1\sqrt{R_{n+1,J}}\rho_{n,J}^{-1}J\tilde\sigma_{n,J}^{[n]}
\sqrt{L_{n+1,J}}^{-1}.
\end{equation}
Then $\{w\in\D:\det\rho_{n,J}(w)\ne0\}=\{w\in\D:\det\tau_{n,J}(w)\ne0\}$
holds, and the identity
\begin{equation}\label{nd-L2-1A}
\chi_{n,J}=\cE_1\sqrt{R_{n+1,J}}^{-1}\tilde\pi_{n,J}^{[n]}J\tau_{n,J}^{-1}\sqrt{L_{n+1,J}}
\end{equation}
is satisfied.
\end{lem}

\begin{pf}
First we note that, in the case $n\ge1$, equations (\ref{nd-E5})--(\ref{nd-E8}) hold true where
$V_n^\Box$ and $W_n^\Box$ are given by (\ref{nd-E0}).
Thus, from \cite[Propositions 2.19 and 2.20]{FKRS} we get $\{w\in\D:\det\rho_{n,J}(w)\ne0\}=\{w\in\D:\det\tau_{n,J}(w)\ne0\}$. 
Further, \cite[Remark 2.5 and Theorems 2.7. and 2.8]{FKRS} imply the identities
\begin{equation}\label{nd-L2-2}
\tau_{n,J}\pi_{n,J}=\sigma_{n,J}\rho_{n,J} \ \ \mbox{and}\ \ 
\tilde\pi_{n,J}^{[n]}\tilde\tau_{n,J}^{[n]}=\tilde\rho_{n,J}^{[n]}\tilde\sigma_{n,J}^{[n]}.
\end{equation}
Moreover, \cite[Proposition 2.11]{FKRS} yields
\begin{equation}\label{nd-L2-3}
\tilde\rho_{n,J}^{[n]}(w)J\rho_{n,J}(w)-\tilde\pi_{n,J}^{[n]}(w)J\pi_{n,J}(w)=w^nR_{n+1,J}
\end{equation}
and
\begin{equation}\label{nd-L2-4}
\tau_{n,J}(w)J\tilde\tau_{n,J}^{[n]}(w)-\sigma_{n,J}(w)J\tilde\sigma_{n,J}^{[n]}(w)=w^nL_{n+1,J}
\end{equation}
for each $w\in\C$.
Now let $w\in\D$ be such that $\det\rho_{n,J}(w)\ne0$.
Then using (\ref{nd-L2-2}), (\ref{nd-L2-3}), and (\ref{nd-L2-4}) we obtain
\begin{align*}
&w^nR_{n+1,J}[\rho_{n,J}(w)]^{-1}J\tilde\sigma_{n,J}^{[n]}(w) \\
&= \Big(\tilde\rho_{n,J}^{[n]}(w)J\rho_{n,J}(w)-\tilde\pi_{n,J}^{[n]}(w)J\pi_{n,J}(w)\Big)
[\rho_{n,J}(w)]^{-1}J\tilde\sigma_{n,J}^{[n]}(w)
\displaybreak[0] \\
&= \tilde\pi_{n,J}^{[n]}(w)J[\tau_{n,J}(w)]^{-1}
    \Big(\tau_{n,J}(w)J\tilde\tau_{n,J}^{[n]}(w)-\sigma_{n,J}(w)J\tilde\sigma_{n,J}^{[n]}(w)\Big)
\displaybreak[0] \\
&= w^n\tilde\pi_{n,J}^{[n]}(w)J[\tau_{n,J}(w)]^{-1}L_{n+1,J}
\end{align*}
and therefore (\ref{nd-L2-1A}).
\qed\end{pf}

\begin{rem}\label{nd-R3}
Let $p,q\in\N$, and let $K\in\C^{p\times q}$. Then, using \cite[Remark 1.1.2, Lemma 1.1.13]{DFK} and 
the singular value decomposition of $K$, it is readily checked that $\det(I_p+KX)\ne0$ holds for each $X\in\K_{q\times p}$ if and only if $K$ is 
strictly contractive.
\end{rem}

Now we are going to determine 
the Weyl matrix balls associated with a finite strict $J$-Potapov sequence.

\begin{thm}\label{nd-P1}
Let $J$ be an $m\times m$ signature matrix, let $n\in\N_0$, 
and let $(A_j)_{j=0}^n$ be a strict $J$-Potapov sequence. 
Further, let $\chi_{n,J}$ be given by (\ref{nd-L2-1}), and
let $w\in\D$. Then:
\begin{enumerate}
\item[(a)] The following statements are equivalent:
 \begin{enumerate}
 \item[(i)] Every $f\in\cP_{J,0}[\D,(A_j)_{j=0}^n]$ is holomorphic at $w$.
 \item[(ii)] The condition $\det\rho_{n,J}(w)\ne0$ holds, and the matrix
 $\chi_{n,J}(w)$  is strictly contractive.
 \end{enumerate}
Moreover, the set 
\begin{equation}\label{nd-P1-DefHn}
\mathbb H^{(n)}:=\bigcap\limits_{f\in\cP_{J,0}[\D,(A_j)_{j=0}^n]}\mathbb H_f
\end{equation}
is open in $\D$.
In particular, there exists a positive real number $r$ such that the inclusion
$\{v\in\D:|v|<r\}\sq\mathbb H^{(n)}$ holds.
\item[(b)] Suppose that (i) is fulfilled. 
 Then the matrices 
\begin{equation}\label{nd-P1-Phi}
\Phi_{n,J}(w):=[\tau_{n,J}(w)]^*L_{n+1,J}^{-1}\tau_{n,J}(w)
    -|w|^2J[\tilde\pi_{n,J}^{[n]}(w)]^*R_{n+1,J}^{-1}\tilde\pi_{n,J}^{[n]}(w)J
\end{equation}
and 
\begin{equation}\label{nd-P1-Psi}
\Psi_{n,J}(w):=\rho_{n,J}(w)R_{n+1,J}^{-1}[\rho_{n,J}(w)]^*
    -|w|^2J\tilde\sigma_{n,J}^{[n]}(w)L_{n+1,J}^{-1}[\tilde\sigma_{n,J}^{[n]}(w)]^*J
\end{equation}
are both positive Hermitian, and the identity
\begin{equation}\label{nd-P1-1}
 \{f(w):f\in\cP_{J,0}[\D,(A_j)_{j=0}^n]\}
 =\gK\Big(\cM_{n,J}(w);|w|^{n+1}\sqrt{\cL_{n,J}(w)},\sqrt{\cR_{n,J}(w)}\Big) 
\end{equation}
holds true where
\begin{equation}\label{nd-P1-M}
\cM_{n,J}(w)\!:=\![\Phi_{n,J}(w)]^{-1}\Big([\tau_{n,J}(w)]^*L_{n+1,J}^{-1}\sigma_{n,J}(w)
    -|w|^2J[\tilde\pi_{n,J}^{[n]}(w)]^*R_{n+1,J}^{-1}\tilde\rho_{n,J}^{[n]}(w)J\Big), 
\end{equation}   
\begin{equation}\label{nd-P1-LR}
 \cL_{n,J}(w):=[\Phi_{n,J}(w)]^{-1}, \qquad\mbox{and}\qquad
 \cR_{n,J}(w):=[\Psi_{n,J}(w)]^{-1}.
\end{equation}
Moreover, the matrix $\cM_{n,J}(w)$ admits the representation
\begin{equation}\label{nd-P1-1A}
 \cM_{n,J}(w)\!=\!\Big(\pi_{n,J}(w)R_{n+1,J}^{-1}[\rho_{n,J}(w)]^*
    -|w|^2J\tilde\tau_{n,J}^{[n]}(w)L_{n+1,J}^{-1}[\tilde\sigma_{n,J}^{[n]}(w)]^*J\Big)
     [\Psi_{n,J}(w)]^{-1}.
\end{equation}     
\item[(c)] Suppose that (i) is fulfilled. Then the identity 
$\det\cL_{n,J}(w)=\det\cR_{n,J}(w)$
holds true.
\end{enumerate}
\end{thm}

\begin{pf}
(a) First we note that in the case $\det\rho_{n,J}(w)\ne0$ we have
\begin{align}\label{nd-P1-2}
&\det
\big(wJ\tilde\sigma_{n,J}^{[n]}(w)\sqrt{L_{n+1,J}}^{-1}S(w)+\rho_{n,J}(w)\sqrt{R_{n+1,J}}^{-1}\big)
\nonumber\\
&=\det\rho_{n,J}(w)\cdot\det\sqrt{R_{n+1,J}}^{-1}\cdot\det[I_m+\chi_{n,J}(w)S(w)]
\end{align}
for each $S\in\cS_{m\times m}(\D)$.
Now suppose that (i) is satisfied. Then Theorem \ref{nd-T1} implies
\begin{equation}\label{nd-P1-3}
\det\big(wJ\tilde\sigma_{n,J}^{[n]}(w)\sqrt{L_{n+1,J}}^{-1}S(w)+\rho_{n,J}(w)\sqrt{R_{n+1,J}}^{-1} 
 \big)\ne0
\end{equation}
for each $S\in\cS_{m\times m}(\D)$. In particular, $\det\rho_{n,J}(w)\ne0$ follows.
Hence, (\ref{nd-P1-2}) and (\ref{nd-P1-3}) imply $\det[I_m+\chi_{n,J}(w)S(w)]\ne0$ 
for each $S\in\cS_{m\times m}(\D)$. Thus, Remark \ref{nd-R3} yields (ii).
Conversely, now let (ii) be satisfied. 
Then from Remark \ref{nd-R3} we get $\det[I_m+\chi_{n,J}(w)S(w)]\ne0$ for every $S\in\cS_{m\times m}(\D)$.
In view of $\det\rho_{n,J}(w)\ne0$ and (\ref{nd-P1-2}) it follows therefore that (\ref{nd-P1-3})
holds for each $S\in\cS_{m\times m}(\D)$. Taking into account Theorem \ref{nd-T1} we obtain (i).
Thus, (i) and (ii) are equivalent.
Furthermore, since $\chi_{n,J}$ is continuous in the open subset $\{v\in\D:\det\rho_{n,J}(v)\ne0\}$
of $\D$, the equivalence of (i) and (ii) implies that the set $\mathbb H^{(n)}$ is open in $\D$.
In particular, because of $0\in\mathbb H^{(n)}$,
there is an $r>0$ such that 
$\{v\in\D:|v|<r\}\sq\mathbb H^{(n)}$ holds.

(b) In view of $\{f(0):f\in\cP_{J,0}[\D,(A_j)_{j=0}^n]\}=\{A_0\}$, the case $w=0$ is trivial.
Now let $w\ne0$.
From part (a) we know that (ii) is valid. Hence we have
$$\Psi_{n,J}(w)=\rho_{n,J}(w)\sqrt{R_{n+1,J}}^{-1}\Big(I-\chi_{n,J}(w)[\chi_{n,J}(w)]^*\Big)
  \Big(\rho_{n,J}(w)\sqrt{R_{n+1,J}}^{-1}\Big)^*>0.$$
Similarly, using Lemma \ref{nd-L2}, we obtain
$\Phi_{n,J}(w)>0$. 
Taking into account (i) and Theorem \ref{nd-T1}, we get (\ref{nd-P1-3}) for each 
$S\in\cS_{m\times m}(\D)$ and, consequently, 
$\K_{m\times m}\sq\cQ_{\gC_{n,J}(w)}^{(m,m)}$.
Further, Theorem \ref{nd-T1} provides us 
\begin{equation}\label{nd-P1-4}
\{f(w):f\in\cP_{J,0}[\D,(A_j)_{j=0}^n]\}=\cS_{\gC_{n,J}(w)}^{(m,m)}(\K_{m\times m}).
\end{equation}
According to \cite[Lemma 5.4]{FKRS}, the identity
$\gD_{n,J}(w)U_{mm}\gC_{n,J}(w)=w^{n+1}U_{mm}$ holds where
\begin{equation}\label{nd-P1-Umm}
U_{mm}:=\begin{pmatrix}0_{m\times m}&\ I_m\\-I_m&\ 0_{m\times m}\end{pmatrix}.
\end{equation}
Because of $w\ne0$ this implies $\det\gC_{n,J}(w)\ne0$ and, in view of (\ref{nd-E10}),
\begin{equation}\label{nd-P1-5}
 [\gC_{n,J}(w)]^{-1}=\frac1{w^{n+1}} \begin{pmatrix}
   \sqrt{L_{n+1,J}}^{-1}\tau_{n,J}(w) & \ \ \ -\sqrt{L_{n+1,J}}^{-1}\sigma_{n,J}(w) \\
   -w\sqrt{R_{n+1,J}}^{-1}\tilde\pi_{n,J}^{[n]}(w)J 
               & \ \ \ w\sqrt{R_{n+1,J}}^{-1}\tilde\rho_{n,J}^{[n]}(w)J
   \end{pmatrix}. 
\end{equation}
Let $W:=[\gC_{n,J}(w)]^{-*}j_{mm}[\gC_{n,J}(w)]^{-1}$ and $V:=\gC_{n,J}(w)j_{mm}[\gC_{n,J}(w)]^*$.
Let the block partitions of $W$ and $V$ be given as in (\ref{nd-L1-1})
where $W_{11}$ and $V_{11}$ are $m\times m$ blocks. Having in mind condition (ii) and Lemma \ref{nd-L1}
we see that the matrices $W_{11}$ and $V_{22}$ are nonsingular. Moreover, 
using (\ref{nd-P1-5}) and (\ref{nd-E9}) it is readily checked that
\begin{equation}\label{nd-P1-5A}
\cM_{n,J}(w)=-W_{11}^{-1}W_{12}, \quad |w|^{2(n+1)}\cL_{n,J}(w)=W_{11}^{-1}, \quad\mbox{and}\quad
\cR_{n,J}(w)=-V_{22}^{-1}
\end{equation}
hold true. Thus, condition (ii) and Lemma \ref{nd-L1} provide us 
\begin{equation}\label{nd-P1-6}
\cS_{\gC_{n,J}(w)}^{(m,m)}(\K_{m\times m})=
\gK\Big(\cM_{n,J}(w);|w|^{n+1}\sqrt{\cL_{n,J}(w)},\sqrt{\cR_{n,J}(w)}\Big) 
\end{equation}
and $\cM_{n,J}(w)=V_{12}V_{22}^{-1}$. It is easy to see that the latter identity implies (\ref{nd-P1-1A}). Finally, (\ref{nd-P1-4}) and (\ref{nd-P1-6}) yield (\ref{nd-P1-1}).

(c) According to \cite[Lemma 3.11]{FKR2} the equation $\det L_{n+1,J}=\det R_{n+1,J}$ holds. 
Thus, in the case $w=0$ the assertion of (c) is obvious. Now suppose $w\ne0$.
Then using $\det\gC_{n,J}(w)\ne0$, (ii), Lemma \ref{nd-L1}, (\ref{nd-P1-5A}), 
Remark \ref{nd-RL1det}, and Lemma \ref{nd-L1A} we obtain
$$\frac{\det(|w|^{2(n+1)}\cL_{n,J}(w))}{\det\cR_{n,J}(w)}=|\det\gC_{n,J}(w)|^2=|w|^{2(n+1)m}$$
and, consequently, (c).
\qed\end{pf}

From now on, whenever some $\kappa\in\N_0\cup\{+\infty\}$ and a strict $J$-Potapov sequence
$(A_j)_{j=0}^\kappa$ are given, then, for each $n\in\N_{0,\kappa}$, let the set
$\mathbb H^{(n)}$ and 
the matrix-valued functions $\cM_{n,J}:\mathbb H^{(n)}\rightarrow\C^{m\times m}$,
$\cL_{n,J}:\mathbb H^{(n)}\rightarrow\C^{m\times m}$, and
$\cR_{n,J}:\mathbb H^{(n)}\rightarrow\C^{m\times m}$ be given by (\ref{nd-P1-DefHn}),
(\ref{nd-P1-Phi}), (\ref{nd-P1-Psi}), (\ref{nd-P1-M}), and (\ref{nd-P1-LR}), respectively.

\vspace{2ex}
A closer look at formula (\ref{nd-P1-1}) shows that there occurs an appropriate normalization of the left semi-radius of the matrix ball under consideration. This type of normalization originates in V.K. Dubovoj's paper \cite{Dub2}.

\begin{rem}\label{nd-R1}
Let $J$ be an $m\times m$ signature matrix, and let $A_0$ be a strictly $J$-contractive matrix.
Then a straightforward calculation yields 
$$I-\Big(\sqrt{R_{1,J}}JA_0^*\sqrt{L_{1,J}}^{-1}\Big)^*\Big(\sqrt{R_{1,J}}JA_0^*\sqrt{L_{1,J}}^{-1}\Big)=\sqrt{L_{1,J}}J\sqrt{L_{1,J}}.$$ 
In the case $J\ne I_m$ this implies
$\|\sqrt{R_{1,J}}JA_0^*\sqrt{L_{1,J}}^{-1}\|>1$.
Taking into account Theorem \ref{nd-P1} we see that 
$\mathbb H^{(0)}$ coincides with the set of all 
$w\in\D$ for which the matrix $w\sqrt{R_{1,J}}JA_0^*\sqrt{L_{1,J}}^{-1}$ is strictly contractive,
i.e., $\mathbb H^{(0)}=\{w\in\C:|w|<r_0\}$ where
\begin{equation}\label{nd-R1-1}
r_0:=\left\{\begin{array}{cl}1,&\mbox{ if } J=I_m\\ 
           \|\sqrt{R_{1,J}}JA_0^*\sqrt{L_{1,J}}^{-1}\|^{-1},&\mbox{ if }J\ne I_m.\end{array}\right.
\end{equation}
In particular, if $f$ is some matrix-valued function belonging to $\cP_{J,0}(\D)$ such that
$A_0:=f(0)$ is strictly $J$-contractive, then $f$ is holomorphic in the open disk
$\{w\in\C:|w|<r_0\}$.
Further, it is readily checked that in the case $J\ne I_m$ the identity 
$r_0=\|\sqrt{R_{1,J}}^{-1}A_0^*J\sqrt{L_{1,J}}\|^{-1}$ holds.
\end{rem}

\begin{rem}\label{nd-R2}
Let $J$ be an $m\times m$ signature matrix, let $n\in\N_0$, 
and let $(A_j)_{j=0}^n$ be a strict $J$-Potapov sequence. 
Then Remark \ref{nd-R1} provides us in particular 
$\{w\in\C:|w|<r_0\}\sq\mathbb H^{(n)}$ where
$r_0$ is given by (\ref{nd-R1-1}). 
\end{rem}

The following result was inspired by \cite[Proposition 7.2]{FKL2}. However, the strategy of the following proof of
Proposition \ref{nd-L6} is completely different from the proof of \cite[Proposition 7.2]{FKL2}, which is essentielly 
based on using Christoffel-Darboux formulas.

\begin{prop}\label{nd-L6}
Let $J$ be an $m\times m$ signature matrix, let $n\in\N_0$, 
and let $(A_j)_{j=0}^{n+1}$ be a strict $J$-Potapov sequence. Let $w\in\mathbb H^{(n)}$. 
Then:
\begin{enumerate}
\item[(a)] $\cL_{n,J}(w)\ge\cL_{n+1,J}(w)$ and $\cR_{n,J}(w)\ge\cR_{n+1,J}(w)$.
\item[(b)] Let $K_{n+1,J}$ be given by (\ref{nd-DefSP}), and let $\chi_{n,J}$ be defined by (\ref{nd-L2-1}). Then the following statements are equivalent:
\begin{enumerate}
\item[(i)] $\cL_{n,J}(w)=\cL_{n+1,J}(w)$.
\item[(ii)] $\cR_{n,J}(w)=\cR_{n+1,J}(w)$.
\item[(iii)] $\chi_{n,J}(w)=-K_{n+1,J}^*$.
\end{enumerate}
Moreover, if $w\ne0$, then (i) is equivalent to
\begin{enumerate}
\item[(iv)] $\cM_{n+1,J}(w)=f_{c,n+1}(w)$,
\end{enumerate}
where $f_{c,n+1}$ denotes the $J$-central $J$-Potapov function corresponding to $(A_j)_{j=0}^{n+1}$.
\end{enumerate}
\end{prop}

\begin{pf}
First we note that, in view of $\mathbb H^{(n)}\sq\mathbb H^{(n+1)}$ and part (a) of Theorem \ref{nd-P1}, the complex matrices 
$\cL_{n+1,J}(w)$, $\cR_{n+1,J}(w)$, $\cM_{n+1,J}(w)$, and $\chi_{n,J}(w)$ are well-defined.
Further, let
$G_{n+1,J}$ and $H_{n+1,J}$ be defined via (\ref{nd-E12}) and (\ref{nd-E13}), respectively.
From \cite[Proposition 4.1]{FKR2} we know that $K_{n+1,J}$ is strictly contractive and that the identities
$L_{n+2,J}=\sqrt{L_{n+1,J}}(I-K_{n+1,J}K_{n+1,J}^*)\sqrt{L_{n+1,J}}$
and
$R_{n+2,J}=\sqrt{R_{n+1,J}}(I-K_{n+1,J}^*K_{n+1,J})\sqrt{R_{n+1,J}}$
hold. Having this in mind, we obtain
\begin{align}\label{nd-L6-1}
& [H_{n+1,J}(w)]^*j_{mm}H_{n+1,J}(w) \nonumber\\
&= \begin{pmatrix} I & K_{n+1,J}^* \\ K_{n+1,J} & I \end{pmatrix} 
   \begin{pmatrix} |w|^2(I-K_{n+1,J}^*K_{n+1,J})^{-1} & 0 \\ 
                      0 & -(I-K_{n+1,J}K_{n+1,J}^*)^{-1}\end{pmatrix} \nonumber\\
&\qquad\cdot    \begin{pmatrix} I & K_{n+1,J}^* \\ K_{n+1,J} & I \end{pmatrix} . 
\end{align}
Taking into account (\ref{nd-L6-1}) and the identities 
$K_{n+1,J}(I-K_{n+1,J}^*K_{n+1,J})^{-1}=(I-K_{n+1,J}K_{n+1,J}^*)^{-1}K_{n+1,J}$
as well as 
$K_{n+1,J}^*(I-K_{n+1,J}K_{n+1,J}^*)^{-1}=(I-K_{n+1,J}^*K_{n+1,J})^{-1}K_{n+1,J}^*$,
a straightforward calculation yields
\begin{align}\label{nd-L6-2}
&j_{mm}-[H_{n+1,J}(w)]^*j_{mm}H_{n+1,J}(w) \nonumber\\
&=(1-|w|^2)
  \begin{pmatrix} (I-K_{n+1,J}^*K_{n+1,J})^{-1} & (I-K_{n+1,J}^*K_{n+1,J})^{-1}K_{n+1,J}^* \\
    K_{n+1,J}(I-K_{n+1,J}^*K_{n+1,J})^{-1} & \ \ \ K_{n+1,J}(I-K_{n+1,J}^*K_{n+1,J})^{-1}K_{n+1,J}^*
  \end{pmatrix}  \nonumber\\
&= (1-|w|^2)(I,K_{n+1,J}^*)^*(I-K_{n+1,J}^*K_{n+1,J})^{-1}(I,K_{n+1,J}^*).
\end{align}
Analoguosly it can be checked that 
\begin{equation}\label{nd-L6-3}
j_{mm}-G_{n+1,J}(w)j_{mm}[G_{n+1,J}(w)]^*
 =(1-|w|^2)(I,K_{n+1,J})^*(I-K_{n+1,J}K_{n+1,J}^*)^{-1}(I,K_{n+1,J})
\end{equation}
holds true.
Further, because of (\ref{nd-E11}) the equations
\begin{equation}\label{nd-L6-4}
\gC_{n,J}(w)G_{n+1,J}(w)=\gC_{n+1,J}(w) \quad\mbox{and}\quad H_{n+1,J}(w)\gD_{n,J}(w)=\gD_{n+1,J}(w)
\end{equation}
are fulfilled. For $k\in\{n,n+1\}$ let $\Phi_{k,J}(w)$ and $\Psi_{k,J}(w)$ be given by
(\ref{nd-P1-Phi}) and (\ref{nd-P1-Psi}), respectively. It is readily checked that 
\begin{equation}\label{nd-L6-5}
\Phi_{k,J}(w)=-(0_{m\times m},I_m)[\gD_{k,J}(w)]^*j_{mm}\gD_{k,J}(w)(0_{m\times m},I_m)^*
\end{equation}
and 
\begin{equation}\label{nd-L6-6}
\Psi_{k,J}(w)=-(0_{m\times m},I_m)\gC_{k,J}(w)j_{mm}[\gC_{k,J}(w)]^*(0_{m\times m},I_m)^*
\end{equation}
hold for each $k\in\{n,n+1\}$.
Using (\ref{nd-L6-5}), (\ref{nd-L6-4}), (\ref{nd-L6-2}), and Lemma \ref{nd-L2} we obtain
\begin{align}\label{nd-L6-7}
&\Phi_{n+1,J}(w)-\Phi_{n,J}(w) \nonumber\\
&=
(0_{m\times m},I_m)[\gD_{n,J}(w)]^*\Big(j_{mm}-[H_{n+1,J}(w)]^*j_{mm}H_{n+1,J}(w)\Big)
  \gD_{n,J}(w)(0_{m\times m},I_m)^* \nonumber
\displaybreak[0]\\
&=
(1-|w|^2)(0_{m\times m},I_m)[\gD_{n,J}(w)]^*
  (I,K_{n+1,J}^*)^*(I-K_{n+1,J}^*K_{n+1,J})^{-1}   \nonumber\\
& \quad\cdot  (I,K_{n+1,J}^*) \gD_{n,J}(w)(0_{m\times m},I_m)^* \nonumber
\displaybreak[0]\\
&=
(1-|w|^2)
 \Big(w\sqrt{R_{n+1,J}}^{-1}\tilde\pi_{n,J}^{[n]}(w)J
    +K_{n+1,J}^*\sqrt{L_{n+1,J}}^{-1}\tau_{n,J}(w)\Big)^* \nonumber\\
&\quad\cdot   \big(I-K_{n+1,J}^*K_{n+1,J}\big)^{-1}
 \Big(w\sqrt{R_{n+1,J}}^{-1}\tilde\pi_{n,J}^{[n]}(w)J
     +K_{n+1,J}^*\sqrt{L_{n+1,J}}^{-1}\tau_{n,J}(w)\Big)
  \nonumber
\displaybreak[0]\\
&=
(1-|w|^2)[\tau_{n,J}(w)]^*\sqrt{L_{n+1,J}}^{-1}\Big(\chi_{n,J}(w)+K_{n+1,J}^*\Big)^*
 \big(I-K_{n+1,J}^*K_{n+1,J}\big)^{-1}  \nonumber\\
&\quad\cdot   \Big(\chi_{n,J}(w)+K_{n+1,J}^*\Big)\sqrt{L_{n+1,J}}^{-1}\tau_{n,J}(w)
\end{align}
and, analoguosly,
\begin{align}\label{nd-L6-8}
&\Psi_{n+1,J}(w)-\Psi_{n,J}(w) \nonumber\\
&= (1-|w|^2)\rho_{n,J}(w)\sqrt{R_{n+1,J}}^{-1}\Big(\chi_{n,J}(w)+K_{n+1,J}^*\Big)
 \big(I-K_{n+1,J}K_{n+1,J}^*\big)^{-1} \nonumber\\
&\qquad\cdot
 \Big(\chi_{n,J}(w)+K_{n+1,J}^*\Big)^*\sqrt{R_{n+1,J}}^{-1}[\rho_{n,J}(w)]^*.
\end{align}
In view of Theorem \ref{nd-P1} we have $0<\cL_{k,J}(w)=[\Phi_{k,J}(w)]^{-1}$
and $0<\cR_{k,J}(w)=[\Psi_{k,J}(w)]^{-1}$ for $k\in\{n,n+1\}$.
Further, because of part (a) of Theorem \ref{nd-P1} and Lemma \ref{nd-L2}, the matrices 
$\rho_{n,J}(w)$ and $\tau_{n,J}(w)$ are both nonsingular. Consequently, since
$(I-K_{n+1,J}^*K_{n+1,J})^{-1}$ and $(I-K_{n+1,J}K_{n+1,J}^*)^{-1}$ are both positive Hermitian,
(\ref{nd-L6-7}) and (\ref{nd-L6-8}) imply (a) and the equivalence of (i), (ii), and (iii).

Now let $w\ne0$. Because of $w\in\mathbb H^{(n)}$, Theorem \ref{nd-T1} yields 
$\K_{m\times m}\sq\cQ_{\gC_{n,J}(w)}^{(m,m)}$. In particular, $K_{n+1,J}\in\cQ_{\gC_{n,J}(w)}^{(m,m)}$. Hence, in view of (\ref{nd-DefSP}),
Corollary \ref{CPF-C1} (in combination with (\ref{nd-E5})--(\ref{nd-E8})) implies
\begin{equation}\label{nd-L6-9}
f_{c,n+1}(w)=\cS_{\gC_{n,J}(w)}^{(m,m)}(K_{n+1,J}).
\end{equation}
Furthermore, from the block representation of $\gC_{n+1,J}(w)j_{mm}[\gC_{n+1,J}(w)]^*$
and Theorem \ref{nd-P1} (see formula (\ref{nd-P1-1A})) we infer that
$0_{m\times m}\in\cQ_{\gC_{n+1,J}(w)j_{mm}[\gC_{n+1,J}(w)]^*}^{(m,m)}$ and
\begin{equation}\label{nd-L6-10}
\cM_{n+1,J}(w)=\cS_{\gC_{n+1,J}(w)j_{mm}[\gC_{n+1,J}(w)]^*}^{(m,m)}(0_{m\times m})
\end{equation}
are satisfied.
Because of $w\in\mathbb H^{(n)}$, part (a) of Theorem \ref{nd-P1} provides us
$0_{m\times m}\in\cQ_{[\gC_{n,J}(w)]^*}^{(m,m)}$. Thus,
\begin{equation}\label{nd-L6-11}
\cS_{[\gC_{n,J}(w)]^*}^{(m,m)}(0_{m\times m})=[\chi_{n,J}(w)]^*.
\end{equation}
Using well-known properties of linear-fractional transformations of matrices (see, e.g., 
\cite[Proposition 1.6.3]{DFK}, we conclude from (\ref{nd-L6-10}), (\ref{nd-L6-4}),
and (\ref{nd-L6-11}) that $[\chi_{n,J}(w)]^*\in\cQ_{\gC_{n,J}(w)G_{n+1,J}(w)j_{mm}[G_{n+1,J}(w)]^*}^{(m,m)}$ 
and
\begin{align}\label{nd-L6-12}
\cM_{n+1,J}(w)
&=\cS_{\gC_{n,J}(w)G_{n+1,J}(w)j_{mm}[G_{n+1,J}(w)]^*[\gC_{n,J}(w)]^*}^{(m,m)}(0_{m\times m}) \nonumber\\
&=\cS_{\gC_{n,J}(w)G_{n+1,J}(w)j_{mm}[G_{n+1,J}(w)]^*}^{(m,m)}([\chi_{n,J}(w)]^*).
\end{align}
Applying again \cite[Proposition 1.6.3]{DFK},
we see that all terms in the equations (\ref{nd-L6-13}) below are 
well-defined and
\begin{align}\label{nd-L6-13}
&\cS_{G_{n+1,J}(w)j_{mm}[G_{n+1,J}(w)]^*}^{(m,m)}(-K_{n+1,J})
=\cS_{G_{n+1,J}(w)j_{mm}}^{(m,m)}\big(\cS_{[G_{n+1,J}(w)]^*}^{(m,m)}(-K_{n+1,J})\big) \nonumber\\
&=\cS_{G_{n+1,J}(w)j_{mm}}^{(m,m)}(0_{m\times m})
=\cS_{G_{n+1,J}(w)}^{(m,m)}(0_{m\times m})
=K_{n+1,J}.
\end{align}
In view of Lemma \ref{nd-L1A}, the matrix $\gC_{n,J}(w)$ is nonsingular. Hence
$\cS_{\gC_{n,J}(w)}^{(m,m)}$ is injective and $\big(\cS_{\gC_{n,J}(w)}^{(m,m)}\big)^{-1}\linebreak[0]=\linebreak[0]\cS_{[\gC_{n,J}(w)]^{-1}}^{(m,m)}$
(see, e.g., \cite[Proposition 1.6.2]{DFK}). 
Consequently, (\ref{nd-L6-9}), (\ref{nd-L6-12}) and \cite[Proposition 1.6.3]{DFK} imply
that (iv) holds if and only if $[\chi_{n,J}(w)]^*\in\cQ_{G_{n+1,J}(w)j_{mm}[G_{n+1,J}(w)]^*}^{(m,m)}$
and
$$\cS_{G_{n+1,J}(w)j_{mm}[G_{n+1,J}(w)]^*}^{(m,m)}([\chi_{n,J}(w)]^*)=K_{n+1,J}.$$
Due to (\ref{nd-L6-4}) and Lemma \ref{nd-L1A}, the matrix $G_{n+1,J}(w)j_{mm}[G_{n+1,J}(w)]^*$ is nonsingular. Consequently, the map $\cS_{G_{n+1,J}(w)j_{mm}[G_{n+1,J}(w)]^*}^{(m,m)}$ is injective.
Thus, taking into account additionally (\ref{nd-L6-13}) we infer that (iv) holds if and only if
(iii) is satisfied.
This completes the proof.
\qed\end{pf}

\section{Interrelations between the Weyl matrix balls connected with strict $J$-Potapov sequences
and their Potapov-Ginzburg-associated Schur sequences}

In \cite[Section 5]{FKR1} and \cite[Section 6]{FKR2} we have shown that via the concept of the $J$-PG transform there is an intimate connection between $J$-Potapov sequences and $m\times m$ Schur sequences. In the following, we will continue these studies by focussing on interrelations between the parameters of the Weyl matrix balls associated with a strict $J$-Potapov sequence and the corresponding Schur sequence.

We start with some notations. Let $\kappa\in\N_0\cup\{+\infty\}$. In the sequel, 
if $(A_j)_{j=0}^\kappa$ is a strict $J$-Potapov sequence,
then we will continue to use the notations
$\pi_{n,J}$, $\rho_{n,J}$, $\sigma_{n,J}$, $\tau_{n,J}$, $\gC_{n,J}$, and $\gD_{n,J}$ 
introduced in (\ref{nd-E1})--(\ref{nd-E4}), (\ref{nd-E9}), and (\ref{nd-E10}), respectively,
for each $n\in\N_{0,\kappa}$.
Furthermore, whenever some
strict $m\times m$ Schur sequence $(B_j)_{j=0}^\kappa$ is given, then we will assign the following matrices and
matrix polynomials to the sequence $(B_j)_{j=0}^\kappa$ by specifying 
the corresponding settings for $J$-Potapov sequences in the case $J=I_m$.
Let 
\begin{equation}\label{nd-E0-S}
\pi_0(w):=B_0, \quad \rho_0(w):=I_m, \quad \sigma_0(w):=B_0, \quad\mbox{and}\quad \tau_0(w):=I_m
\end{equation}
for each $w\in\C$. If $\kappa\ge1$, then let for every $n\in\N_{1,\kappa}$ the matrix polynomials
$\pi_n$, $\rho_n$, $\sigma_n$, and $\tau_n$ be given by
\begin{equation}\label{nd-E1-S}
\pi_n(w):=B_0+w e_{n-1,m}(w)P_{n-1}^{-1}y_n^{(B)}, 
\end{equation}
\begin{equation}\label{nd-E2-S}
\rho_n(w):=I_m+w e_{n,m}(w)(S_{n-1}^{(B)})^*P_{n-1}^{-1}y_n^{(B)}, 
\end{equation}
\begin{equation}\label{nd-E3-S}
\sigma_n(w):=z_n^{(B)}Q_{n-1}^{-1}w\eps_{n-1,m}(w)+B_0,
\end{equation}
and
\begin{equation}\label{nd-E4-S}
\tau_n(w):=z_n^{(B)}Q_{n-1}^{-1}(S_{n-1}^{(B)})^*w\eps_{n-1,m}(w)+I_m 
\end{equation}
for each $w\in\C$, where $P_{n-1}:=I-S_{n-1}^{(B)}(S_{n-1}^{(B)})^*$ and $Q_{n-1}:=I-(S_{n-1}^{(B)})^*S_{n-1}^{(B)}$.
(Here the matrices $S_{n-1}^{(B)}$, $y_n^{(B)}$, and $z_n^{(B)}$ are defined as in (\ref{NrS1}) and
(\ref{Nrynzn}), respectively, with the sequence $(B_j)_{j=0}^n$ instead of $(A_j)_{j=0}^n$.)
Furthermore, let 
\begin{equation}\label{nd-E4LR0-S}
l_1:=I_m-B_0B_0^* \quad\mbox{and}\quad r_1:=I_m-B_0^*B_0 
\end{equation}
and, in the case $\kappa\ge1$, 
\begin{equation}\label{nd-E4LR-S}
l_{n+1}:=l_1-z_n^{(B)}Q_{n-1}^{-1}(z_n^{(B)})^*
\ \ \mbox{and}\ \ 
r_{n+1}:=r_1-(y_n^{(B)})^*P_{n-1}^{-1}y_n^{(B)}
\end{equation}
for each $n\in\N_{1,\kappa}$.
Note that, for each $n\in\N_{0,\kappa}$, the matrices $l_{n+1}$ and $r_{n+1}$ are positive Hermitian if $(B_j)_{j=0}^\kappa$ is a strict $m\times m$ Schur sequence. 
Further, for each $n\in\N_{0,\kappa}$ and each
$w\in\C$, let the matrix polynomials $\cC_n$ and $\cD_n$ be defined by
\begin{equation}\label{nd-E9-S}
\cC_n(w):=\begin{pmatrix} 
        w \tilde\tau_n^{[n]}(w) & \ \pi_n(w) \\
        w \tilde\sigma_n^{[n]}(w) &\  \rho_n(w) \end{pmatrix}
        \begin{pmatrix} \sqrt{l_{n+1}}^{-1} & 0_{m\times m} \\
                         0_{m\times m} & \sqrt{r_{n+1}}^{-1} \end{pmatrix}
\end{equation}
and
\begin{equation}\label{nd-E10-S}
\cD_n(w):=\begin{pmatrix} \sqrt{r_{n+1}}^{-1} & 0_{m\times m} \\
                         0_{m\times m} & \sqrt{l_{n+1}}^{-1} \end{pmatrix}
       \begin{pmatrix} 
        w \tilde\rho_n^{[n]}(w) &
                 \  w \tilde\pi_n^{[n]}(w) \\
       \sigma_n(w)  & \tau_n(w) \end{pmatrix}.
\end{equation}

The special choice $J=I_m$ in Theorem \ref{nd-P1} leads to the following well-known
description of the Weyl matrix balls corresponding to a finite strict $m\times m$ Schur sequence
(see \cite[Theorem 3.9.2]{DFK}).

\begin{prop}\label{nd-P1A}
Let $n\!\in\!\N_0$, and let $(B_j)_{j=0}^n$ be a strict $m\times m$ Schur sequence. 
Denote by $\cS_{m\times m}[\D,\!(B_j)_{j=0}^n]$
the set of all $f\in\cS_{m\times m}(\D)$ satisfying $\frac{f^{(j)}(0)}{j!}=B_j$ for each $j\in\N_{0,n}$.
Further, let $w\in\D$. 
Then the matrices 
\begin{equation}\label{nd-P1A-Phi}
\Phi_n(w):=[\tau_n(w)]^*l_{n+1}^{-1}\tau_n(z)
    -|w|^2[\tilde\pi_n^{[n]}(w)]^*r_{n+1}^{-1}\tilde\pi_n^{[n]}(w)
\end{equation}
and 
\begin{equation}\label{nd-P1A-Psi}
\Psi_n(w):=\rho_n(w)r_{n+1}^{-1}[\rho_n(w)]^*
    -|w|^2\tilde\sigma_n^{[n]}(w)l_{n+1}^{-1}[\tilde\sigma_n^{[n]}(w)]^*
\end{equation}
are both positive Hermitian, and the identity
$$ \{f(w):f\in\cS_{m\times m}[\D,(B_j)_{j=0}^n]\}
 =\gK\big(\cM_n(w);|w|^{n+1}\sqrt{\cL_n(w)},\sqrt{\cR_n(w)}\big) $$
holds true where
\begin{equation}\label{nd-P1A-M}
\cM_n(w):=[\Phi_n(w)]^{-1}\Big([\tau_n(w)]^*l_{n+1}^{-1}\sigma_n(w)
    -|w|^2[\tilde\pi_n^{[n]}(w)]^*r_{n+1}^{-1}\tilde\rho_n^{[n]}(w)\Big), 
\end{equation}
\begin{equation}\label{nd-P1A-LR}
 \cL_n(w):=[\Phi_n(w)]^{-1}, \qquad\mbox{and}\qquad
 \cR_n(w):=[\Psi_n(w)]^{-1}.
\end{equation}
Moreover, the matrix $\cM_n(w)$ admits the representation
$$\cM_n(w)=\Big(\pi_n(w)r_{n+1}^{-1}[\rho_n(w)]^*
    -|w|^2\tilde\tau_n^{[n]}(w)l_{n+1}^{-1}[\tilde\sigma_n^{[n]}(w)]^*\Big)
     [\Psi_n(w)]^{-1}.$$
\end{prop}

From now on, whenever some $\kappa\in\N_0\cup\{+\infty\}$ and a strict $m\times m$ Schur sequence
$(B_j)_{j=0}^\kappa$ are given, then, for each $n\in\N_{0,\kappa}$, let
the matrix-valued functions $\cM_n:\D\rightarrow\C^{m\times m}$,
$\cL_n:\D\rightarrow\C^{m\times m}$, and
$\cR_n:\D\rightarrow\C^{m\times m}$ be given by 
(\ref{nd-P1A-Phi}), (\ref{nd-P1A-Psi}), (\ref{nd-P1A-M}), and (\ref{nd-P1A-LR}), respectively.

\vspace{2ex}
For the convenience of the reader, we will now recall the notion of the $J$-Potapov-Ginzburg transform of a sequence of matrices, which was introduced in \cite[Section 5]{FKR1}. 
In the sequel, whenever some $m\times m$ signature matrix $J$ is given, then let the orthoprojection matrices
$\mathbf P_J$ and $\mathbf Q_J$ be defined by (\ref{eind-R3-1}).

Now let $J$ be an $m\times m$ signature matrix, and for each $n\in\N_0$, let the $(n+1)m\times(n+1)m$ signature matrix
$J_{[n]}$ be given by (\ref{JnDef}).
Let $\kappa\in\N_0\cup\{+\infty\}$, and let $(A_j)_{j=0}^\kappa$ be a sequence of complex $m\times m$ matrices with $\det (\mathbf Q_JA_0+\mathbf P_J)\neq 0$. Then 
$\det (\mathbf Q_{J_{[n]}}S_n^{(A)}+\mathbf P_{J_{[n]}})\neq 0$ holds for each $n\in\N_{0,\kappa}$, and 
there is a unique sequence $(B_j)_{j=0}^\kappa$ of complex $m\times m$ matrices such that $\left(\mathbf P_{J_{[n]}}S_n^{(A)}+\mathbf Q_{J_{[n]}}\right)
 \left(\mathbf Q_{J_{[n]}}S_n^{(A)}+\mathbf P_{J_{[n]}}\right)^{-1}=S_n^{(B)}$ for each $n\in\N_{0,\kappa}$ (see \cite[Definition 5.9 and Proposition 5.11]{FKR1}). This sequence $(B_j)_{j=0}^\kappa$ is said to be the \emph{$J$-Potapov-Ginzburg transform} (short: \emph{$J$-PG transform}) of $(A_j)_{j=0}^\kappa$.

Now we formulate the announced connection between $J$-Potapov sequences and Schur sequences.

\begin{prop}\label{R516}
Let $J$ be an $m\times m$ signature matrix, and let $\kappa\in\N_0\cup\{+\infty\}$.
\begin{enumerate}
 \item[(a)] If $(A_j)_{j=0}^\kappa$ is a $J$-Potapov sequence (respectively, a strict $J$-Potapov sequence), then $\det (\mathbf Q_JA_0\linebreak[0]+\linebreak[0]\mathbf P_J)\neq 0$, and the $J$-PG transform $(B_j)_{j=0}^\kappa$ of $(A_j)_{j=0}^\kappa$ is an $m\times m$ Schur sequence (respectively, a strict $m\times m$ Schur sequence) which fulfills $\det (\mathbf Q_JB_0+\mathbf P_J)\neq 0$. Furthermore, $(A_j)_{j=0}^\kappa$ is the $J$-PG transform of $(B_j)_{j=0}^\kappa$.
\item[(b)] If $(B_j)_{j=0}^\kappa$ is an $m\times m$ Schur sequence (respectively, a strict $m\times m$ Schur sequence) such that $\det (\mathbf Q_JB_0+\mathbf P_J)\neq 0$, then the $J$-PG transform $(A_j)_{j=0}^\kappa$ of $(B_j)_{j=0}^\kappa$ is a $J$-Potapov sequence (respectively, a strict $J$-Potapov sequence), and $(B_j)_{j=0}^\kappa$ is the $J$-PG transform of $(A_j)_{j=0}^\kappa$.
\end{enumerate}
\end{prop}

A proof of Proposition \ref{R516} is given in \cite[Propositions 5.16 and 5.17]{FKR1}.

In order to derive the desired interrelations between the parameters of the Weyl matrix balls associated with a strict $J$-Potapov sequence and the corresponding strict Schur sequence, we need the following lemmas.

\begin{lem}\label{nd-L3}
Let $p,q\in\N$, let $A$ be a nonsingular matrix from $\C^{(p+q)\times(p+q)}$, 
Suppose that $\K_{p\times q}\sq\cQ_A^{(p,q)}$ and 
$\cS_A^{(p,q)}(\K_{p\times q})=\K_{p\times q}$ are fulfilled.
Then there exist a positive real number $\lambda$ and a $j_{pq}$-unitary matrix $U$ such that
$A=\lambda U$ is satisfied. If, additionally, $\cS_A^{(p,q)}(0_{p\times q})=0_{p\times q}$ holds, then
there are unitary matrices $U_1\in\C^{p\times p}$ and $U_2\in\C^{p\times p}$ such that
$A=\lambda\,\diag(U_1,U_2)$ is fulfilled.
\end{lem}

\begin{pf}
Let
the block partition of $A$ be given as in (\ref{nd-L1-1}) with $p\times p$ block $A_{11}$.
In view of $0_{p\times q}\in\K_{p\times q}\sq\cQ_A^{(p,q)}$ we have $\det A_{22}\ne0$. Hence 
$\det(I_q+A_{22}^{-1}A_{21}X)
\ne0$ holds for each 
$X\in\K_{p\times q}$. Thus, Remark \ref{nd-R3} yields that $A_{22}^{-1}A_{21}$ is strictly contractive. Let $V:=Aj_{pq}A^*$, and let the block partition of $V$ be given as in 
(\ref{nd-L1-1}), where $V_{11}$ is a $p\times p$ block. Then Lemma \ref{nd-L1} provides us
$\det V_{22}\ne0$ and $\K_{p\times q}=\cS_A^{(p,q)}(\K_{p\times q})=\gK(M;\sqrt L,\sqrt R)$,
where $M:=V_{12}V_{22}^{-1}$, $L:=V_{11}-V_{12}V_{22}^{-1}V_{21}>0_{p\times p}$, and
$R:=-V_{22}^{-1}>0_{q\times q}$. Thus, 
taking into account $\K_{p\times q}=\gK(0_{p\times q};I_p,I_q)$ and
using well-known properties of matrix balls (\cite[Theorems 1.1 and 1.3]{Sm}, see also \cite[Corollary 1.5.1 and Theorem 1.5.2]{DFK}), we obtain 
$M=0_{p\times q}$, $L=\rho I_p$, and $R=\frac1\rho I_q$ with some positive real number $\rho$.
This implies $V=\rho j_{pq}$.
Setting $\lambda:=\sqrt\rho$ and $U:=\frac1\lambda A$, we get therefore 
$Uj_{pq}U^*=j_{pq}$. Consequently, $U$ is $j_{pq}$-unitary. Thus, the first part of the assertion is verified. Now suppose $\cS_A^{(p,q)}(0_{p\times q})=0_{p\times q}$. Then it follows immediately that
$A_{12}=0_{p\times q}$ holds. In view of $A^*j_{pq}A=\lambda^2j_{pq}$
it is then readily checked that the second part of the assertion is true.
\qed\end{pf}

\begin{lem}\label{nd-L4}
Let $p,q\in\N$, let $A,B\in\C^{(p+q)\times(p+q)}$ be such that $|\det A|=|\det B|\ne0$.
Suppose that $\K_{p\times q}\sq\cQ_A^{(p,q)}\cap\cQ_B^{(p,q)}$ and 
$\cS_A^{(p,q)}(\K_{p\times q})=\cS_B^{(p,q)}(\K_{p\times q})$ are fulfilled.
Then the matrix $B^{-1}A$ is $j_{pq}$-unitary. If, additionally, 
$\cS_A^{(p,q)}(0_{p\times q})=\cS_B^{(p,q)}(0_{p\times q})$  holds, then
there are unitary matrices $U_1\in\C^{p\times p}$ and $U_2\in\C^{p\times p}$ such that
$B^{-1}A=\diag(U_1,U_2)$ is fulfilled.
\end{lem}

\begin{pf}
Using well-known properties
of linear fractional transformations of matrices (see, e.g., \cite[Propositions 1.6.2 and 1.6.3]{DFK}),
it follows $\K_{p\times q}\sq\cQ_{B^{-1}A}^{(p,q)}$ and
$\cS_{B^{-1}A}^{(p,q)}(\K_{p\times q})=\cS_{B^{-1}}^{(p,q)}(\cS_A^{(p,q)}(\K_{p\times q})) \linebreak[0]
 = \linebreak[0] \K_{p\times q}$.
Thus, from Lemma \ref{nd-L3} we infer that there are a positive real number $\lambda$ and a 
$j_{pq}$-unitary matrix $U$ such that $B^{-1}A=\lambda U$ is fulfilled. In view of $|\det U|=1$, this implies
$\lambda^{p+q}=\frac{|\det A|}{|\det B|\cdot|\det U|}=1$, i.e., $B^{-1}A$ is $j_{pq}$-unitary.
Now suppose 
$\cS_A^{(p,q)}(0_{p\times q})=\cS_B^{(p,q)}(0_{p\times q})$. Then
$\cS_{B^{-1}A}^{(p,q)}(0_{p\times q})=0_{p\times q}$ follows, and Lemma \ref{nd-L3} 
yields the remaining assertion.
\qed\end{pf}

\begin{lem}\label{nd-L4A}
Let $J$ be an $m\times m$ signature matrix, let $n\in\N_0$, 
and let $(A_j)_{j=0}^n$ be a strict $J$-Potapov sequence. 
Let $\fC:\D\rightarrow\C^{2m\times 2m}$ be a holomorphic matrix-valued function, and let
$\fC=\Big(\begin{array}{cc} \fC_{11} & \fC_{12} \\ \fC_{21} & \fC_{22} \end{array}\Big)$
be the block decomposition of $\fC$ into $m\times m$ blocks. Suppose that
there is some $w_0\in\mathbb H^{(n)}$ satisfying $\K_{m\times m}\sq\cQ_{\fC(w_0)}^{(m,m)}$,
that
\begin{equation}\label{nd-L4A-2}
 \cP_{J,0}[\D,(A_j)_{j=0}^n]=
  \big\{(\fC_{11}S+\fC_{12})(\fC_{21}S+\fC_{22})^{-1}:S\in\cS_{m\times m}(\D)\big\}, 
\end{equation}
and that $|\det\fC(w)|=|\det\gC_{n,J}(w)|$ holds for each $w\in\D$. Then there exists 
a $j_{mm}$-unitary matrix $U$ such that 
\begin{equation}\label{nd-L4A-2A}
 \fC=\hat\gC_{n,J}U 
\end{equation}
is satisfied, where $\hat\gC_{n,J}$ denotes the restriction of $\gC_{n,J}$ onto $\D$. If, additionally, $\fC_{12}\fC_{22}^{-1}$ coincides with the $J$-central $J$-Potapov 
function corresponding to $(A_j)_{j=0}^n$, then $U=\diag(U_1,U_2)$ holds with some 
unitary $m\times m$ matrices $U_1$ and $U_2$.
\end{lem}

\begin{pf}
Because of $\K_{m\times m}\sq\cQ_{\fC(w_0)}^{(m,m)}$, we have
$\det\fC_{22}(w_0)\ne0$ and $\det(I+\fC_{22}^{-1}(w_0)\fC_{21}(w_0)K)\ne0$ for each 
$K\in\K_{m\times m}$. Since $\mathbb H^{(n)}$ is open in $\D$ (see Theorem \ref{nd-P1}),
Remark \ref{nd-R3} and a continuity argument yield the existence of a positive real number $r$ such that $K(w_0,r):=\{w\in\C:|w-w_0|<r\}\sq\mathbb H^{(n)}$ and 
$\K_{m\times m}\sq\cQ_{\fC(w)}^{(m,m)}$ for each 
$w\in K(w_0,r)$. In view of Theorem \ref{nd-T1} and (\ref{nd-L4A-2}), this implies
$\K_{m\times m}\sq\cQ_{\gC_{n,J}(w)}^{(m,m)}\cap\cQ_{\fC(w)}^{(m,m)}$ and
$$ \cS_{\gC_{n,J}(w)}^{(m,m)}(\K_{m\times m})=\{f(w):f\in\cP_{J,0}[\D,(A_j)_{j=0}^n]\}
 =\cS_{\fC(w)}^{(m,m)}(\K_{m\times m}) $$
for each $w\in K(w_0,r)$. Furthermore, from Lemma \ref{nd-L1A} we get $|\det\fC(w)|=|\det\gC_{n,J}(w)|\ne0$ for each $w\in K(w_0,r)\setminus\{0\}$.
Let $h:=\hat\gC_{n,J}^{-1}\fC$. 
An application of Lemma \ref{nd-L4} provides us then that $h(w)$ is a $j_{mm}$-unitary matrix for each $w\in K(w_0,r)\setminus\{0\}$. Thus, using well-known properties of Potapov functions (see, e.g., \cite[Corollary 2.4.1]{DFK})
we can conclude that $h$ is a constant matrix-valued function in $\D$ with $j_{mm}$-unitary value, i.e.,
there is a $j_{mm}$-unitary matrix $U$ such that (\ref{nd-L4A-2A}) is satisfied. 
Now assume that $\fC_{12}\fC_{22}^{-1}$ coincides with the $J$-central $J$-Potapov 
function corresponding to $(A_j)_{j=0}^n$. Then Theorem \ref{LRQ-T1} (in combination with (\ref{nd-E5})--(\ref{nd-E8})) implies 
$\cS_{\fC(w)}^{(m,m)}(0_{m\times m})
 =\fC_{12}(w)\fC_{22}^{-1}(w)=\cS_{\gC_{n,J}(w)}^{(m,m)}(0_{m\times m})$
for each $w\in K(w_0,r)$. Because of (\ref{nd-L4A-2A}), 
Lemma \ref{nd-L4} yields finally the existence of unitary $m\times m$ matrices 
$U_1$ and $U_2$ such that $U=\diag(U_1,U_2)$ is satisfied.
\qed\end{pf}

%

In the sequel, if some $m\times m$ signature matrix $J$ is given, then let
$$ \cA_J:=\begin{pmatrix}\mathbf P_J&\ \mathbf Q_J\\ \mathbf Q_J&\ \mathbf P_J\end{pmatrix}
\quad\mbox{and}\quad
\cB_J:=\begin{pmatrix}-\mathbf P_J&\ \mathbf Q_J\\ \mathbf Q_J&\ -\mathbf P_J\end{pmatrix}. $$

Now let some strict $J$-Potapov sequence
$(A_j)_{j=0}^n$ be given, and let $(B_j)_{j=0}^n$ be its $J$-PG transform.
Then we will study the interrelations between the matrix polynomials
$\gC_{n,J}$ and $\gD_{n,J}$ defined by (\ref{nd-E9}) and (\ref{nd-E10}) on the one hand,
and the matrix polynomials
$\cC_n$ and $\cD_n$ defined by (\ref{nd-E9-S}) and (\ref{nd-E10-S}) on the other hand.
Here we will use the well-known notion of \emph{central $m\times m$ Schur functions}. Observe that, if an $n\in\N_0$ and
an $m\times m$ Schur sequence $(B_j)_{j=0}^n$ are given, then the central $m\times m$ Schur function corresponding to $(B_j)_{j=0}^n$ is just the $I_m$-central $I_m$-Potapov function corresponding to (the $I_m$-Potapov sequence) $(B_j)_{j=0}^n$.

\begin{prop}\label{nd-PolPGT}
Let $J$ be an $m\times m$ signature matrix, let $n\in\N_0$, 
and let $(A_j)_{j=0}^n$ be a strict $J$-Potapov sequence. 
Let $(B_j)_{j=0}^n$ be the $J$-PG transform of $(A_j)_{j=0}^n$.
Then $(B_j)_{j=0}^n$ is a strict $m\times m$ Schur sequence 
and the matrices 
$U_1:=\sqrt{L_{n+1,J}}(B_0\fQ_J-\fP_J)^*\sqrt{l_{n+1}}^{-1}$
and 
$U_2:=\sqrt{R_{n+1,J}}(\fQ_JB_0+\fP_J)\sqrt{r_{n+1}}^{-1}$
are unitary. Moreover,
\begin{equation}\label{nd-PolPGT-0} 
\cA_J\cC_n=\gC_{n,J}U \qquad\mbox{and}\qquad \cD_n\cB_J=V^*\gD_{n,J}
\end{equation}
holds, where $U:=\diag(-U_1,U_2)$ and $V:=\diag(-U_2,U_1)$.
\end{prop}

\begin{pf}
In view of Proposition \ref{R516}, $(B_j)_{j=0}^n$ is a strict $m\times m$ Schur sequence.
Denote by $\hat\cC_n$ the restriction of $\cC_n$ onto $\D$.
Let
$ \cA_J\hat\cC_n=\begin{pmatrix}\hat\cC_{n,J}^{[11]} &\ \hat\cC_{n,J}^{[12]} \\ 
                    \hat\cC_{n,J}^{[21]} &\ \hat\cC_{n,J}^{[22]}   \end{pmatrix}   $
be the block partition of $\cA_J\hat\cC_n$ into $m\times m$ blocks. 
Due to \cite[Theorem 7.4]{FKR1}, we have 
$\K_{m\times m}\sq\cQ_{\cA_J\cC_n(0)}^{(m,m)}$ and
$$ \cP_{J,0}[\D,(A_j)_{j=0}^n]
 =\big\{(\hat\cC_{n,J}^{[11]} S+\hat\cC_{n,J}^{[12]})
   (\hat\cC_{n,J}^{[21]} S+\hat\cC_{n,J}^{[22]})^{-1} : S\in\cS_{m\times m}(\D)\big\}. $$
An application of Lemma \ref{nd-L1A} yields
$\det(\cC_n(w))=w^{(n+1)m}=\det(\gC_{n,J}(w))$ for each $w\in\C$. Because of $\cA_J^2=I$
we get $|\det\cA_J|=1$ and therefore $|\det(\cA_J\hat\cC_n(w))|=|\det(\hat\gC_{n,J}(w))|$ 
for each $w\in\D$. 
Furthermore, let $f_{c,n}$ be the $J$-central $J$-Potapov function corresponding to $(A_j)_{j=0}^n$,
and let $g_{c,n}$ be the central $m\times m$ Schur function corresponding to $(B_j)_{j=0}^n$.
Then a combination of \cite[Proposition 6.9, Remark 6.7]{FKR2} and \cite[Theorem 1.2]{FK3}
provides us $f_{c,n}(w)\in\cQ_{\cA_J}^{(m,m)}$, $0_{m\times m}\in\cQ_{\cC_n(w)}^{(m,m)}$, 
and $\cS_{\cA_J}^{(m,m)}(f_{c,n}(w))=g_{c,n}(w)=\cS_{\cC_n(w)}^{(m,m)}(0_{m\times m})$
for each $w\in\mathbb H_{f_{c,n}}$.
Thus, using \cite[Propositions 1.6.2 and 1.6.3]{DFK} we obtain
$0_{m\times m}\in\cQ_{\cA_J\cC(w)}^{(m,m)}$ and 
$f_{c,n}(w)=\cS_{\cA_J\cC(w)}^{(m,m)}(0_{m\times m})$
for each $w\in\mathbb H_{f_{c,n}}$, i.e., $f_{c,n}=\hat\cC_{n,J}^{[12]}(\hat\cC_{n,J}^{[22]})^{-1}$
holds. Hence an application of Lemma \ref{nd-L4A} provides us the existence of some unitary
$m\times m$ matrices $U_1^\Box$ and $U_2^\Box$ such that 
\begin{equation}\label{nd-PolPGT-1}
\cA_J\cC_n=\gC_{n,J}\diag(U_1^\Box,U_2^\Box).
\end{equation}
Comparing the right lower $m\times m$ blocks and the left upper $m\times m$ blocks 
in (\ref{nd-PolPGT-1}) we get
\begin{equation}\label{nd-PolPGT-3}
(\fQ_J\pi_n+\fP_J\rho_n)\sqrt{r_{n+1}}^{-1}=\rho_{n,J}\sqrt{R_{n+1,J}}^{-1}U_2^\Box. 
\end{equation}
and
$(\fP_J\tilde\tau_n^{[n]}+\fQ_J\tilde\sigma_n^{[n]})\sqrt{l_{n+1}}^{-1}
 =J\tilde\tau_{n,J}^{[n]}\sqrt{L_{n+1,J}}^{-1}U_1^\Box$.
The latter equation implies
\begin{equation}\label{nd-PolPGT-4}
\sqrt{l_{n+1}}^{-1}(\tau_n\fP_J-\sigma_n\fQ_J)=(U_1^\Box)^*\sqrt{L_{n+1,J}}^{-1}\tau_{n,J}.
\end{equation}
In view of $\tau_{n,J}(0)=\rho_{n,J}(0)=\tau_n(0)=\rho_n(0)=I_m$ and
$\sigma_n(0)=\pi_n(0)=B_0$, we obtain therefore from (\ref{nd-PolPGT-4}) and (\ref{nd-PolPGT-3}) 
that the identities
$U_1^\Box=-U_1$ and $U_2^\Box=U_2$ hold.
Consequently, $U_1$ and $U_2$ are unitary, and (\ref{nd-PolPGT-1}) implies the first equation
stated in (\ref{nd-PolPGT-0}).

Further, taking into account \cite[Lemma 5.4]{FKRS}, 
we get
\begin{equation}\label{nd-L5-8}
\gD_{n,J}(w)U_{mm}\gC_{n,J}(w)=w^{n+1}U_{mm} \quad\mbox{and}\quad
 \cD_n(w)U_{mm}\cC_n(w)=w^{n+1}U_{mm}
\end{equation}
for each $w\in\C$, where $U_{mm}$ is given by (\ref{nd-P1-Umm}).
Obviously $V^*=U_{mm}U^{-1}U_{mm}$ holds. 
Further, a short calculation yields $\cB_JU_{mm}\cA_J=-U_{mm}$ and $U_{mm}^2=-I_{2m}$.
Having this in mind, it can be readily checked that (\ref{nd-L5-8}) and the first equation in (\ref{nd-PolPGT-0}) imply the second equation stated in (\ref{nd-PolPGT-0}). We omit the details.
Thus, the proof is complete.
\qed\end{pf}

Now we are able to derive the desired interrelations between the 
parameters of the Weyl matrix balls associated with a strict $J$-Potapov sequence $(A_j)_{j=0}^n$ 
and its $J$-PG transform.

\begin{prop}\label{nd-P1B}
Let $J$ be an $m\times m$ signature matrix, let $n\in\N_0$, and let $(A_j)_{j=0}^n$ be a strict $J$-Potapov sequence. Furthermore, let $(B_j)_{j=0}^n$ be the $J$-PG transform
of $(A_j)_{j=0}^n$. Then $(B_j)_{j=0}^n$ is a strict $m\times m$ Schur sequence and, 
for each $w\in\mathbb H^{(n)}$, the identities
\begin{equation}\label{nd-P1B-3}
\cL_{n,J}(w)
\!=\!\!\Big(\!(\cM_{n}(w)\mathbf Q_J-\mathbf P_{\!\!J})^*[\cL_{n}(w)]^{-1}
  (\cM_{n}(w)\mathbf Q_J-\mathbf P_{\!\!J})
          -|w|^{2(n+1)}\mathbf Q_J\cR_{n}(w)\mathbf Q_J\!\Big)^{-1}, 
\end{equation}
\begin{equation}\label{nd-P1B-4}
\cR_{n,J}(w)
\!=\!\Big(\!(\mathbf Q_J\cM_{n}(w)+\mathbf P_{\!\!J})[\cR_{n}(w)]^{-1}
   (\mathbf Q_J\cM_{n}(w)+\mathbf P_{\!\!J})^*
      -|w|^{2(n+1)}\mathbf P_{\!\!J}\cL_{n}(w)\mathbf P_{\!\!J}\!\Big)^{-1},
\end{equation}
\begin{equation}\label{nd-P1B-1}
\cM_{n,J}(w)
\!=\!\Big(\!(\mathbf P_{\!\!J}\cM_{n}(w)+\mathbf Q_J)[\cR_{n}(w)]^{-1}
  (\mathbf Q_J\cM_{n}(w)+\mathbf P_{\!\!J})^*
      -|w|^{2(n+1)}\mathbf P_{\!\!J}\cL_{n}(w)\mathbf Q_J\!\Big) \cR_{n,J}(w),
\end{equation}
and
\begin{equation}\label{nd-P1B-2}
\cM_{n,J}(w)
\!=\! 
  \cL_{n,J}(w) \Big(\!(\cM_{n}(w)\mathbf Q_J-\mathbf P_{\!\!J})^*[\cL_{n}(w)]^{-1}
   (\mathbf Q_J-\cM_{n}(w)\mathbf P_{\!\!J})
       +|w|^{2(n+1)}\mathbf Q_J\cR_{n}(w)\mathbf P_{\!\!J}\!\Big)
\end{equation}
hold true.
\end{prop}

\begin{pf}
Proposition \ref{R516} shows that $(B_j)_{j=0}^n$ is a strict $m\times m$ Schur sequence.
In the case $w=0$ formulas (\ref{nd-P1B-3}) and (\ref{nd-P1B-4}) are an immediate
consequence of \cite[Proposition 6.4]{FKR2}, wheras 
(\ref{nd-P1B-1}) and (\ref{nd-P1B-2}) follow from 
(\ref{nd-P1B-3}), (\ref{nd-P1B-4}),
and \cite[formula (2.6) and Remark 2.1]{FKR1}.
Now suppose $w\in\mathbb H^{(n)}\setminus\{0\}$.
Because of Lemma \ref{nd-L1A} we have $\det\gC_{n,J}(w)\ne0$.
Let $W:=[\gC_{n,J}(w)]^{-*}j_{mm}[\gC_{n,J}(w)]^{-1}$
and $V:=\gC_{n,J}(w)j_{mm}[\gC_{n,J}(w)]^*$, and let the block partitions of $W$ and $V$ be given as in (\ref{nd-L1-1}) with $m\times m$ blocks $W_{11}$ and $V_{11}$. In the same way as in the proof of Theorem \ref{nd-P1} we see then that $W_{11}$ and $V_{22}$ are nonsingular and that 
(\ref{nd-P1-5A}) is fulfilled. 
Furthermore, from part (a) of Theorem \ref{nd-P1} and Lemma \ref{nd-L1} (see formula (\ref{nd-L1-3}))
we get $V_{22}^{-1}=W_{22}-W_{21}W_{11}^{-1}W_{12}$. 
Hence, using (\ref{nd-P1-5A}) and \cite[Lemma 1.1.7]{DFK}, we obtain
\begin{align}\label{nd-P1B-5}
&[\gC_{n,J}(w)]^{-*}j_{mm}[\gC_{n,J}(w)]^{-1} \nonumber\\
& =\begin{pmatrix} I & \ \ -\cM_{n,J}(w) \\ 0 & I \end{pmatrix}^*
    \begin{pmatrix}|w|^{-2(n+1)}[\cL_{n,J}(w)]^{-1} & 0 \\ 0 & -\cR_{n,J}(w) \end{pmatrix}
    \begin{pmatrix} I & \ \ -\cM_{n,J}(w) \\ 0 & I \end{pmatrix}. 
\end{align}
Since $(B_j)_{j=0}^n$ is a strict $I_m$-Potapov sequence, (\ref{nd-P1B-5}) implies in particular
\begin{align}\label{nd-P1B-6}
&[\cC_n(w)]^{-*}j_{mm}[\cC_n(w)]^{-1} \nonumber\\
& =\begin{pmatrix} I & \ \ -\cM_n(w) \\ 0 & I \end{pmatrix}^*
    \begin{pmatrix}|w|^{-2(n+1)}[\cL_n(w)]^{-1} & 0 \\ 0 & -\cR_n(w) \end{pmatrix}
    \begin{pmatrix} I & \ \ -\cM_n(w) \\ 0 & I \end{pmatrix}.
\end{align}
Furthermore, from Proposition \ref{nd-PolPGT} it follows that there is a $j_{mm}$-unitary matrix $U$
such that $\cA_J\cC_n(w)=\gC_{n,J}(w)U$ is satisfied. In particular,
$U$ is invertible and $U^{-1}$ is $j_{mm}$-unitary as well. Consequently, in view of $\cA_J^2=I$,
\begin{align}\label{nd-P1B-7}
[\gC_{n,J}(w)]^{-*}j_{mm}[\gC_{n,J}(w)]^{-1}
&=[\gC_{n,J}(w)]^{-*}U^{-*}j_{mm}U^{-1}[\gC_{n,J}(w)]^{-1} \nonumber\\
&=\cA_J[\cC_{n}(w)]^{-*}j_{mm}[\cC_{n}(w)]^{-1}\cA_J
\end{align}
holds true. Combining (\ref{nd-P1B-5}), (\ref{nd-P1B-6}), and (\ref{nd-P1B-7}), we get
\begin{align}\label{nd-P1B-8}
& \begin{pmatrix} I & 0 \\ -(\cM_{n,J}(w))^* & \ \ I \end{pmatrix}
    \begin{pmatrix}[\cL_{n,J}(w)]^{-1} & 0 \\ 0 & -|w|^{2(n+1)}\cR_{n,J}(w) \end{pmatrix}
    \begin{pmatrix} I & \ \ -\cM_{n,J}(w) \\ 0 & I \end{pmatrix} \nonumber\\
&= \begin{pmatrix} (\mathbf P_J-\cM_{n}(w)\mathbf Q_J)^* & \ \ \ \mathbf Q_J \\
               (\mathbf Q_J-\cM_{n}(w)\mathbf P_J)^*  & \ \ \ \mathbf P_J \end{pmatrix}
    \begin{pmatrix}[\cL_{n}(w)]^{-1} & 0 \\ 0 & -|w|^{2(n+1)}\cR_{n}(w) \end{pmatrix} \nonumber\\
& \quad\quad\quad\cdot
   \begin{pmatrix} \mathbf P_J-\cM_{n}(w)\mathbf Q_J & \ \ \ \mathbf Q_J-\cM_{n}(w)\mathbf P_J \\ 
                     \mathbf Q_J & \mathbf P_J \end{pmatrix}.
\end{align}
Comparing the left upper $m\times m$ blocks and the right upper $m\times m$ blocks in 
(\ref{nd-P1B-8}), we obtain (\ref{nd-P1B-3}) and (\ref{nd-P1B-2}). Equations (\ref{nd-P1B-4}) and (\ref{nd-P1B-1})
can be proven similarly.
\qed\end{pf}

Note that Proposition \ref{nd-L6} includes in particular the case of strict $m\times m$
Schur sequences. This fact suggests the following observation.

\begin{rem}\label{ndPGT-R1}
Let $J$ be an $m\times m$ signature matrix, let $n\in\N$, let $(A_j)_{j=0}^n$ be a strict $J$-Potapov sequence, and let $w\in\mathbb H^{(n)}$. Furthermore, let $(B_j)_{j=0}^n$ be the $J$-PG transform
of $(A_j)_{j=0}^n$. Denote by $f_{c,n}$ the $J$-central $J$-Potapov function corresponding
to $(A_j)_{j=0}^n$, and let $g_{c,n}$ be the central $m\times m$ Schur function corresponding
to the strict $m\times m$ Schur sequence $(B_j)_{j=0}^n$. 
From \cite[Proposition 6.9]{FKR2} and \cite[Proposition 3.4]{FKR1} we know that $\det(\fQ_Jg_{c,n}(w)+\fP_J)\ne0$
and $f_{c,n}(w)=(\fP_Jg_{c,n}(w)+\fQ_J)(\fQ_Jg_{c,n}(w)+\fP_J)^{-1}$ hold.
Hence, if $\cL_{n-1,J}(w)=\cL_{n,J}(w)$ and $\cL_{n-1}(w)=\cL_n(w)$ are satisfied, then
Proposition \ref{nd-L6} and the (trivial) identities $\cM_{n,J}(0)=f_{c,n}(0)$ and
$\cM_n(0)=g_{c,n}(0)$ imply
$$\cM_{n,J}(w)=(\fP_J\cM_n(w)+\fQ_J)(\fQ_J\cM_n(w)+\fP_J)^{-1}.$$
\end{rem}

%
%

\section{Limit behaviour of the Weyl matrix balls associated with a nondegenerate $J$-Potapov function}

In this section we are going to study the limit behaviour of the parameters of the Weyl matrix balls
associated with a given nondegenerate $J$-Potapov function. 

In the following, if $f$ is some $m\times m$ matrix-valued function which is holomorphic at $0$,
then the sequence $(A_j)_{j=0}^\infty$ in the Taylor series representation (\ref{NrTS})
of $f$ around the origin will be shortly called the Taylor coefficient sequence of $f$.

Let $J$ be an $m\times m$ signature matrix. 
A $J$-Potapov function $f$ in $\D$ will be called \emph{nondegenerate} if it belongs to $\cP_{J,0}(\D)$ and
if the Taylor coefficient sequence of $f$ is a strict $J$-Potapov sequence. 
In the sequel, 
we will write $\cP_{J,0,\infty}(\D)$ for the class of all nondegenerate $J$-Potapov functions in $\D$.
Furthermore, whenever some $f\in\cP_{J,0,\infty}(\D)$ is given, we will use the following notations.
For every $n\in\N_0$, let  
\begin{equation}\label{ndLim-DefHfn}
\mathbb H_f^{(n)}:=\bigcap\limits_{\varphi\in\cP_{J,0}[\D,(A_j)_{j=0}^n]}\mathbb H_\varphi,
\end{equation}
where $(A_j)_{j=0}^\infty$ is the Taylor coefficient sequence of $f$,
and let the matrix-valued functions 
$\cM_{n,J}:\mathbb H_f^{(n)}\rightarrow\C^{m\times m}$,
$\cL_{n,J}:\mathbb H_f^{(n)}\rightarrow\C^{m\times m}$, and
$\cR_{n,J}:\mathbb H_f^{(n)}\rightarrow\C^{m\times m}$ be given by 
(\ref{nd-P1-Phi}), (\ref{nd-P1-Psi}), (\ref{nd-P1-M}), and (\ref{nd-P1-LR}), respectively.

\begin{lem}\label{ndLim-L1}
Let $J$ be an $m\times m$ signature matrix, and let $f\in\cP_{J,0,\infty}(\D)$. Then, for each $n\in\N_0$, the set $\mathbb H_f^{(n)}$ is open in $\D$ with $0\in\mathbb H_f^{(n)}$ and 
$\mathbb H_f^{(n)}\sq\mathbb H_f^{(n+1)}$. Furthermore, 
$\mathbb H_f=\bigcup\limits_{n=0}^\infty\mathbb H_f^{(n)}$ holds true.
\end{lem}

\begin{pf}
Taking into account (\ref{ndLim-DefHfn}) and part (a) of Theorem \ref{nd-P1} we see that, 
for each $n\in\N_0$, the set 
$\mathbb H_f^{(n)}$ is an open subset of $\D$ with
$0\in\mathbb H_f^{(n)}$ and that $\mathbb H_f^{(n)}\sq\mathbb H_f^{(n+1)}$ holds for every $n\in\N_0$.
Obviously, $\bigcup\limits_{n=0}^\infty\mathbb H_f^{(n)}\sq\mathbb H_f$ is valid.
Now let $w\in\mathbb H_f$, and assume that $w\not\in\bigcup\limits_{n=0}^\infty\mathbb H_f^{(n)}$.
Let $(A_j)_{j=0}^\infty$ be the Taylor coefficient sequence of $f$. 
Then, for each $n\in\N_0$, there is a $\varphi_n\in\cP_{J,0}[\D,(A_j)_{j=0}^n]$ such that
$w\not\in\mathbb H_{\varphi_n}$, which contradicts Proposition \ref{CPF-P2}. Thus, the assertion follows.
\qed\end{pf}

\begin{prop}\label{nd-P2}
Let $J$ be an $m\times m$ signature matrix, and let $f\in\cP_{J,0,\infty}(\D)$. 
Further, let $w\in\mathbb H_f$. 
In view of Lemma \ref{ndLim-L1}, let $n_0\in\N_0$ be such that $w\in\mathbb H_f^{(n)}$
for each $n\in\N_{n_0,\infty}$.
Then:
\begin{enumerate}
\item[(a)] The sequences $(\cL_{n,J}(w))_{n=n_0}^\infty$ and 
$(\cR_{n,J}(w))_{n=n_0}^\infty$ are both monotonously nonincreasing.
In particular, the limits
\begin{equation}\label{nd-LRinf}
\cL_J(w):=\lim\limits_{n\to\infty}\cL_{n,J}(w) \quad\mbox{and}\quad
\cR_J(w):=\lim\limits_{n\to\infty}\cR_{n,J}(w)
\end{equation}
exist and are both nonnegative Hermitian. Moreover, $\det\cL_J(w)=\det\cR_J(w)$ is satisfied.
\item[(b)] The sequence $(\cM_{n,J}(w))_{n=n_0}^\infty$ converges and
$\lim\limits_{n\to\infty}\cM_{n,J}(w)=f(w)$ holds.
\end{enumerate}
\end{prop}

\begin{pf}
(a) The assertion  of (a) is an immediate consequence of Theorem \ref{nd-P1} and Proposition \ref{nd-L6}.

(b) Part (a) implies in particular
\begin{equation}\label{ndlim-P3-1}
\lim\limits_{n\to\infty}|w|^{n+1}\sqrt{\cL_{n,J}(w)}=0_{m\times m}.
\end{equation}
For each $n\in\N_{n_0,\infty}$, let 
$\gK_n(w):=\gK\big(\cM_{n,J}(w);|w|^{n+1}\sqrt{\cL_{n,J}(w)},\sqrt{\cR_{n,J}(w)}\big)$.
Then Theorem \ref{nd-P1} provides us $\gK_{n+1}(w)\sq\gK_n(w)$ for each $n\in\N_{n_0,\infty}$.
Hence, using a well-known result from the theory of matrix balls (\cite{Sm}, see also, e.g., 
\cite[Theorem 1.5.3]{DFK}), we obtain that the sequence $(\cM_{n,J}(w))_{n=n_0}^\infty$ converges
to some complex $m\times m$ matrix $\cM_J(w)$ and, in view of (a) and (\ref{ndlim-P3-1}), that
$\bigcap\limits_{n=n_0}^\infty\gK_n(w)
 =\gK\big(\cM_J(w);0_{m\times m},\sqrt{\cR_J(w)}\big)=\{\cM_J(w)\}$.
Finally, since Theorem \ref{nd-P1} implies $f(w)\in\bigcap\limits_{n=n_0}^\infty\gK_n(w)$,
we get $f(w)=\cM_J(w)$.
\qed\end{pf}

Denote by $\cS_{m\times m,\infty}(\D)$ the class of all 
nondegenerate $m\times m$ Schur functions in $\D$, i.e., the class of all
$g\in\cS_{m\times m}(\D)$ whose Taylor coefficient sequence 
is a strict $m\times m$ Schur sequence. Obviously, $\cS_{m\times m,\infty}(\D)=\cP_{I_m,0,\infty}(\D)$
holds true. In the sequel, whenever some $g\in\cS_{m\times m,\infty}(\D)$ is given, then let,
for every $n\in\N_0$, the matrix-valued functions 
$\cM_n:\D\rightarrow\C^{m\times m}$,
$\cL_n:\D\rightarrow\C^{m\times m}$, and
$\cR_n:\D\rightarrow\C^{m\times m}$ be given by 
(\ref{nd-P1A-Phi}), (\ref{nd-P1A-Psi}), (\ref{nd-P1A-M}), and (\ref{nd-P1A-LR}), respectively, where $(B_j)_{j=0}^\infty$ is the Taylor coefficient sequence of $g$.

Choosing $J=I_m$, we reobtain from the above proposition the corresponding well-known statement
for $m\times m$ Schur functions in $\D$ (see \cite[Theorems 3.11.2 and 5.6.1]{DFK}), i.e., if $g\in\cS_{m\times m,\infty}(\D)$, then the sequences $(\cM_n(w))_{n=0}^\infty$, 
$(\cL_n(w))_{n=0}^\infty$, and $(\cR_n(w))_{n=0}^\infty$
converge for each $w\in\D$. Moreover, for each $w\in\D$
\begin{equation}\label{ndlim-SchM}
\lim\limits_{n\to\infty}\cM_n(w)=g(w)
\end{equation}
holds and the limits
\begin{equation}\label{nd-LRinfB}
\cL(w):=\lim\limits_{n\to\infty}\cL_n(w) \quad\mbox{and}\quad
\cR(w):=\lim\limits_{n\to\infty}\cR_n(w)
\end{equation}
are nonnegative Hermitian.

Let us now consider some $m\times m$ signature matrix $J$ and
a function $f\in\cP_{J,0,\infty}(\D)$. According to \cite[Proposition 6.3]{FKR1}, the function
\begin{equation}\label{ndLim-PGT}
g:=(\mathbf P_Jf+\mathbf Q_J)(\mathbf Q_Jf+\mathbf P_J)^{-1}
\end{equation}
is well-defined and belongs to $\cS_{m\times m,\infty}(\D)$.
We will now describe the interrelations between
the limit semi-radius functions $\cL_J$ and $\cR_J$ associated with $f$, which are given by
(\ref{nd-LRinf}), and the limit semi-radius functions $\cL$ and $\cR$ associated with $g$, which are defined in (\ref{nd-LRinfB}).

\begin{prop}\label{ndlim-P4}
Let $J$ be an $m\times m$ signature matrix, let $f\in\cP_{J,0,\infty}(\D)$, and let 
the functions $\cL_J:\mathbb H_f\rightarrow\C^{m\times m}$ and 
$\cR_J:\mathbb H_f\rightarrow\C^{m\times m}$ be given by (\ref{nd-LRinf}).
Further, let $g$ be defined by (\ref{ndLim-PGT}), and let 
the functions $\cL:\D\rightarrow\C^{m\times m}$ and 
$\cR:\D\rightarrow\C^{m\times m}$ be given by (\ref{nd-LRinfB}).
Then, for each $w\in\mathbb H_f$, the matrices $g(w)\mathbf Q_J-\mathbf P_J$ and
$\mathbf Q_Jg(w)+\mathbf P_J$ are nonsingular, and
\begin{equation}\label{ndLim-P4-1}
\cL_J(w)
 =\big(g(w)\mathbf Q_J-\mathbf P_J\big)^{-1}\cL(w)\big(g(w)\mathbf Q_J-\mathbf P_J\big)^{-*}
\end{equation}
and
\begin{equation}\label{ndLim-P4-2}
\cR_J(w)
 =\big(\mathbf Q_Jg(w)+\mathbf P_J\big)^{-*}\cR(w)\big(\mathbf Q_Jg(w)+\mathbf P_J\big)^{-1}
\end{equation}
hold true.
\end{prop}

\begin{pf}
Let $w\in\mathbb H_f$.
From \cite[Proposition 3.4 and Remark 2.1]{FKR1} it follows that the 
matrices $g(w)\mathbf Q_J-\mathbf P_J$ and $\mathbf Q_Jg(w)+\mathbf P_J$ are nonsingular.
Let $(A_j)_{j=0}^\infty$ (respectively, $(B_j)_{j=0}^\infty$) be the Taylor coefficient sequence
of $f$ (respectively, $g$). 
Then from \cite[Remark 6.1 and Proposition 5.11]{FKR1} we infer that, for each $n\in\N_0$,  $(B_j)_{j=0}^n$ is the $J$-PG transform of $(A_j)_{j=0}^n$. 
In view of Lemma \ref{ndLim-L1}, let $n_0\in\N$ be such that $w\in\mathbb H_f^{(n)}$ for each $n\in\N_{n_0,\infty}$. Now let $n\in\N_{n_0,\infty}$. 
Then Proposition \ref{nd-P1B} shows that (\ref{nd-P1B-3}) and (\ref{nd-P1B-4}) hold. 
Because of Theorem \ref{nd-P1} and Proposition \ref{nd-P1A}
the matrix $H_{n,J}(w):=[\cL_{n,J}(w)]^{-1}+|w|^{2(n+1)}\mathbf Q_J\cR_n(w)\mathbf Q_J$
is positive Hermitian. Hence (\ref{nd-P1B-3}) implies
\begin{align}\label{ndLim-P4-3}
&\cL_{n,J}(w)\Big(I+|w|^{2(n+1)}\mathbf Q_J\cR_n(w)\mathbf Q_J\cL_{n,J}(w)\Big)^{-1} \nonumber\\
&=[H_{n,J}(w)]^{-1} 
= \big((\cM_n(w)\mathbf Q_J-\mathbf P_J\big)^{-1}\cL_n(w)\big(\cM_n(w)\mathbf Q_J-\mathbf P_J\big)^{-*}.
\end{align}
Having in mind (\ref{ndlim-SchM}), formula (\ref{ndLim-P4-1}) follows from (\ref{ndLim-P4-3})
by taking limits for $n\to\infty$. Using (\ref{nd-P1B-4}), equation (\ref{ndLim-P4-2}) can be shown analogously.
\qed\end{pf}

In view of Proposition \ref{ndlim-P4}, the following statement follows easily from the corresponding 
well-known result for $m\times m$ Schur functions in $\D$.


\begin{prop}\label{ndlim-P5}
Let $J$ be an $m\times m$ signature matrix, let $f\in\cP_{J,0,\infty}(\D)$, and let 
the functions $\cL_J:\mathbb H_f\rightarrow\C^{m\times m}$ and 
$\cR_J:\mathbb H_f\rightarrow\C^{m\times m}$ be given by (\ref{nd-LRinf}).
Then the functions $\rank\cL_J$ and $\rank\cR_J$ are constant in $\mathbb H_f$.
\end{prop}

\begin{pf}
According to \cite[Proposition 6.3]{FKR1}, the function $g$ given by (\ref{ndLim-PGT})
is well-defined and belongs to $\cS_{m\times m,\infty}(\D)$. Let 
the functions $\cL:\D\rightarrow\C^{m\times m}$ and 
$\cR:\D\rightarrow\C^{m\times m}$ be given by (\ref{nd-LRinfB}).
Then \cite[Theorem 3.11.2]{DFK} shows that the functions 
$\rank\cL$ and $\rank\cR$ are constant in $\D$.
Thus, an application of Proposition \ref{ndlim-P4} yields the assertion.
\qed\end{pf}

The proof of Proposition \ref{ndlim-P5} is mainly based on \cite[Theorem 3.11.2]{DFK}. A closer analysis of the proof of \cite[Theorem 3.11.2]{DFK}
shows that the assertion is a consequence of a famous result due to Orlov \cite{Or} on the limit behaviour of the parameters of a sequence 
of nested matrix balls (see also \cite[Section 2.5]{DFK}).

\begin{rem}\label{ndlim-R1}
Let $J$ be an $m\times m$ signature matrix, let $f\in\cP_{J,0,\infty}(\D)$, and let 
the functions $\cL_J:\mathbb H_f\rightarrow\C^{m\times m}$ and 
$\cR_J:\mathbb H_f\rightarrow\C^{m\times m}$ be given by (\ref{nd-LRinf}).
Then Proposition \ref{nd-P2} shows in particular that $\rank\cL_J(0)=m$ if and only if 
$\rank\cR_J(0)=m$. The trivial example of a constant function $f$ defined on $\D$ with some
strictly $J$-contractive value $A_0$ yields $\cL_J(0)=J-A_0JA_0^*$ and $\cR_J(0)=J-A_0^*JA_0$.
Hence there exists an $f\in\cP_{J,0,\infty}(\D)$ with $\rank\cL_J(0)=\rank\cR_J(0)=m$.
\end{rem} 

\begin{rem}\label{ndlim-R2}
Let $J$ be an $m\times m$ signature matrix, let $n\in\N_0$, and let $(A_j)_{j=0}^n$ be a strict $J$-Potapov sequence.
Let the matrix polynomial $\gC_{n,J}$ be given by (\ref{nd-E9}), and let $\chi_{n,J}$ be defined by (\ref{nd-L2-1}).
Let $w\in\mathbb H^{(n)}$. In view of part (a) of Theorem \ref{nd-P1}, $\chi_{n,J}(w)$ is a (well-defined) strictly 
contractive matrix. For each $v\in\D$ with $-[\chi_{n,J}(w)]^*\in\cQ_{\gC_{n,J}(v)}^{(m,m)}$ let
\begin{equation}\label{ndlim-R2-1}
f_w(v):=\cS_{\gC_{n,J}(v)}^{(m,m)}\big(-[\chi_{n,J}(w)]^*\big)
\end{equation}
Then Theorem \ref{nd-T1} implies that formula (\ref{ndlim-R2-1}) defines a matrix-valued function $f_w$ meromorphic
in $\D$ with $\mathbb H_{f_w}=\{v\in\D:-[\chi_{n,J}(w)]^*\in\cQ_{\gC_{n,J}(v)}^{(m,m)}\}$ and that 
$f_w\in\cP_{J,0}[\D,(A_j)_{j=0}^n]$ holds. 
Let $A_{n+1}:=M_{n+1,J}-\sqrt{L_{n+1,J}}[\chi_{n,J}(w)]^*\sqrt{R_{n+1,J}}$. Since $-[\chi_{n,J}(w)]^*$ 
is strictly contractive, \cite[Theorem 3.9]{FKR2} yields that $(A_j)_{j=0}^{n+1}$ is a strict $J$-Potapov sequence.
Moreover, from Corollary \ref{CPF-C1} and (\ref{nd-E5})--(\ref{nd-E8}) we can conclude that $f_w$ is the $J$-central
$J$-Potapov function corresponding to $(A_j)_{j=0}^{n+1}$. In particular, \cite[Proposition 5.3]{FKR2} implies 
$f_w\in\cP_{J,0,\infty}(\D)$. 
\end{rem}

The following two results were inspired by \cite[Theorem 7.3, Corollary 7.4]{FKL2}.

\begin{prop}\label{ndlim-P6}
Let $J$ be an $m\times m$ signature matrix, let $n\in\N_0$, and let $(A_j)_{j=0}^n$ be a strict $J$-Potapov sequence.
Let $f\in\cP_{J,0}[\D,(A_j)_{j=0}^n]\cap\cP_{J,0,\infty}(\D)$, and let $w\in\mathbb H_f^{(n)}$.
\begin{enumerate}
\item[(a)] The following statements are equivalent:
\begin{enumerate}
\item[(i)] $f$ coincides with the matrix-valued function $f_w$ defined in Remark \ref{ndlim-R2}.
\item[(ii)] $\cL_{k,J}(w)=\cL_{k+1,J}(w)$ for each $k\in\N_{n,\infty}$.
\item[(iii)] $\cR_{k,J}(w)=\cR_{k+1,J}(w)$ for each $k\in\N_{n,\infty}$.
\item[(iv)] $\cL_{n,J}(w)=\cL_J(w)$.
\item[(v)] $\cR_{n,J}(w)=\cR_J(w)$.
\end{enumerate}
\item[(b)] Let $w\ne0$. Further, denote by $(A_j)_{j=0}^\infty$ the Taylor coefficient sequence of $f$. 
Then (i) holds if and only if $\cM_{k+1,J}(w)=f_{c,k+1}(w)$ for each $k\in\N_{n,\infty}$,
where $f_{c,k+1}$ stands for the $J$-central $J$-Potapov function corresponding to $(A_j)_{j=0}^{k+1}$.
\end{enumerate}
\end{prop}

\begin{pf}
For each $k\in\N_0$, let the matrix polynomials $\pi_{k,J}$, $\rho_{k,J}$, $\sigma_{k,J}$, $\tau_{k,J}$, $\gC_{k,J}$
and the matrix-valued function $\chi_{k,J}$ be defined by (\ref{nd-E1})--(\ref{nd-E4}), (\ref{nd-E9}), and
(\ref{nd-L2-1}), with respect to the strict $J$-Potapov sequence $(A_j)_{j=0}^k$. 
Further, for each $k\in\N$, let $K_{k,J}$ be given by (\ref{nd-DefSP}), and let the 
matrix polynomial $G_{k,J}$ be defined by (\ref{nd-E12}). 
Part (a) of Theorem \ref{nd-P1} yields $\det\rho_{k,J}(w)\ne0$ for each $k\in\N_{n,\infty}$. Thus, we obtain
$w\in\mathbb H_{\chi_{k,J}}$, $0_{m\times m}\in\cR_{\gC_{k,J}(w)}^{(m,m)}$, 
and $\chi_{k,J}(w)=\cT_{\gC_{k,J}(w)}^{(m,m)}(0_{m\times m})$ for each $k\in\N_{n,\infty}$. Hence, using (\ref{nd-E11}),
(\ref{nd-E12}), and well-known properties of linear-fractional transformations of matrices
(see, e.g., \cite[Proposition 1.6.3]{DFK}), we get $\chi_{k,J}(w)\in\cR_{G_{k+1,J}(w)}^{(m,m)}$ and
\begin{align}\label{ndlim-P6-1}
& \chi_{k+1,J}(w)
= \cT_{\gC_{k+1,J}(w)}^{(m,m)}(0_{m\times m}) = \cT_{\gC_{k,J}(w)G_{k+1,J}(w)}^{(m,m)}(0_{m\times m})   \nonumber\\
&=\cT_{G_{k+1,J}(w)}^{(m,m)}\big(\cT_{\gC_{k,J}(w)}^{(m,m)}(0_{m\times m})\big) 
=\cT_{G_{k+1,J}(w)}^{(m,m)}(\chi_{k,J}(w)) \nonumber\\
&= w\sqrt{R_{k+2,J}}\sqrt{R_{k+1,J}}^{-1}\Big(\chi_{k,J}(w)K_{k+1,J}+I\Big)^{-1}\Big(\chi_{k,J}(w)+K_{k+1,J}^*\Big)\!
  \sqrt{L_{k+1,J}}\sqrt{L_{k+2,J}}^{-1} 
\end{align}
for each $k\in\N_{n,\infty}$.

(i)$\Rightarrow$(ii): In view of (i), Remark \ref{ndlim-R2} yields 
\begin{equation}\label{ndlim-P6-2}
A_{n+1}=M_{n+1,J}-\sqrt{L_{n+1,J}}[\chi_{n,J}(w)]^*\sqrt{R_{n+1,J}}
\end{equation}
and therefore $\chi_{n,J}(w)=-K_{n+1,J}^*$. Moreover, from Remark \ref{ndlim-R2} we know that $f$ is the $J$-central $J$-Potapov function corresponding to $(A_j)_{j=0}^{n+1}$. Thus, \cite[Proposition 4.1]{FKR2}
implies $K_{k+1,J}=0_{m\times m}$ for each $k\in\N_{n+1,\infty}$. Taking into account (\ref{ndlim-P6-1}),
we get inductively $\chi_{k+1,J}(w)=0_{m\times m}$ for each $k\in\N_{n,\infty}$. Thus, we have shown that 
$\chi_{k,J}(w)=-K_{k+1,J}^*$ for each $k\in\N_{n,\infty}$. Therefore, Proposition \ref{nd-L6} yields (ii).

(ii)$\Rightarrow$(i): In view of (ii), Proposition \ref{nd-L6} provides us $\chi_{k,J}(w)=-K_{k+1,J}^*$ for each $k\in\N_{n,\infty}$. Thus, from (\ref{ndlim-P6-1}) we obtain $\chi_{k+1,J}(w)=0_{m\times m}$ and hence
$K_{k+2,J}=0_{m\times m}$ for each $k\in\N_{n,\infty}$. Consequently, \cite[Proposition 4.1]{FKR2}
implies that $f$ is the $J$-central $J$-Potapov function corresponding to $(A_j)_{j=0}^{n+1}$. Further, because
of $\chi_{n,J}(w)=-K_{n+1,J}^*$ we see that (\ref{ndlim-P6-2}) holds. Therefore, Remark \ref{ndlim-R2} yields (i).

Furthermore, the equivalences (ii)$\Leftrightarrow$(iii), (ii)$\Leftrightarrow$(iv), and
(iii)$\Leftrightarrow$(v) follow immediately from Proposition \ref{nd-L6}. Thus, part (a) is shown.
Finally, part (b) is an immediate consequence of part (a) and Proposition \ref{nd-L6}. 
\qed\end{pf}

\begin{cor}\label{ndlim-C1}
Let $J$ be an $m\times m$ signature matrix, let $n\in\N_0$, and let $(A_j)_{j=0}^n$ be a strict $J$-Potapov sequence.
Further, let $w\in\mathbb H^{(n)}$, and let $f_w$ be defined as in Remark \ref{ndlim-R2}. Then
$$\cM_{k,J}(w)=f_w(w)$$
holds for each $k\in\N_{n,\infty}$, where the matrix-valued function 
$\cM_{k,J}:\mathbb H_{f_w}^{(k)}\rightarrow\C^{m\times m}$ is defined with respect to $f_w$.
\end{cor}

\begin{pf}
The case $w=0$ is obvious. 
Now let $w\ne0$. Because of Remark \ref{ndlim-R2} we have $f_w\in\cP_{J,0}[\D,(A_j)_{j=0}^n]$.
Thus, Theorem \ref{nd-P1} provides us (\ref{nd-P1-1A}). Having in mind 
(\ref{ndlim-R2-1}), (\ref{nd-E9}), (\ref{nd-L2-1}), and (\ref{nd-P1-1A}),
a straightforward calculation yields
\begin{equation}\label{ndlim-C1-1}
f_w(w)=\cS_{\gC_{n,J}(w)}^{(m,m)}\big(-[\chi_{n,J}(w)]^*\big)=\cM_{n,J}(w).
\end{equation}
Thus, the assertion follows from (\ref{ndlim-C1-1}), Remark \ref{ndlim-R2}, and Proposition \ref{ndlim-P6}.
\qed\end{pf}

\end{document}